\documentclass[11pt,a4paper]{article}

\usepackage{amsmath,amssymb}
\newcommand{\e}{\varepsilon}
\newcommand{\va}{\varphi}
\newcommand{\D}{\Delta}

\newcommand{\n}{\nabla}
\newcommand{\N}{\frac{N}{2}}

\newcommand{\NN}{\frac{N}{p}}
\newcommand{\q}{q_{l}}
\newcommand{\ui}{u_{l}}
\newcommand{\g}{\int_{\mathbb{R}^{N}}}
\newcommand{\p}{\partial}

\newcommand{\R}{\mathbb{R}}
\newcommand{\h}{\hookrightarrow}

\newcommand{\de}{\delta}

\newtheorem{definition}{Definition}
\newtheorem{theorem}{Theorem}
\newtheorem{proposition}{Proposition}
\newtheorem{corollaire}{Corollary}

\newtheorem{notation}{Notation}
\newtheorem{remarka}{Remark}

\newtheorem{lem}{Lemma}
\newtheorem{lemme}{Lemma}

\title{Existence of global strong solution for Korteweg system with large infinite energy initial data}
\author{Boris Haspot  \thanks{Ceremade, Umr Cnrs 7534, Universit\'e Paris Dauphine, place du Mar\' echal De Lattre De Tassigny 75775 Paris cedex 16 (France), haspot@ceremade.dauphine.fr}}
\date{}
\begin{document}
%\tableofcontents
\maketitle
%Dans le second theoreme il faut vraiment que $q$ avec index $1$ pour boucler les estimations sur la pression.\\
%Comment passer de l'equation $\ln\rho$ a l'equation en $\rho$. Dans notre cas on a une explosion forte en norme $L^{\infty}$ disons polynomiale a cause de $BMO$.\\
%Dans mon resultat sur Navier-Stokes il semble qu'en cuttant une pression a l'infini, telle que $P(\rho)$ soit dans $L^{\infty}$ fait marcher les choses, a expliquer du moins pour avoir des gains d'integrabilite (pas tout a fait car quand $\rho$ grand cela peut pŽter, il faudrait la cuter sur des $\rho$ grand, disons $P(\rho)$ (que se passe t-il si $P(\rho)$ est $L^{\infty}$, a t-on un gain d'integrabilite, a voir...)
%In passing our result allow to consider initial data with discontinuous interfaces (up our knowledge in the literature the results of existence of strong solutions  require generally continuous initial density), it implies also that we can work with initial data space of infinite energy.  
\begin{abstract}
This work is devoted to the
study of the initial boundary value problem for a general
isothermal model of capillary fluids derived by J.E Dunn and
J.Serrin (1985)  (see \cite{fDS}), which can be used as a phase transition model. We will prove the existence of local and global (under a condition of smallness on the initial data) strong solutions  with discontinuous initial density when $\ln\rho_{0}$ belongs in the Besov space $B^{\N}_{2,\infty}(\R^{N})$. Our result relies on the fact that the density  can be written as the sum of the solution $\rho_{L}$  associated to linear system and a remainder density $\bar{\rho}$ which is more regular than $\rho_{L}$ by taking into account the regularizing effects induced on the bilinear convection term. The main difficulty concerns the proof of new estimate of maximum principle type for the linear system associated to the Korteweg system, the proof is based on a characterization of the Besov space in terms of the semi group associated to this linear system. Let also point out that we prove the existence of global strong solution with a smallness hypothesis which is subcritical in terms of the scaling of the equations, it allows us to exhibit a family of large energy initial data for the scaling of the equations providing global strong solution. In particular for the first time up our knowledge we show the existence of global strong solution for some large energy initial data when $N=2$.\\
We finish this paper by introducing the notion of quasi-solutions for the Korteweg's system  (a tool which has been developed in the framework of the compressible Navier-Stokes equations \cite{arxiv,arxiv1,hal,cras1,cras2}) which enables us to improve the previous result and to obtain the existence of global strong solution with large initial velocity in $B^{\N-1}_{2,\infty}$. As a corollary, we get global existence (and uniqueness) for highly compressibleKorteweg system when $N\geq2$. It means that for any large initial data (under an irrotational condition on the initial velocity) we have the existence of global strong solution provided that the pressure is sufficiently highly compressible.
%We obtain also  a result of ill-posedness for Korteweg system and we derive a new blow-up criterion which is the second main result of this paper. More precisely we show that if we control only the vacuum (i.e $\frac{1}{\rho}\in L^{\infty}_{T}(\dot{B}^{0}_{N+\e,1}(\R^{N}))$ with $\e>0$ ) then we can extend the strong solutions in finite time.% It extends previous results obtained for compressible equations.
\end{abstract}
\section{Introduction}
We are concerned with compressible fluids endowed with internal
capillarity. The model we consider  originates from the XIXth
century work by Van der Waals and Korteweg \cite{VW,fK} and was
actually derived in its modern form in the 1980s using the second
gradient theory, see for instance \cite{fDS,fJL,fTN}. The first investigations begin with the Young-Laplace theory which claims that the phases are separated by a hypersurface and that the jump in the pressure across the hypersurface is proportional to the curvature of the hypersurface. The main difficulty consists in describing the location and the movement of the interfaces.\\
Another major problem is to understand whether the interface behaves as a discontinuity in the state space (sharp interface) or whether the phase boundary corresponds to a more regular transition (diffuse interface, DI).
The diffuse interface models have the advantage to consider only one set of equations in a single spatial domain (the density takes into account the different phases) which considerably simplifies the mathematical and numerical study (indeed in the case of sharp interfaces, we have to treat a problem with free boundary).\\
Another approach corresponds to determine equilibrium solutions which classically consists in the minimization of the free energy functional.
Unfortunately this minimization problem has an infinity of solutions, and many of them are physically wrong (some details are given later). In order to overcome this difficulty, Van der Waals in the XIX-th century was the first to add a term of capillarity to select the physically correct solutions, modulo the introduction of a diffuse interface. This theory is widely accepted as a thermodynamically consistent model for equilibria.\\
Korteweg-type models are based on an extended version of
nonequilibrium thermodynamics, which assumes that the energy of the
fluid not only depends on standard variables but also on the
gradient of
the density. Alternatively, another way to penalize the high density variations consists in applying a zero order but non-local operator to the density gradient ( see \cite{9Ro}, \cite{5Ro}, \cite{Rohdehdr}). For more results on non local Korteweg system, we refer also to \cite{Cha, CH,CH2,Has1,Has4}.\\
Let us now consider a fluid of density $\rho\geq 0$, velocity field $u\in\R^{N}$, we are now interested in the following
compressible capillary fluid model, which can be derived from a Cahn-Hilliard like free energy (see the
pioneering work by J.- E. Dunn and J. Serrin in \cite{fDS} and also in
\cite{fA,fC,fGP}).
The conservation of mass and of momentum write:
\begin{equation}
\begin{cases}
\begin{aligned}
&\frac{\p}{\p t}\rho+{\rm div}(\rho u)=0,\\
&\frac{\p}{\p t}(\rho u)+{\rm div}(\rho
u\otimes u)-\rm div(2\mu(\rho) D (u))-\n\big(\lambda(\rho)){\rm div}u\big)+\n P(\rho)={\rm div}K,
\end{aligned}
\end{cases}
\label{3systeme}
\end{equation}
where the Korteweg tensor read as following:
\begin{equation}
{\rm div}K
=\n\big(\rho\kappa(\rho)\D\rho+\frac{1}{2}(\kappa(\rho)+\rho\kappa^{'}(\rho))|\n\rho|^{2}\big)
-{\rm div}\big(\kappa(\rho)\n\rho\otimes\n\rho\big).
\label{divK}
\end{equation}
$\kappa$ is the coefficient of capillarity and is a regular function. The term
${\rm div}K$  allows to describe the variation of density at the interfaces between two phases, generally a mixture liquid-vapor. $P$ is a general increasing pressure term.
%In the sequel to simplify, $\kappa(\rho)$ will take the form $\kappa(\rho)=\rho^{\alpha}$
%with $\alpha\in\R$.\\
$D (u)=(\n u+^{t}\n u)$ defines the stress tensor, $\mu$ and $\lambda$ are the two Lam\'e viscosity coefficients depending on the density $\rho$ and satisfying:
$$\mu>0\;\;\mbox{and}\;\;2\mu+N\lambda\geq0.$$
In what follows, we are interested in investigating the existence of global strong solution for the system (\ref{3systeme}) when
we authorize  jump in the pressure across the interfaces. In order to obtain a such result, we are going to show the  existence of strong solutions in critical space (it means in spaces as large as possible) for the scaling of the equations with an initial density which is not necessary continuous (let us point out that it will be one of the main interest of this paper).\\
Before detailing this program, let us recall briefly the classical energy estimates for the system (\ref{3systeme}). Let $\bar{\rho}>0$ be a constant reference density (in what follows, we shall assume that $\bar{\rho}=1$),   and let $\Pi $ be defined
by:
$$\Pi(s)=s\biggl(\int^{s}_{\bar{\rho}}\frac{P(z)}{z^{2}}dz-\frac{P(\bar{\rho})}{\bar{\rho}}\biggl),$$
so that $P(s)=s\Pi^{'}(s)-\Pi(s)\, ,\,\Pi^{'}(\bar{\rho})=0$. Multiplying the equation of momentum conservation in the system
(\ref{3systeme}) by $u$ and integrating by parts over $(0,t)\times\R^{N}$,
we obtain the following
estimate:
\begin{equation}
\begin{aligned}
&\int_{\R^{N}}\big(\frac{1}{2}\rho
|u|^{2}+(\Pi(\rho)-\Pi(\bar{\rho}))+\frac{1}{2}\kappa(\rho)|\nabla\rho|^{2}\big)(t)dx
+2\int_{0}^{t}\int_{\R^{N}}\mu(\rho)
|D(u)|^{2}dxdt\\
&\hspace{0,5cm}+\int_{0}^{t}\int_{\R^{N}}\lambda(\rho)
({\rm div}u)^{2}dxdt\leq\int_{\R^{N}}\big(\rho_{0}|u_{0}|^{2}+(\Pi(\rho_{0})-\Pi(\bar{\rho}))
+\frac{1}{2}\kappa(\rho_{0})|\nabla\rho_{0}|^{2}\big)dx.
\label{3inegaliteenergie1}
\end{aligned}
\end{equation}
It follows that assuming that the initial total energy is finite:
$${\cal E}_{0}=\int_{\R^{N}}\big(\rho_{0}|u_{0}|^{2}+(\Pi(\rho_{0})-\Pi(\bar{\rho}))
+\frac{\kappa(\rho_{0})}{2}|\nabla\rho_{0}|^{2}\big)dx<+\infty\,,$$ then we
have the a priori following bounds when $P(\rho)=a\rho^{\gamma}$ with $\gamma>1$ and when $\nu\mu(\rho)\leq2\mu(\rho)+N\lambda(\rho)\leq\frac{1}{\nu}\mu(\rho)$ with $\nu>0$:
$$(\rho-1)\in L^{\infty}(L^{\gamma}_{2}),\;\;\mbox{and}\;\;\rho |u|^{2}\in L^{\infty}((0,+\infty),L^{1}(\R^{N})),$$
$$\sqrt{\kappa(\rho)}\n\rho\in L^{\infty}((0,+\infty),L^{2}(\R^{N}))^{N},\;\;\mbox{and}\;\;\sqrt{\mu(\rho)}D u,\in
L^{2}((0,+\infty)\times \R^{N})^{N^{2}}.$$
We refer to \cite{fL2} for the definition of the Orlicz spaces. Let us now recall the notion of scaling for the Korteweg's system (\ref{3systeme}). Such an
approach is now classical for incompressible Navier-Stokes equation
and yields local well-posedness (or global well-posedness for small
data) in spaces with minimal regularity.
In our situation we can
easily check that, if $(\rho,u)$ solves (\ref{3systeme}), then
$(\rho_{\lambda},u_{\lambda})$ solves also this system:
$$\rho_{\lambda}(t,x)=\rho(\lambda^{2}t,\lambda x)\,,\,u_{\lambda}(t,x)=\lambda u(\lambda^{2}t,\lambda x)$$
provided the pressure laws $P$ have been changed into
$\lambda^{2}P$.
\begin{definition}
We say that a functional space is critical with respect to the
scaling of the equation if the associated norm is invariant under
the transformation:
$$(\rho,u)\longrightarrow(\rho_{\lambda},u_{\lambda})$$
(up to a constant independent of $\lambda$).
\label{scal}
\end{definition}
This suggests us to choose initial data
$(\rho_{0},u_{0})$ in spaces whose norm is invariant for
all $\lambda>0$ by the transformation
$(\rho_{0},u_{0})\longrightarrow(\rho_{0}(\lambda\cdot),\lambda
u_{0}(\lambda\cdot)).$ A natural candidate is the Besov space (see the section \ref{section2} for some definitions of Besov spaces)
$B^{N/2}_{2,\infty}\times (B^{N/2-1}_{2,\infty})^{N}$,
however since $B^{N/2}_{2,\infty}$ is not included in $L^{\infty}$, we cannot
expect to get a priori $L^{\infty}$ estimate on the density when
$\rho_{0}\in B^{N/2}_{2,\infty}$ (in particular it makes the study of the non linear term delicate since it appears impossible to use composition theorems). Typically except in the case of specific pressure we are not able to deal with this last one. An other difficulty concerns the control of the vacuum or more precisely of the $L^{\infty}$ norm of $\frac{1}{\rho}$, indeed it is crucial to avoid vacuum in order to involve the parabolicity of the momentum equation.\\
That is why an other candidate of initial data space would be $(B^{N/2}_{2,\infty}\cap L^{\infty})\times (B^{N/2-1}_{2,\infty})^{N}$, it will be typically in this last type of space that we shall prove some result of global strong solution by considering specific physical coefficients. Indeed controlling the $L^{\infty}$ norm of the density for general physical coefficients  is generally an hard task in fluid mechanics, even if the Korteweg system allows to obtain regularizing effects on the density. This is the reason why in the literature the authors consider initial density which are in Banach spaces imbedded in $L^{\infty}$ and such that the regularity is propagate via the parabolic structure of the system.\\
Let us briefly mention that the existence of strong solutions for $N\geq2$ is known since the works by H. Hattori and D. Li \cite{fH1,fH2}. R. Danchin and B. Desjardins in \cite{fDD} improve this result by working in critical spaces for the scaling of the equations, more precisely the initial data $(\rho_{0}-1,\rho_{0}u_{0})$
belong to $B^{\frac{N}{2}}_{2,1}\times B^{\frac{N}{2}1}_{2,1}$ (as mentioned previously, let us emphasize on the fact that $B^{\frac{N}{2}}_{2,1}$ is embedded in $L^{\infty}$and it  plays a crucial role in order to control the vacuum and the $L^{\infty}$ norm of the density). In \cite{3MK}, M. Kotschote showed the existence of strong solution for the isothermal model in bounded domain by using Dore\^a-Venni Theory and $\mathcal{H}^{\infty}$ calculus.
In \cite{fH1},  we generalize the results of \cite{fDD} in the case of non isothermal Korteweg system with physical coefficients depending on the density and the temperature.%\\
\subsection{Mathematical results}
We now are going to state our main results. As we explained previously, one of the main difficulty to obtain strong solutions for Korteweg's system in very general Besov spaces consists in dealing with the $L^{\infty}$ control the vacuum $\frac{1}{\rho}$ and of the density $\rho$ without simply propagating the regularity of a Banach space embedded in $L^{\infty}$. To do this we shall precisely understand the structure of the equations by developing new arguments related to the maximum principle.\\
Before explaining how to prove such results,  let us start with observing specific structure on the Korteweg system when $\kappa(\rho)=\frac{\kappa}{\rho}$ with $\kappa>0$, indeed it is possible to rewrite the system (\ref{3systeme}) in a simple way by introducing an effective velocity, see the system (\ref{3systeme2}) (we refer also to \cite{J, Hprepa} where the authors prove the existence of global weak solutions). More precisely as we shall see in system (\ref{NHV1}) with such capillarity coefficients we are able to write the Korteweg system with a density unknown depending only on $\ln\rho$ (web refer to (\ref{NHV1}) for more details).%. Indeed in the most general case, we have to consider parabolic equations for the velocity and for the density with variable coefficients  what requires a least a control on the $L^{\infty}$ norms of $\rho$ and on $\frac{1}{\rho}$ in  $L^{\infty}$ in order to provide regularizing effects on the density and the velocity.  To avoid these main restrictions, we will work with a specific choice on the capillarity coefficients, more precisely we are going to work with 
\begin{remarka}
Let us give some explanations on this choice of capillarity $\kappa(\rho)=\frac{\kappa}{\rho}$, indeed this regime flows  exhibits particular phenomena in the case of the compressible Korteweg Euler system (which is called quantum compressible Euler system when $\kappa(\rho)=\frac{\kappa}{\rho}$). At least heuristically, the system is equivalent via the Madelung transform to the Gross-Pitaevskii equations which are globally well-posed for large initial data in dimension $N=1,2,3$ (we refer to \cite{Gerard}). One of the main difficulty to pass from Gross-Pitaevskii to Quantic Euler consists in dealing with the vacuum. This is one of the reasons why the mathematical community is interested in building solitons for this type of problem (one of the main other reasons correspond to give a negative answer to the problem of scattering and consequently to study the stability of the soliton, we refer to \cite{PG}). We would also like to mention very interesting results of global weak solutions for the compressible quantic Euler equation with a regime $\kappa(\rho)= \frac{1}{\rho}$ due to Antonelli and Marcati (see \cite{Antonelli}). To finish let also mention that we prove in \cite{CH1} that it exits global strong solutions in one dimension for the system (\ref{3systeme}) when $\kappa(\rho)=\frac{\mu^{2}}{\rho}$ and $\mu(\rho)=\mu\rho$, $\lambda(\rho)=0$ and furthermore these solutions converge to a global weak entropy solution of the compressible Euler system when $\mu$ goes to $0$. It shows in particular that the Korteweg system is relevant to select the physical solution of the compressible Euler system via a viscosity-capillarity vanishing process.
\end{remarka}
%In the sequel we are going to study the case of the following capillarity $\kappa(\rho)=\frac{\kappa}{\rho}$ and with the following viscosity $\mu(\rho)=\mu\rho$.
When $\kappa(\rho)=\frac{\kappa}{\rho}$, we rewrite the capillarity tensor as follows (see the appendix for more details on the computations):
$$
\begin{aligned}
K(\rho)=&\kappa\rho(\n\D(\ln\rho)+\frac{1}{2}\n(|\n\ln\rho|^{2})).
\end{aligned}
$$
%We refer to the appendix in section \ref{appendix} for more details on the formal computations on the tensor $K$. 
We now want to consider the eulerian form of the system (\ref{3systeme}) when $\kappa(\rho)=\frac{\kappa}{\rho}$ with $\kappa>0$, we obtain then by dividing by $\rho$ the momentum equation the following system:
\begin{equation}
\begin{aligned}
\begin{cases}
&\p_{t}\ln\rho+u\cdot\n\ln\rho+{\rm div}u=0,\\
&\p_{t}u+u\cdot\n u-\frac{1}{\rho}{\rm div}(2\mu(\rho)D u)-\frac{1}{\rho}
\n(\lambda(\rho){\rm div}u)+\n F(\rho)=\kappa\n\D(\ln\rho)\\
&\hspace{10cm}+\frac{\kappa}{2}\n(|\n\ln\rho|^{2}),\\
&(\ln\rho,u)_{\ t=0}=(\ln\rho_{0},u_{0}).
\end{cases}
\end{aligned}
\label{NHV}
\end{equation}
with $F(\rho)$ defined by $\frac{F^{'}(\rho)}{\rho}=P^{'}(\rho)$. In the sequel we will use the following definition.
\begin{definition}
We now set:
$$q=\ln\rho.$$
\end{definition}
One can now state the main results of the paper. In the first  theorem we prove the existence of global strong solution for (\ref{NHV}) with {\it small} initial data and of strong solution in finite time with large initial data when we choose specific viscosity coefficients and pressure terms. More precisely we shall deal with the shallow water viscosity coefficients:%We also prove the existence of strong solution in finite time with large initial data. The interest of this result will be to consider initial data with discontinuous interfaces.\\
%To do this, we will  assume that: 
$$\mu(\rho)=\mu\rho,\; \lambda(\rho)=\lambda\rho,$$
 with $\mu>0$ and $2\mu+\lambda> 0$. It leads to the following system:
\begin{equation}
\begin{aligned}
\begin{cases}
&\p_{t}q+u\cdot\n q+{\rm div}u=0,\\
&\p_{t}u+u\cdot\n u-\mu\D u-
2\mu\n q\cdot D(u)-(\lambda+\mu)\n {\rm div}u-\lambda{\rm div}u \n q+\n F(\rho)\\
&\hspace{8cm}=\kappa\n\D q+\frac{\kappa}{2}\n(|\n q|^{2}),\\
&(q,u)_{\ t=0}=(\ln\rho_{0},u_{0}).
\end{cases}
\end{aligned}
\label{NHV1}
\end{equation}
Let us now state our main results, we refer to the section \ref{section2} for the definitions of the Besov space and the Hybrid Besov spaces.
%with ${\cal A}u=\mu\D u-(\lambda+\mu)\n{\rm div}u$.
\begin{theorem}
\label{ftheo1}Let $N\geq2$. Assume that $\mu(\rho)=\mu\rho$, $\lambda(\rho)=\lambda\rho$ with $\mu>0$, $2\mu+\lambda> 0$ and $P(\rho)=K\rho$ with $K>0$. We also suppose that:
$$q_{0}\in B^{\N}_{2,\infty}\;\;\mbox{and}\;\;u_{0}\in B^{\N-1}_{2,\infty}.$$
There exists  a time $T$ such that (\ref{NHV1}) has a unique  solution  $(q,u)$ on $(0,T)$ with:
$$
\begin{aligned}
&q\in\widetilde{L}^{\infty}_{T}(B^{\N}_{2,\infty})\cap\widetilde{L}_{T}^{1}(B^{\N+2}_{2,\infty}),\;\;\mbox{and}\;\;u\in\widetilde{L}_{T}^{\infty}(B^{\N-1}_{2,\infty})
\cap\widetilde{L}_{T}^{1}(B^{\N+1}_{2,\infty}).
%&\hspace{7cm}\mbox{and}\;\;\mathbb{Q}u\in\widetilde{L}^{\infty}(B^{\N-1}_{2,1})\cap\widetilde{L}^{1}(B^{\N+1}_{2,1}).
\end{aligned}
$$
Furthermore it exists $\e_{0}$ such that if in addition $q_{0}\in B^{\N-1}_{2,\infty}$ and:
\begin{equation}
\|q_{0}\|_{\widetilde{B}^{\N-1,\N}_{2,\infty}}+\|u_{0}\|_{B^{\N-1}_{2,\infty}}\leq \e_{0}.
\label{petit}
\end{equation}
then the solution $(q,u)$ is global with:
\begin{equation}
\begin{aligned}
&q\in\widetilde{L}^{\infty}(\widetilde{B}^{\N-1,\N}_{2,\infty})\cap\widetilde{L}^{1}(\widetilde{B}^{\N+1,\N+2}_{2,\infty}),\;\;\mbox{and}\;\;u\in\widetilde{L}^{\infty}(B^{\N-1}_{2,\infty})
\cap\widetilde{L}^{1}(B^{\N+1}_{2,\infty}).
%&\hspace{7cm}\mbox{and}\;\;\mathbb{Q}u\in\widetilde{L}^{\infty}(B^{\N-1}_{2,1})\cap\widetilde{L}^{1}(B^{\N+1}_{2,1}).
\end{aligned}
\label{impra}
\end{equation}
If in addition of the hypothesis (\ref{petit}) we assume that $(q_{0},u_{0})$ belongs in $\widetilde{B}^{\N-1,\N}_{2,2}\times B^{\N-1}_{2,2}$ and that $q_{0}\in L^{\infty}$ then there exists a unique solution $(\rho,u)$ of the system (\ref{3systeme}) with $q=\ln\rho$. Furthermore for any $T>0$ it exists $C_{T}>0$ depending on $T$ such that:
$$\|\frac{1}{\rho}\|_{L^{\infty}_{T}(L^{\infty})}+\|\rho\|_{L^{\infty}_{T}(L^{\infty})}\leq C_{T}.$$
In addition it implies that for any $T>0$:
\begin{equation}
\begin{aligned}
&q\in\widetilde{L}^{\infty}_{T}(B^{\N}_{2,2})\cap\widetilde{L}^{1}(\widetilde{B}^{\N+1,\N+2}_{2,2}),\;\;\mbox{and}\;\;u\in\widetilde{L}_{T}^{\infty}(B^{\N-1}_{2,2})
\cap\widetilde{L}^{1}(B^{\N+1}_{2,2}).
&\hspace{7cm}\mbox{and}\;\;\mathbb{Q}u\in\widetilde{L}^{\infty}(B^{\N-1}_{2,1})\cap\widetilde{L}^{1}(B^{\N+1}_{2,1}).
\end{aligned}
\label{impra}
\end{equation}
\end{theorem}
%\begin{remarka}
%We would like to mention that for the existence of strong solution in finite time, we could obtain the same kind of results in non homogeneous Besov space. We now want also to emphasize on the fact that when we are working with critical Besov spaces for the initial data, we generally are not able to give an estimate of the time of existence $T$. However when we work with subcritical non homogeneous Besov space for the initial data, we are able to show that the time of existence for system 8\ref{NHV1}) verifies the following inequality:
%$$T^{'}\geq \frac{C}{(1+\|u_{0}\|_{\dot{B}^{\NN-1+\e}_{p,\infty}}+\|\ln\rho_{0}\|_{\dot{B}^{\NN+\e}_{p,\infty}})^{\beta}},$$
%with $1\leq p<2N$ and $C,\beta$ depending on the physical coefficients and on the dimension. To do this, it is just a matter to show how $(q_{L},u_{L})$ are small in the critical Besov spaces in function of the time and of the initial data. For more details we refer to the proof of theorem \ref{ftheo1}.
%\end{remarka}
\begin{remarka}
We want to point out that in the first part of the theorem \ref{ftheo1}, we solve the system (\ref{NHV1}) and not the system (\ref{3systeme}). Indeed we do not assume any control on $q_{0}$ in $L^{\infty}$, it means that we have no information on the vacuum of the density. So it is not clear that a solution from (\ref{NHV1}) is also a solution from (\ref{3systeme}) when there is vacuum. This result proves in a certain way that the good variable to consider is not the density $\rho$ but rather $\ln\rho$. Let us emphasize on the fact that this result is the first (up our knowledge) for compressible system where we are able to work with such general critical initial data (in particular $u_{0}\in B^{\N-1}_{2,\infty}$ what is a classical regularity for incompressible Navier-Stokes equations, see \cite{CP}).\\ %Generally the results of global strong solution for compressible system involve initial density and velocity in $B^{\N}_{2,1}\times B^{\N-1}_{2,1}$ in order to ensure a Lipschitz control on the velocity and a $L^{\infty}$ bound on the density (see \cite{fDD}).In a certain sense, this result shows hat the Korteweg system with such viscosity and capillary coefficients is close from the incompressible Navier-Stokes equations in terms of critical regularity for the initial data.\\
In the second part of theorem we provide additional hypothesis on the density since $q_{0}=\ln\rho_{0}$ belongs also in $L^{\infty}\cap B^{\N}_{22}$ and $u_{0}\in B^{\N-1}_{22}$, in counter part we are able to show that with such initial data we have existence of strong solution for the "real" Korteweg system (\ref{3systeme}). The main difficulty consists in proving $L^{\infty}$ estimates on $q$ by splitting the solution $q$ in the sum of the solution of the linearized system $q_{L}$ and a remainder term $\bar{q}$ taking into account the non linear terms. We combine some maximum principle results on $q_{L}$ and regularizing effect on $\bar{q}$ on the third index of the Besov space to show that $q$ is well controlled in $L^{\infty}$. This part is quite technical and involves new ideas for developing some results of maximum principle type, to do this we show in particular a characterization of Besov spaces in terms of the semi group associated to the linear system related to (\ref{3systeme}).\\
\\
Let us finish by observing that in this theorem we authorize that the initial density is not necessary continuous as in \cite{fDD}, this is one of the main motivation of this paper.
Indeed we can deal with initial density which have some jump, we observe that the density is immediately regularized in the sense that $\rho$ is in $C^{\infty}((0,T),\R^{N})$ for any $T>0$. This is due to the fact that the interfaces are diffuse.
\label{remcrucim}
\end{remarka}
\begin{remarka}
The second interest of this theorem it to obtain the existence of global strong solution with subcritical smallness on the initial data. More precisely we assume that $(q_{0},u_{0})$ belongs in $(\widetilde{B}^{\N-1,\N}_{2,2}\cap L^{\infty})\times B^{\N-1}_{2,2}$ but we require only a smallness assumption in $\widetilde{B}^{\N-1,\N}_{2,\infty}\times B^{\N-1}_{2,\infty}$. In particular it allows us to prove the existence of global strong solution with large initial data in the energy space when $N=2$, this is up our knowledge the first result of global strong solution with large initial energy data for the Korteweg system when $N=2$.\\
We give a example of such initial data in the corollary \ref{corotech} where $u_{0 ,\e,l_{0}}$ can be chose as small as possible in $B^{\N-1}_{2,\infty}$ but very \textit{large} in $B^{\N-1}_{2,2}$. Let us recall that when $N=2$, $B^{\N-1}_{2,2}$ corresponds to $L^{2}$ which is the energy space for $u_{0}$ (see (\ref{3inegaliteenergie1})).
\\[3mm]
An other interesting point is that compared with \cite{fDD}, we can choose initial density with a large $L^{\infty}$ norm.
\end{remarka}
\begin{remarka}
We now want to point out the specificity  of the different physical coefficients, it is typically the case for the pressure and the viscosity coefficients where $P(\rho)=K\rho$ and $\mu(\rho)=\mu\rho$. Indeed under this form we can check that the variable $\ln\rho$ appears naturally everywhere but also that there is no non-linear terms depending on the density such that it would be necessary to have $L^{\infty}$ control in order to deal with them (typically by using a composition theorem). Let us also mention that it is not possible to obtain the same type of result for compressible Navier-Stokes equations as we have no regularizing effects on the density. %Indeed in this case it would be not possible to deal with the following term $\n\ln\rho\cdot\n u$ and to prove that it remains in $\widetilde{L}^{1}(B^{\N-1}_{2,\infty})$.\\
%In the following theorem we are going to extend the previous analysis to general physical coefficients.%we will prove a result of ill-posedness when we consider general pressure. Indeed in these case we need to control the $L^{\infty}$ norm of the density to deal with the pressure term of type $P(\rho)=\rho^{\gamma}$ with $\gamma$ large enough. In the case where $\gamma$ is not so large, we could extend the theorem \ref{ftheo1} by taking advantage of the regularizing effect to estimate $\n\rho^{\gamma}$ in $\widetilde{L}^{1}(B^{\N-1}_{2,\infty})$. However if we do not assume enough regularity on the divergence of the velocity, we loss the $L^{\infty}$ control on the norm of the density, it will prove in the third theorem.
\end{remarka}
\begin{remarka}
We would like to mention that we could easily extend the result of strong solution in finite time to the framework of Besov spaces constructed on general $L^{p}$ spaces when the initial data verify:
$$q_{0}\in B^{\NN}_{p,\infty}\;\;\mbox{and}\;\;u_{0}\in B^{\NN-1}_{p,\infty},$$ 
Concerning the existence of global strong solution we would have to assume that $q_{0}$ is in $\widetilde{B}^{\N-1\NN}_{2,p,\infty}$ and $u_{0}\in \widetilde{B}^{\N-1,\NN-1}_{2,p,\infty}$, the main difficulty would concern the proof of estimate for the linear system associated to the Korteweg system with such initial data (we refer to \cite{arma} for a such study in the framework of the compressible Navier Stokes equations). 
\end{remarka}
The previous theorem takes into account in an essential way the structure of the viscosity, capillary and pressure coefficients, indeed it allows us to assume a smallness hypothesis only on initial data in $\widetilde{B}^{\N-1,\N}_{2,\infty}\times B^{\N-1}_{2,\infty}$ in order to obtain global strong solution. We would like to extend this result to general physical coefficients and in particular dealing with the case of the constant capillary coefficient in order to improve the results of \cite{fDD}. Let us state a new theorem dealing with  general pressure and with capillary coefficient of the form $\kappa(\rho)=\frac{\kappa}{\rho}$ or $\kappa$ with $\kappa>0$.
\begin{theorem}
\label{ftheo2}Let $N\geq2$. Assume that $\mu(\rho)=\mu\rho$ or $\mu$, $\lambda(\rho)=\lambda\rho$ or $\lambda$ with $\mu>0$, $2\mu+\lambda> 0$, $\kappa(\rho)=\frac{\kappa}{\rho}$ or $\kappa$ with $\kappa>0$ and $P$ a regular function such that $P^{'}(1)>0$. Furthermore  we suppose that $\rho_{0}=1+h_{0}$:
$$h_{0}\in \widetilde{B}^{\N-1,\N}_{2,2}\cap B^{\N-2}_{2,1}\cap L^{\infty}\;\;\mbox{and}\;\;u_{0}\in B^{\N-2}_{2,1}\cap B^{\N-1}_{2,2}.$$
%There exists  a time $T$ such that (\ref{3systeme}) has a unique  solution  $(\rho,u)$ on $(0,T)$ with:
There exists $\e_{0}$ such that if:
$$\|h_{0}\|_{\widetilde{B}^{\N-1,\N}_{2,2}\cap B^{\N-2}_{2,1}\cap L^{\infty}}+\|u_{0}\|_{B^{\N-2}_{2,1}\cap B^{\N-1}_{2,2}}\leq \e_{0}.$$
then it exists a global unique solution $(\rho,u)$ of the system (\ref{3systeme}) with $\rho=1+h$ and:
\begin{equation}
\begin{aligned}
&h\in\widetilde{L}^{\infty}(\widetilde{B}^{\N-1,\N}_{2,2}\cap B^{\N-2}_{2,1})\cap\widetilde{L}^{1}(\widetilde{B}^{\N+1,\N+2}_{2,2}\cap B^{\N}_{2,1})\cap L^{\infty}(L^{\infty})\\
&\;\;\;\;\;\;\;\;\;\;\mbox{and}\;\;u\in\widetilde{L}^{\infty}(B^{\N-1}_{2,2}\cap B^{\N-2}_{2,1})
\cap\widetilde{L}^{1}(B^{\N+1}_{2,2}\cap B^{\N}_{2,1}).
%&\hspace{7cm}\mbox{and}\;\;\mathbb{Q}u\in\widetilde{L}^{\infty}(B^{\N-1}_{2,1})\cap\widetilde{L}^{1}(B^{\N+1}_{2,1}).
\end{aligned}
\label{aimpra}
\end{equation}
%If in addition with the previous smallness hypothesis we assume that $(q_{0},u_{0})$ belongs in $B^{\N}_{2,2}\times B^{\N-1}_{2,2}$ and that $q_{0}\in L^{\infty}$ and $(\D)^{-1}{\rm div}u_{0}\in  L^{\infty}$ then there exists a unique solution $(\rho,u)$ of the system (\ref{3systeme}) with $q=\ln\rho$ and for any $T$ it exists $C_{T}$ such that:
%$$\|\frac{1}{\rho}\|_{L^{\infty}_{T}(L^{\infty})}+\|\rho\|_{L^{\infty}_{T}(L^{\infty})}\leq C_{T}.$$
\end{theorem}
\begin{remarka}
We would like to mention that this theorem extend the results of \cite{fDD} in terms of rough regularity on the initial data. Here we assume only that $u_{0}$ belongs in $B^{\N-1}_{2,2}$ instead $B^{\N-1}_{2,1}$, the main task as in the previous theorem \ref{ftheo1} consists in getting a control of the density in $L^{\infty}$ without assuming $h_{0}\in B^{\N}_{2,1}$ as in \cite{fDD}. To do this we have to study precisely the linearized system of (\ref{3systeme}) and involving some results of maximum principle type, in addition the solution are the sum of solution for the linear system associated to (\ref{3systeme}) and remainder terms taking into account the non linearities. These remainders terms benefit from regularizing effects on the third index of the Besov space which will ensure a $L^{\infty}$ control on the density. For more details we refer to the proof of theorem \ref{ftheo2}.\\
Let us mention that for the first time up my knowledge we can choose initial density which are not continuous as in \cite{fDD}, in particular we authorize some jump on the initial density which is immediately regularized. 
 \end{remarka}
\begin{remarka}
Let us point out that in the previous theorem \ref{ftheo2} we ask additional regularity in low frequencies on $(h_{0},u_{0})$, indeed we assume that $(h_{0},u_{0})$ are in $B^{\N-2}_{2,1}\times (B^{\N-2}_{2,1})^{N}$. This is essentially due to the fact that we have shown \textit{maximum principle} for the system $(N)$ p18 (which does not take into account the low frequencies) and not for the system $(N1)$ p 20 (see the proposition \ref{maxfinal}). We think that it would be possible to extend the proposition \ref{maxfinal} to the system $(N1)$. In particular it will allow to avoid this additional regularity on the initial data in low frequencies. For more explanations we refer to the remark \ref{remtecg}.
\end{remarka}
We are now interested in dealing with the specific case $\kappa(\rho)=\frac{\mu^{2}}{\rho}$ and $\mu(\rho)=\mu\rho$, $\lambda(\rho)=0$ where we can exhibit a specific structure on the system (\ref{3systeme}), we obtain from system (\ref{NHV1}) (when $\lambda=0$) the following simplified model by assuming that $v=u+\mu\n\ln\rho$ (we refer to the appendix for more details on the computation or \cite{Hprepa})):
\begin{equation}
\begin{cases}
\begin{aligned}
&\p_{t}\rho-\mu\D\rho=-{\rm div}(\rho v),\\
&\rho\p_{t}v +\rho u\cdot\n v-\rm div(\mu\rho\,\n v)+\n P(\rho)=0,
\end{aligned}
\end{cases}
\label{3systeme2}
\end{equation}
which is equivalent to the following system if we control the vacuum on the density $\rho$ with $q=\ln\rho$:
\begin{equation}
\begin{cases}
\begin{aligned}
&\p_{t}q-\mu\D q+v\cdot\n q=-{\rm div}v+\mu|\n q|^{2},\\
&\p_{t}v +u\cdot\n v-\mu\D v-\mu\n q\cdot\n v+\n F(\rho)=0,
\end{aligned}
\end{cases}
\label{3systeme2ab}
\end{equation}
In this particular case we are able to use a new tool developed in \cite{cras1,cras2,arxiv,arxiv1,hal} called the quasi-solutions. More precisely we can check that it exists a particular solution  of the following system (where we have canceled out the pressure $P$):
\begin{equation}
\begin{cases}
\begin{aligned}
&\p_{t}\rho-\mu\D\rho=-{\rm div}(\rho v),\\
&\rho\p_{t}v +\rho u\cdot\n v-\rm div(\mu\rho\,\n v)=0,
\end{aligned}
\end{cases}
\label{a3systeme2}
\end{equation}
Indeed we verify that $(\rho_{1},-\mu\n\ln\rho_{1})$ is a particular solution of (\ref{a3systeme2}) if the density $\rho_{1}$ verifies the following heat equation:
\begin{equation}
\p_{t}\rho_{1}-\mu\D\rho_{1}=0.
\label{eqchaleur}
\end{equation}
The idea will consists in working around this quasi-solution which allows us to prove the existence of global strong solution with small initial data in subcritical norms. Indeed we are going to search solution under the form:
$$q=\ln\rho=\ln\rho_{1}+h_{2}\;\;\mbox{with}\;\rho_{1}=1+h_{1}\;\;\mbox{and}\;\;u=-\mu\n\ln\rho_{1}+u_{2}.$$ 
We deduce from (\ref{3systeme2ab}) that $(h_{2},u_{2})$ verifies the following system:
\begin{equation}
\begin{cases}
\begin{aligned}
&\p_{t}h_{2}+{\rm div}u_{2}-\mu\n\ln\rho_{1}\cdot\n h_{2}+u_{2}\cdot\n\ln\rho_{1}=F(h_{2},u_{2}),\\[2mm]
&\p_{t}u_{2}-\mu\D u_{2}-\mu\n{\rm div}u_{2}-\kappa\n\D h_{2}+K\n h_{2}-2\mu\n \ln\rho_{1}\cdot D u_{2}-2\mu\n h_{2}\cdot D u_{1}\\
&\hspace{4cm}+u_{1}\cdot\n u_{2}+u_{2}\cdot\n u_{1}-\mu^{2}\n(\n\ln\rho_{1}\cdot\n h_{2})=G(h_{2},u_{2}),\\
&(h_{2}(0,\cdot),u_{2}(0,\cdot))=(h_{0}^{2},u_{0}^{2}).
\end{aligned}
\end{cases}
\label{fond3systeme2}
\end{equation}
with:
\begin{equation}
\begin{aligned}
&F(h_{2},u_{2})=-u_{2}\cdot\n h_{2},\\
&G(h_{2},u_{2})=-u_{2}\cdot\n u_{2}+2\mu\n h_{2}\cdot D u_{2}-K\n\ln\rho_{1}+\frac{\mu^{2}}{2}\n(|\n h_{2}|^{2}).
\end{aligned}
\label{restetech}
\end{equation}
Let us state our main theorem for the system (\ref{fond3systeme2}).
\begin{theorem}
\label{ftheo3}Let $N\geq2$. Assume that $\mu(\rho)=\mu\rho$, $\kappa(\rho)=\frac{\mu^{2}}{\rho}$ and $\lambda(\rho)=0$ with $\mu>0$ and $P(\rho)=K\rho$  with $K>0$. Furthermore  we suppose that $u_{0}=-\mu\n[\ln\rho_{1}^{0}]+u_{2}^{0}$ and $\ln\rho_{0}=\ln(\rho^{0}_{1})+h_{2}^{0}$ such that $\rho_{1}^{0}=1+h_{1}{0}$ and it exits $c_{1}>0$ such that $\rho_{1}^{0}\geq c_{1}>0$. In addition we assume that:
$$h_{1}^{0}\in \widetilde{B}^{\N-2,\N}_{2,1},\;\;h_{2}^{0}\in \widetilde{B}^{\N-1,\N}_{2,1}\;\;\mbox{and}\;\;u_{2}{0}\in B^{\N-1}_{2,1}.$$
%Assuming that $\rho_{0}\geq c>0$ there exists  a time $T$ such that (\ref{3systeme}) has a unique  solution  $(\rho,u)$ on $(0,T)$ with:
%$$
%\begin{aligned}
%&\rho\in\widetilde{L}^{\infty}_{T}(B^{\N}_{2,\infty})\cap\widetilde{L}_{T}^{1}(B^{\N+2}_{2,\infty})\cap L^{\infty}_{T}(L^{\infty}),\;\;\mbox{and}\;\;u\in\widetilde{L}_{T}^{\infty}(B^{\N-1}_{2,\infty})
%\cap\widetilde{L}_{T}^{1}(B^{\N+1}_{2,\infty}).
%&\hspace{7cm}\mbox{and}\;\;\mathbb{Q}u\in\widetilde{L}^{\infty}(B^{\N-1}_{2,1})\cap\widetilde{L}^{1}(B^{\N+1}_{2,1}).
%\end{aligned}
%$$
Furthermore it exists $C>0$, $\e_{0}$ (depending on $h^{0}_{1}$),  and two regular function $g$, $g_{1}$ such that if:
\begin{equation}
\begin{aligned}
&C g(\|(\rho^{0}_{1},\frac{1}{\rho^{0}_{1}})\|_{L^{\infty}})\|h^{0}_{1}\|_{
B^{\N-2}_{2,1}} \exp\big(Cg_{1}(\|(\rho^{0}_{1},\frac{1}{\rho^{0}_{1}})\|_{L^{\infty}})\|h^{0}_{1}\|_{
B^{\N}_{2,1}})\big)\leq\frac{1}{2}
,\\
&\|h_{2}^{0}\|_{\widetilde{B}^{\N-1,\N}_{2,1}}+\|u_{2}^{0}\|_{B^{\N-1}_{2,1}}\leq \e_{0}.
\end{aligned}
\label{crucinitial}
\end{equation}
then it exists a global unique solution $(\rho,u)$ of the system (\ref{3systeme}) with: $u=-\mu\n\ln\rho_{1}+u_{2}$ and $\ln\rho=\ln\rho_{1}+h_{2}$ with $\rho_{1}=1+h_{1}$ verifying the following system:
$$
\begin{cases}
\begin{aligned}
&\p_{t}\rho_{1}-\mu\D\rho_{1}=0,\\
&\rho_{1}(0,\cdot)=\rho_{1}^{0}=1+h_{1}^{0}.
\end{aligned}
\end{cases}
$$
Furthermore we have:
\begin{equation}
\begin{aligned}
&h_{2}\in\widetilde{L}^{\infty}(\widetilde{B}^{\N-1,\N}_{2,1})\cap\widetilde{L}^{1}(\widetilde{B}^{\N+1,\N+2}_{2,1})\;\;\mbox{and}\;\;u_{2}\in\widetilde{L}^{\infty}(B^{\N-1}_{2,1})
\cap\widetilde{L}^{1}(B^{\N+1}_{2,1}).
%&\hspace{7cm}\mbox{and}\;\;\mathbb{Q}u\in\widetilde{L}^{\infty}(B^{\N-1}_{2,1})\cap\widetilde{L}^{1}(B^{\N+1}_{2,1}).
\end{aligned}
\label{aimpra}
\end{equation}
%If in addition with the previous smallness hypothesis we assume that $(q_{0},u_{0})$ belongs in $B^{\N}_{2,2}\times B^{\N-1}_{2,2}$ and that $q_{0}\in L^{\infty}$ and $(\D)^{-1}{\rm div}u_{0}\in  L^{\infty}$ then there exists a unique solution $(\rho,u)$ of the system (\ref{3systeme}) with $q=\ln\rho$ and for any $T$ it exists $C_{T}$ such that:
%$$\|\frac{1}{\rho}\|_{L^{\infty}_{T}(L^{\infty})}+\|\rho\|_{L^{\infty}_{T}(L^{\infty})}\leq C_{T}.$$
\end{theorem}
\begin{remarka}
Let us mention than the main interest of this theorem is to prove the existence of global strong solution with large initial data for the scaling of the equation which is completely new up our knowledge. Indeed it suffices to choose $h_{1}^{0}(x)=\va(\lambda x)$ with $\va\in \widetilde{B}^{\N-2,\N}_{2,1}$ such that $1+\va\geq c>0$ we then verify easily that:
\begin{equation}
\begin{cases}
\begin{aligned}
&\|h_{1}^{0}\|_{ B^{\N}_{2,1}}=\|\va\|_{B^{\N}_{2,1}},\\
&\|\rho_{1}^{0}\|_{L^{\infty}}=\|1+\va\|_{L^{\infty}},\\
&\|\frac{1}{\rho_{1}^{0}}\|_{L^{\infty}}=\|\frac{1}{1+\va}\|_{L^{\infty}},\\
&\|h_{1}^{0}\|_{ B^{\N-2}_{2,1}}=\frac{1}{\lambda^{2}}\|\va\|_{B^{\N-2}_{2,1}}.
\end{aligned}
\end{cases}
\label{genial}
\end{equation}
It implies that  $h_{1}^{0}$ verifies (\ref{crucinitial}) by choosing $\lambda$ large enough. In particular it implies that by taking $\va$ large in $B^{\N}_{2,\infty}$ our initial density $h_{1}^{0}$ is large in the critical Besov space $B^{\N}_{2,\infty}$ for the scaling of the equations. It is also possible to choose $\va$ large in $B^{1}_{2,\infty}$ which shows that there is existence of global strong solution for large initial data in the energy space when $N=2$. This is certainly the main interest of this paper.\\% and this is the first time up our knowledge that we can obtain global strong solution with large initial data in the energy space when $N=2$ for the Korteweg system.\\[3mm]
We could also to choose $h_{1}^{0}(x)=\ln(\lambda)\va(\lambda x)$ with $\lambda>0$ which shall still improve the previous estimates in term of \textit{large} initial data in $B^{\N}_{2,2}$.\\% let us also emphasize that we can take $\va$ large in $B^{\N}_{2,\infty}$, in particular it implies that $h^{0}_{1}$ is large in $B^{\N}_{2,\infty}$ however we obtain the existence of global strong solution, it improves from this point of view the results of theorem \ref{ftheo1}.
We would like also to point out that the condition (\ref{crucinitial}) is in the same spirit that smallness condition in some forks of Chemin and Gallagher in \cite{CG1,CG2} for incompressible Navier-Stokes equations. Indeed in these works the authors prove the existence of global strong solution for large initial data in $B^{-1}_{\infty,\infty}$ which is the largest critical space for the Navier-Stokes equations.
\end{remarka}
\begin{remarka}
We think that we could improve the condition on initial density $h_{1}^{0}\in B^{\N-2}_{2,1}$ by working around the quasi solution only in high frequencies and by working in low frequencies directly with $\rho-1$ in the spirit of \cite{arma}. It would be then possible to assume that $h_{1}^{0}$ belongs only to $B^{\N-1}_{2,1}$.\\
Let also point out that in this theorem, by a very accurate study on the linear system associated to the system (\ref{fond3systeme2}) (see the proposition \ref{fcrucialprop}) we could probably improve the initial condition on $h_{1}^{0}$ by assuming only $h_{1}^{0}$ in $B^{\N-2}_{2,1}\cap B^{\N}_{2,2-\e}\cap L^{\infty}$ with $\e>0$ (see the remark
\ref{importtech44}). In particular it would allow to choose initial density not necessary continuous.
\\
It may be also probably possible to work with general regular pressure, it would make the proof technically more difficult.
\end{remarka}
We are going to finish by presenting a result of global strong solution with large initial data when we assume the system (\ref{3systeme}) highly compressible.
\begin{corollaire}
\label{cor3}
Let $N\geq2$. Assume that $\mu(\rho)=\mu\rho$, $\kappa(\rho)=\frac{\mu^{2}}{\rho}$ and $\lambda(\rho)=0$ with $\mu>0$ and $P(\rho)=K\rho$  with $K>0$. Furthermore  we suppose that $u_{0}=-\mu\n[\ln\rho_{1}^{0}]+u_{2}^{0}$ and $\ln\rho_{0}=\ln(\rho^{0}_{1})+h_{2}^{0}$ such that $\rho_{1}^{0}=1+h_{1}^{0}$ and it exits $c_{1}>0$ such that $\rho_{1}^{0}\geq c_{1}>0$. In addition we suppose that:
$$h_{1}^{0}\in \widetilde{B}^{\N-2,\N}_{2,1},\;\;h_{2}^{0}\in \widetilde{B}^{\N-1,\N}_{2,1}\;\;\mbox{and}\;\;u_{2}^{0}\in B^{\N-1}_{2,1}.$$
%Assuming that $\rho_{0}\geq c>0$ there exists  a time $T$ such that (\ref{3systeme}) has a unique  solution  $(\rho,u)$ on $(0,T)$ with:
%$$
%\begin{aligned}
%&\rho\in\widetilde{L}^{\infty}_{T}(B^{\N}_{2,\infty})\cap\widetilde{L}_{T}^{1}(B^{\N+2}_{2,\infty})\cap L^{\infty}_{T}(L^{\infty}),\;\;\mbox{and}\;\;u\in\widetilde{L}_{T}^{\infty}(B^{\N-1}_{2,\infty})
%\cap\widetilde{L}_{T}^{1}(B^{\N+1}_{2,\infty}).
%&\hspace{7cm}\mbox{and}\;\;\mathbb{Q}u\in\widetilde{L}^{\infty}(B^{\N-1}_{2,1})\cap\widetilde{L}^{1}(B^{\N+1}_{2,1}).
%\end{aligned}
%$$
Furthermore it exists $\e_{0}>0$ (depending on $h_{1}^{0}$ and the viscosity coefficient $\mu$) such that for any $K\leq \e_{0}$ it exits $\e_{1}>0$ such that if
\begin{equation}
\begin{aligned}
&\|h_{2}^{0}\|_{\widetilde{B}^{\N-1,\N}_{2,1}}+\|u_{2}^{0}\|_{B^{\N-1}_{2,1}}\leq \e_{1}.
\end{aligned}
\label{bcrucinitial}
\end{equation}
then it exists a global unique solution $(\rho,u)$ of the system (\ref{3systeme}) with: $u=-\mu\n\ln\rho_{1}+u_{2}$ and $\ln\rho=\ln\rho_{1}+h_{2}$ with $\rho_{1}=1+h_{1}$ verifying the following system:
$$
\begin{cases}
\begin{aligned}
&\p_{t}\rho_{1}-\mu\D\rho_{1}=0,\\
&\rho_{1}(0,\cdot)=\rho_{1}^{0}=1+h_{1}^{0}.
\end{aligned}
\end{cases}
$$
Furthermore we have:
\begin{equation}
\begin{aligned}
&h_{2}\in\widetilde{L}^{\infty}(\widetilde{B}^{\N-1,N}_{2,1})\cap\widetilde{L}^{1}(\widetilde{B}^{\N+1,\N+2}_{2,1})\;\;\mbox{and}\;\;u_{2}\in\widetilde{L}^{\infty}(B^{\N-1}_{2,1})
\cap\widetilde{L}^{1}(B^{\N+1}_{2,1}).
%&\hspace{7cm}\mbox{and}\;\;\mathbb{Q}u\in\widetilde{L}^{\infty}(B^{\N-1}_{2,1})\cap\widetilde{L}^{1}(B^{\N+1}_{2,1}).
\end{aligned}
\label{qaimpra}
\end{equation}
%If in addition with the previous smallness hypothesis we assume that $(q_{0},u_{0})$ belongs in $B^{\N}_{2,2}\times B^{\N-1}_{2,2}$ and that $q_{0}\in L^{\infty}$ and $(\D)^{-1}{\rm div}u_{0}\in  L^{\infty}$ then there exists a unique solution $(\rho,u)$ of the system (\ref{3systeme}) with $q=\ln\rho$ and for any $T$ it exists $C_{T}$ such that:
%$$\|\frac{1}{\rho}\|_{L^{\infty}_{T}(L^{\infty})}+\|\rho\|_{L^{\infty}_{T}(L^{\infty})}\leq C_{T}.$$
%If in addition with the previous smallness hypothesis we assume that $(q_{0},u_{0})$ belongs in $B^{\N}_{2,2}\times B^{\N-1}_{2,2}$ and that $q_{0}\in L^{\infty}$ and $(\D)^{-1}{\rm div}u_{0}\in  L^{\infty}$ then there exists a unique solution $(\rho,u)$ of the system (\ref{3systeme}) with $q=\ln\rho$ and for any $T$ it exists $C_{T}$ such that:
%$$\|\frac{1}{\rho}\|_{L^{\infty}_{T}(L^{\infty})}+\|\rho\|_{L^{\infty}_{T}(L^{\infty})}\leq C_{T}.$$
\end{corollaire}
\begin{remarka}
The main interest of this theorem is to prove the existence of global strong solution for any large initial data provided that $K$ is sufficiently small with $P(\rho)=K\rho$. In other terms
 we get global existence (and uniqueness) for highly compressible fluids in any dimension $N\geq 2$. Up my knowledge it is the first time that we have a result of global strong solution with large initial data in dimension $3$ (under a condition of course of high compressibility, which means that $K$ must be small in function of the initial data and in function of $\mu$). Roughly speaking $K$ tends to be very small when $\|h_{1}^{0}\|_{\widetilde{B}^{\N-2,\N}_{2,1}}$ is very large.
\end{remarka}
%Resultat suivant que l'on va enlever!!!!
This article is structured in the following way, first of all we recall
in the section \ref{section2} some definitions an theorems related to the Littlewood-Paley theory. Next we will
concentrate in the section \ref{section4} on the proof of theorem \ref{ftheo1} and \ref{ftheo2}. In section \ref{section4}, we shall prove the theorem \ref{ftheo3} by introducing the notion of quasi solution which will play a crucial role. In section \ref{section5} we show the corollary \ref{cor3}. We postpone in appendix (see section \ref{appendix}) some technical computation on the capillarity tensor and some extensions on the previous results.
\section{Littlewood-Paley theory and Besov spaces}
\label{section2}
Throughout the paper, $C$ stands for a constant whose exact meaning depends on the context. The notation $A\lesssim B$ means
that $A\leq CB$.
For all Banach space $X$, we denote by $C([0,T],X)$ the set of continuous functions on $[0,T]$ with values in $X$.
For $p\in[1,+\infty]$, the notation $L^{p}(0,T,X)$ or $L^{p}_{T}(X)$ stands for the set of measurable functions on $(0,T)$
with values in $X$ such that $t\rightarrow\|f(t)\|_{X}$ belongs to $L^{p}(0,T)$.
Littlewood-Paley decomposition  corresponds to a dyadic
decomposition  of the space in Fourier variables.
We can use for instance any $\varphi\in C^{\infty}(\R^{N})$ and  $\chi\in C^{\infty}(\R^{N})$ ,
supported respectively in
${\cal{C}}=\{\xi\in\R^{N}/\frac{3}{4}\leq|\xi|\leq\frac{8}{3}\}$ and $B(0,\frac{4}{3})$
such that:
$$\sum_{l\in\mathbb{Z}}\varphi(2^{-l}\xi)=1\,\,\,\,\mbox{if}\,\,\,\,\xi\ne 0,$$
and:
$$\chi(\xi)+\sum_{l\in\mathbb{N}}\varphi(2^{-l}\xi)=1\,\,\,\,\mbox{if}\,\,\,\,\forall \xi.$$
Denoting $h={\cal{F}}^{-1}\varphi$, we then define the dyadic
blocks by:
$$\D_{l}u=\varphi(2^{-l}D)u=2^{lN}\int_{\R^{N}}h(2^{l}y)u(x-y)dy\,\,\,\,\mbox{and}\,\,\,S_{l}u=\sum_{k\leq
l-1}\D_{k}u\,.$$ Formally, one can write that:
$$u=\sum_{k\in\mathbb{Z}}\D_{k}u\,.$$
This decomposition is called homogeneous Littlewood-Paley
decomposition. Let us observe that the above formal equality does
not hold in ${\cal{S}}^{'}(\R^{N})$ for two reasons:
\begin{enumerate}
\item The right hand-side does not necessarily converge in
${\cal{S}}^{'}(\R^{N})$.
\item Even if it does, the equality is not
always true in ${\cal{S}}^{'}(\R^{N})$ (consider the case of the polynomials).
\end{enumerate}
For the non homogeneous decomposition, we define the dyadic blocks as follows:
$$
\begin{aligned}
&\dot{\D}_{l}u=0\;\;\mbox{for}\;\;l\leq -2,\\
&\dot{\D}_{-1}u=\chi(D)u,\\
&\dot{\D}_{l}u=\varphi(2^{-l}D)u=2^{lN}\int_{\R^{N}}h(2^{l}y)u(x-y)dy,\;\;\mbox{for}\;\;l\geq 0.
\end{aligned}
$$
Formally, one can write that:
$$u=\sum_{k\in\mathbb{Z}}\D_{k}u\,.$$
This decomposition is called non homogeneous Littlewood-Paley
decomposition.
\subsection{Homogeneous and non homogeneous Besov spaces and first properties}
\begin{definition}
We denote by ${\cal S}_{h}^{'}$ the space of temperate distribution $u$ such that:
$$\lim S_{j}u_{j\rightarrow+\infty}=0\,\,\,\mbox{in}\;\;{\cal S}^{'}.$$
\end{definition}
\begin{definition}
For
$s\in\R,\,\,p\in[1,+\infty],\,\,q\in[1,+\infty],\,\,\mbox{and}\,\,u\in{\cal{S}}^{'}(\R^{N})$
we set:
$$\|u\|_{B^{s}_{p,q}}=(\sum_{l\in\mathbb{Z}}(2^{ls}\|\D_{l}u\|_{L^{p}})^{q})^{\frac{1}{q}}.$$
The homogeneous Besov space $B^{s}_{p,q}$ is the set of  distribution $u$ in ${\cal S}_{h}^{'}$ such that $\|u\|_{B^{s}_{p,q}}<+\infty$.
\end{definition}
%We are going to extend the definition of the Besov space $B^{s}_{p,r}$ by considering logarithmic index of regularity.
%\begin{definition}
%For
%$s\in\R,\,\,p\in[1,+\infty],\,\,q\in[1,+\infty],\,\,\mbox{and}\,\,u\in{\cal{S}}^{'}(\R^{N})$
%we set:
%$$\|u\|_{B^{s+\log}_{p,q}}=(\sum_{l\in\mathbb{Z}}(|l|2^{ls}\|\D_{l}u\|_{L^{p}})^{q})^{\frac{1}{q}}.$$
%The homogeneous Besov space $B^{s}_{p,q}$ is the set of  distribution $u$ in ${\cal S}_{h}^{'}$ such that $\|u\|_{B^{s}_{p,q}}<+\infty$.
%\label{Besovlog}
%\end{definition}
%\begin{definition}
%For
%$s\in\R,\,\,p\in[1,+\infty],\,\,q\in[1,+\infty],\,\,\mbox{and}\,\,u\in{\cal{S}}^{'}(\R^{N})$
%we set:
%$$\|u\|_{\dot{B}^{s}_{p,q}}=(\sum_{l\in\mathbb{Z}}(2^{ls}\|\dot{\D}_{l}u\|_{L^{p}})^{q})^{\frac{1}{q}}.$$
%The non homogeneous Besov space $B^{s}_{p,q}$ is the set of temperate  distribution $u$  such that $\|u\|_{\dot{B}^{s}_{p,q}}<+\infty$.
%\end{definition}
%In the sequel we will give only properties on the homogeneous Besov spaces but the most of then can be generalize for non homogeneous Besov spaces.
\begin{remarka}The above definition is a natural generalization of the
homogeneous Sobolev and H$\ddot{\mbox{o}}$lder spaces: one can show
that $B^{s}_{\infty,\infty}$ is the homogeneous
H$\ddot{\mbox{o}}$lder space $C^{s}$ and that $B^{s}_{2,2}$ is
the homogeneous space $H^{s}$.
\end{remarka}
\begin{proposition}
\label{derivation,interpolation}
The following properties holds:
\begin{enumerate}
\item there exists a constant universal $C$
such that:\\
$C^{-1}\|u\|_{B^{s}_{p,r}}\leq\|\n u\|_{B^{s-1}_{p,r}}\leq
C\|u\|_{B^{s}_{p,r}}.$
\item If
$p_{1}<p_{2}$ and $r_{1}\leq r_{2}$ then $B^{s}_{p_{1},r_{1}}\hookrightarrow
B^{s-N(1/p_{1}-1/p_{2})}_{p_{2},r_{2}}$.
\item Moreover we have the following interpolation inequalities, it exists $C>0$ such that for any $\theta\in]0,1[$ and $s<\widetilde{s}$ we have:
$$
\begin{aligned}
&\|u\|_{B^{\theta s+(1-\theta)\widetilde{s}}_{p,r}}\leq\|u\|_{B^{s}_{p,r}}^{\theta}\|u\|_{B^{\widetilde{s}}_{p,r}}^{1-\theta},\\
&\|u\|_{B^{\theta s+(1-\theta)\widetilde{s}}_{p,1}}\leq\frac{C}{\theta(1-\theta)(\widetilde{s}-s)}\|u\|_{B^{s}_{p,\infty}}^{\theta}\|u\|_{B^{\widetilde{s}}_{p,\infty}}^{1-\theta}.
\end{aligned}
$$
%$B^{s^{'}}_{p,r_{1}}\hookrightarrow B^{s}_{p,r}$ if $s^{'}> s$ or if $s=s^{'}$ and $r_{1}\leq r$.
 %$(B^{s_{1}}_{p,r},B^{s_{2}}_{p,r})_{\theta,r^{'}}=B^{\theta
%s_{1}+(1-\theta)s_{2}}_{p,r^{'}}$.
\end{enumerate}
\label{interpolation}
\end{proposition}
Let now recall a few product laws in Besov spaces coming directly from the paradifferential calculus of J-M. Bony
(see \cite{BJM,BCD}).
\begin{proposition}
\label{produit}
We have the following laws of product:
\begin{itemize}
\item For all $s\in\R$, $(p,r)\in[1,+\infty]^{2}$ we have:
\begin{equation}
\|uv\|_{B^{s}_{p,r}}\leq
C(\|u\|_{L^{\infty}}\|v\|_{B^{s}_{p,r}}+\|v\|_{L^{\infty}}\|u\|_{B^{s}_{p,r}})\,.
\label{2.2}
\end{equation}
\item Let $(p,p_{1},p_{2},r,\lambda_{1},\lambda_{2})\in[1,+\infty]^{2}$ such that:$\frac{1}{p}\leq\frac{1}{p_{1}}+\frac{1}{p_{2}}$,
$p_{1}\leq\lambda_{2}$, $p_{2}\leq\lambda_{1}$, $\frac{1}{p}\leq\frac{1}{p_{1}}+\frac{1}{\lambda_{1}}$ and
$\frac{1}{p}\leq\frac{1}{p_{2}}+\frac{1}{\lambda_{2}}$. We have then the following inequalities:\\
if $s_{1}+s_{2}+N\inf(0,1-\frac{1}{p_{1}}-\frac{1}{p_{2}})>0$, $s_{1}+\frac{N}{\lambda_{2}}<\frac{N}{p_{1}}$ and
$s_{2}+\frac{N}{\lambda_{1}}<\frac{N}{p_{2}}$ then:
\begin{equation}
\|uv\|_{B^{s_{1}+s_{2}-N(\frac{1}{p_{1}}+\frac{1}{p_{2}}-\frac{1}{p})}_{p,r}}\lesssim\|u\|_{B^{s_{1}}_{p_{1},r}}
\|v\|_{B^{s_{2}}_{p_{2},\infty}},
\label{2.3}
\end{equation}
when $s_{1}+\frac{N}{\lambda_{2}}=\frac{N}{p_{1}}$ (resp $s_{2}+\frac{N}{\lambda_{1}}=\frac{N}{p_{2}}$) we replace
$\|u\|_{B^{s_{1}}_{p_{1},r}}\|v\|_{B^{s_{2}}_{p_{2},\infty}}$ (resp $\|v\|_{B^{s_{2}}_{p_{2},\infty}}$) by
$\|u\|_{B^{s_{1}}_{p_{1},1}}\|v\|_{B^{s_{2}}_{p_{2},r}}$ (resp $\|v\|_{B^{s_{2}}_{p_{2},\infty}\cap L^{\infty}}$),
if $s_{1}+\frac{N}{\lambda_{2}}=\frac{N}{p_{1}}$ and $s_{2}+\frac{N}{\lambda_{1}}=\frac{N}{p_{2}}$ we take $r=1$.
\\
If $s_{1}+s_{2}=0$, $s_{1}\in(\frac{N}{\lambda_{1}}-\frac{N}{p_{2}},\frac{N}{p_{1}}-\frac{N}{\lambda_{2}}]$ and
$\frac{1}{p_{1}}+\frac{1}{p_{2}}\leq 1$ then:
\begin{equation}
\|uv\|_{B^{-N(\frac{1}{p_{1}}+\frac{1}{p_{2}}-\frac{1}{p})}_{p,\infty}}\lesssim\|u\|_{B^{s_{1}}_{p_{1},1}}
\|v\|_{B^{s_{2}}_{p_{2},\infty}}.
\label{2.4}
\end{equation}
If $|s|<\NN$ for $p\geq2$ and $-\frac{N}{p^{'}}<s<\NN$ else, we have:
\begin{equation}
\|uv\|_{B^{s}_{p,r}}\leq C\|u\|_{B^{s}_{p,r}}\|v\|_{B^{\NN}_{p,\infty}\cap L^{\infty}}.
\label{2.5}
\end{equation}
\end{itemize}
\end{proposition}
\begin{remarka}
In the sequel $p$ will be either $p_{1}$ or $p_{2}$ and in this case $\frac{1}{\lambda}=\frac{1}{p_{1}}-\frac{1}{p_{2}}$
if $p_{1}\leq p_{2}$, resp $\frac{1}{\lambda}=\frac{1}{p_{2}}-\frac{1}{p_{1}}$
if $p_{2}\leq p_{1}$.
\end{remarka}
\begin{corollaire}
\label{produit2}
Let $r\in [1,+\infty]$, $1\leq p\leq p_{1}\leq +\infty$ and $s$ such that:
\begin{itemize}
\item $s\in(-\frac{N}{p_{1}},\frac{N}{p_{1}})$ if $\frac{1}{p}+\frac{1}{p_{1}}\leq 1$,
\item $s\in(-\frac{N}{p_{1}}+N(\frac{1}{p}+\frac{1}{p_{1}}-1),\frac{N}{p_{1}})$ if $\frac{1}{p}+\frac{1}{p_{1}}> 1$,
\end{itemize}
then we have if $u\in B^{s}_{p,r}$ and $v\in B^{\frac{N}{p_{1}}}_{p_{1},\infty}\cap L^{\infty}$:
$$\|uv\|_{B^{s}_{p,r}}\leq C\|u\|_{B^{s}_{p,r}}\|v\|_{B^{\frac{N}{p_{1}}}_{p_{1},\infty}\cap L^{\infty}}.$$
\end{corollaire}
%For a proof of this proposition see \cite{5DG}. The limit case $s_{1}+s_{2}=t_{1}+t_{2}=0$ in (\ref{5prodinteressant2}) is of interest.
%When $p\geq2$, the following estimate holds true whenever $s$ is in the range $(-\NN,\NN]$ (see e.g. \cite{5RS}):
%\begin{equation}
%\|uv\|_{B^{-\NN}_{p,\infty}}\leq C\|u\|_{B^{s}_{p,1}}\|v\|_{B^{-s}_{p,\infty}}.
%\end{equation}
The study of non stationary PDE's requires space of type $L^{\rho}(0,T,X)$ for appropriate Banach spaces $X$. In our case, we
expect $X$ to be a Besov space, so that it is natural to localize the equation through Littlewood-Payley decomposition. But, in doing so, we obtain
bounds in spaces which are not type $L^{\rho}(0,T,X)$ (except if $r=p$).
We are now going to
define the spaces of Chemin-Lerner in which we will work, which are
a refinement of the spaces
$L_{T}^{\rho}(B^{s}_{p,r})$.
$\hspace{15cm}$
\begin{definition}
Let $\rho\in[1,+\infty]$, $T\in[1,+\infty]$ and $s_{1}\in\R$. We set:
$$\|u\|_{\widetilde{L}^{\rho}_{T}(B^{s_{1}}_{p,r})}=
\big(\sum_{l\in\mathbb{Z}}2^{lrs_{1}}\|\D_{l}u(t)\|_{L^{\rho}(L^{p})}^{r}\big)^{\frac{1}{r}}\,.$$
We then define the space $\widetilde{L}^{\rho}_{T}(B^{s_{1}}_{p,r})$ as the set of temperate distribution $u$ over
$(0,T)\times\R^{N}$ such that %$\lim_{q\rightarrow+\infty}S_{q}u=0$ in ${\cal S}^{'}((0,T)\times\R^{N})$
%and
$\|u\|_{\widetilde{L}^{\rho}_{T}(B^{s_{1}}_{p,r})}<+\infty$.
\end{definition}
We set $\widetilde{C}_{T}(\widetilde{B}^{s_{1}}_{p,r})=\widetilde{L}^{\infty}_{T}(\widetilde{B}^{s_{1}}_{p,r})\cap
{\cal C}([0,T],B^{s_{1}}_{p,r})$.
Let us emphasize that, according to Minkowski inequality, we have:
\begin{equation}\|u\|_{\widetilde{L}^{\rho}_{T}(B^{s_{1}}_{p,r})}\leq\|u\|_{L^{\rho}_{T}(B^{s_{1}}_{p,r})}\;\;\mbox{if}\;\;r\geq\rho
,\;\;\;\|u\|_{\widetilde{L}^{\rho}_{T}(B^{s_{1}}_{p,r})}\geq\|u\|_{L^{\rho}_{T}(B^{s_{1}}_{p,r})}\;\;\mbox{if}\;\;r\leq\rho
.
\label{Minko}
\end{equation}
\begin{remarka}
It is easy to generalize propositions \ref{produit}, \ref{produit2}
to $\widetilde{L}^{\rho}_{T}(B^{s_{1}}_{p,r})$ spaces. The indices $s_{1}$, $p$, $r$
behave just as in the stationary case whereas the time exponent $\rho$ behaves according to H\"older inequality.
\end{remarka}
In the sequel we will need of composition lemma in $\widetilde{L}^{\rho}_{T}(B^{s}_{p,r})$ spaces (we refer to \cite{BCD} for a proof).
\begin{proposition}
\label{composition}
Let $s>0$, $(p,r)\in[1,+\infty]$ and $u\in \widetilde{L}^{\rho}_{T}(B^{s}_{p,r})\cap L^{\infty}_{T}(L^{\infty})$.
\begin{enumerate}
 \item Let $F\in W_{loc}^{[s]+2,\infty}(\R^{N})$ such that $F(0)=0$. Then $F(u)\in \widetilde{L}^{\rho}_{T}(B^{s}_{p,r})$. More precisely there exists a function $C$ depending only on $s$, $p$, $r$, $N$ and $F$ such that:
$$\|F(u)\|_{\widetilde{L}^{\rho}_{T}(B^{s}_{p,r})}\leq C(\|u\|_{L^{\infty}_{T}(L^{\infty})})\|u\|_{\widetilde{L}^{\rho}_{T}(B^{s}_{p,r})}.$$
\item Let $F\in W_{loc}^{[s]+3,\infty}(\R^{N})$ such that $F(0)=0$. Then $F(u)-F^{'}(0)u\in \widetilde{L}^{\rho}_{T}(B^{s}_{p,r})$. More precisely there exists a function $C$ depending only on $s$, $p$, $r$, $N$ and $F$ such that:
$$\|F(u)-F^{'}(0)u\|_{\widetilde{L}^{\rho}_{T}(B^{s}_{p,r})}\leq C(\|u\|_{L^{\infty}_{T}(L^{\infty})})\|u\|^{2}_{\widetilde{L}^{\rho}_{T}(B^{s}_{p,r})}.$$
\end{enumerate}
\end{proposition}
%Here we recall a result of interpolation which explains the link
%of the space $B^{s}_{p,1}$ with the space $B^{s}_{p,\infty}$, see
%\cite{DFourier}.
%\begin{proposition}
%\label{interpolationlog}
%There exists a constant $C$ such that for all $s\in\R$, $\e>0$ and
%$1\leq p<+\infty$,
%$$\|u\|_{\widetilde{L}_{T}^{\rho}(B^{s}_{p,1})}\leq C\frac{1+\e}{\e}\|u\|_{\widetilde{L}_{T}^{\rho}(B^{s}_{p,\infty})}
%\biggl(1+\log\frac{\|u\|_{\widetilde{L}_{T}^{\rho}(B^{s+\e}_{p,\infty})}}
%{\|u\|_{\widetilde{L}_{T}^{\rho}(B^{s}_{p,\infty})}}\biggl).$$ \label{5Yudov}
%\end{proposition}
Now we give some result on the behavior of the Besov spaces via some pseudodifferential operator (see \cite{BCD}).
\begin{definition}
Let $m\in\R$. A smooth function function $f:\R^{N}\rightarrow\R$ is said to be a ${\cal S}^{m}$ multiplier if for all muti-index $\alpha$, there exists a constant $C_{\alpha}$ such that:
$$\forall\xi\in\R^{N},\;\;|\p^{\alpha}f(\xi)|\leq C_{\alpha}(1+|\xi|)^{m-|\alpha|}.$$
\label{smoothf}
\end{definition}
\begin{proposition}
Let $m\in\R$ and $f$ be a ${\cal S}^{m}$ multiplier. Then for all $s\in\R$ and $1\leq p,r\leq+\infty$ the operator $f(D)$ is continuous from $B^{s}_{p,r}$ to $B^{s-m}_{p,r}$.
\label{singuliere}
\end{proposition}
%Actually, in \cite{BCD}, the proposition below is proved for non-homogeneous Besov spaces. The adaptation to homogeneous spaces is straightforward.
Let us now give some estimates for the heat equation:
\begin{proposition}
\label{chaleur} Let $s\in\R$, $(p,r)\in[1,+\infty]^{2}$ and
$1\leq\rho_{2}\leq\rho_{1}\leq+\infty$. Assume that $u_{0}\in B^{s}_{p,r}$ and $f\in\widetilde{L}^{\rho_{2}}_{T}
(B^{s-2+2/\rho_{2}}_{p,r})$.
Let u be a solution of:
$$
\begin{cases}
\begin{aligned}
&\p_{t}u-\mu\D u=f\\
&u_{t=0}=u_{0}\,.
\end{aligned}
\end{cases}
$$
Then there exists $C>0$ depending only on $N,\mu,\rho_{1}$ and
$\rho_{2}$ such that:
$$\|u\|_{\widetilde{L}^{\rho_{1}}_{T}(\widetilde{B}^{s+2/\rho_{1}}_{p,r})}\leq C\big(
 \|u_{0}\|_{B^{s}_{p,r}}+\mu^{\frac{1}{\rho_{2}}-1}\|f\|_{\widetilde{L}^{\rho_{2}}_{T}
 (B^{s-2+2/\rho_{2}}_{p,r})}\big)\,.$$
 If in addition $r$ is finite then $u$ belongs to $C([0,T],B^{s}_{p,r})$.
\end{proposition}
\subsection*{Hybrid Besov spaces}
The homogeneous Besov spaces fail to have nice inclusion properties: owing to the low frequencies, the embedding $B^{s}_{p,1}\hookrightarrow B^{t}_{p,1}$ does not hold for $s>t$. Still, the functions of $B^{s}_{p,1}$ are locally more regular than those of $B^{t}_{p,1}$: for any $\phi\in C^{\infty}_{0}$ and $u\in B^{s}_{p,1}$, the function $\phi u\in B^{t}_{p,1}$. This motivates the definition of Hybrid Besov spaces introduced by R. Danchin  (see a definition in \cite{BCD}) where the growth conditions satisfied by the dyadic blocks and the coefficient of integrability are not the same for low and high frequencies. Hybrid Besov spaces have been used by R. Danchin in order to prove global well-posedness for compressible gases in critical spaces (we refer to \cite{BCD} for an elegant proof of this result). We generalize here a little bit the definition by allowing for different Lebesgue norms in low and high frequencies.
\begin{definition}
 \label{def1.9}
Let  $l_{0}\in\mathbb{N}$, $s,t,\in\R$, $(r,r_{1})\in [1,+\infty]^{2}$ and $(p,q)\in[1,+\infty]$. We set:
$$\|u\|_{\widetilde{B}^{s,t}_{p,q,1}}=\sum_{l\leq l_{0}}2^{ls}\|\D_{l}u\|_{L^{p}}+\sum_{l> l_{0}}2^{lt}\|\D_{l}u\|_{L^{q}},$$
and:
$$\|u\|_{\widetilde{B}^{s,t}_{(p,r),(q,r_{1})}}=\big(\sum_{l\leq l_{0}}(2^{ls}\|\D_{l}u\|_{L^{p}})^{r}\big)^{\frac{1}{r}}+\big(\sum_{l> l_{0}}(2^{lt}\|\D_{l}u\|_{L^{q}})^{r_{1}}\big)^{\frac{1}{r_{1}}}.$$
\end{definition}
\begin{remarka}
When $p=q$ and $r=r_{1}$ we will note to simplify $\widetilde{B}^{s,t}_{(p,r),(p,r)}=\widetilde{B}^{s,t}_{p,r}$.
%It will be important in the sequel to chose $l_{0}$ big enough.
\end{remarka}
\begin{notation}
 We will often use the following notation:
$$u_{BF}=\sum_{l\leq l_{0}}\D_{l}u\;\;\;\mbox{and}\;\;\;u_{HF}=\sum_{l> l_{0}}\D_{l}u.$$
\end{notation}
\begin{remarka}
 We have the following properties:
\begin{itemize}
 \item We have $\widetilde{B}^{s,s}_{p,p,1}=B^{s}_{p,1}$.
\item If $s_{1}\geq s_{3}$ and $s_{2}\geq s_{4}$ then $\widetilde{B}^{s_{3},s_{2}}_{p,q,1}\h \widetilde{B}^{s_{1},s_{4}}_{p,q,1}$.
\end{itemize}
\label{r13}
\end{remarka}
\begin{remarka}
In the sequel we shall often use this hybrid Besov space in order to distinguish the behavior of our solution in low and high frequencies, in particular we would like to mention that we can prove results analogous to propositions \ref{produit} and corollary \ref{produit2}.
\end{remarka}
We shall conclude this section by some example of initial data verifying the theorem \ref{ftheo1} with large energy initial data when $N=2$.More generally we are interested in defining initial data which are small in $\widetilde{B}^{\N-1,\N}_{2,\infty}$ but large in $B^{\N}_{2,2}$. It implies in particular that your initial data $\n\sqrt{\rho_{0}}$ have large energy data when $N=2$, furthermore $\ln\rho_{0}$ is large in $B^{\N}_{2,2}$ which is a scaling invariant space (it improves in particular the results of \cite{fDD} where the initial density is assumed small in $B^{\N}_{2,1}$).\\
Let us start with recalling a classical example of function in $B^{s}_{p,\infty}$, by sake of completness we are going to recall the proof (see also \cite{BCD}). 
\begin{proposition}
\label{homog}
Let $\sigma\in]0,N[$. For any $p\in[1,+\infty]$, the function $|\cdot|^{-\sigma}$ belongs to $B^{\NN-\sigma}_{p,\infty}$.
\end{proposition}
{\bf Proof:} By proposition \ref{interpolation} it suffices to show that $u_{\sigma}=|\cdot|^{-\sigma}$ belongs to $B^{N-\sigma}_{1,\infty}$. Let us introduce a smooth compactly supported function $\chi$ which is identically equal to $1$ near the unit ball and such that $u$ is splitting as follows:
$$u_{\sigma}=u_{0}+u_{1}\;\;\;\mbox{with}\;\;u_{0}(x)=\chi(x)|x|^{-\sigma}\;\;\mbox{and}\;\;u_{1}(x)=(1-\chi(x))|x|^{-\sigma}.$$
Clearly $u_{0}$ is in $L^{1}$ and $u_{1}$ belongs in $L^{q}$ whenever $q>\frac{N}{\sigma}$. The homogeneity of the function $u_{\sigma}$  gives via a change of variable:
\begin{equation}
\begin{aligned}
\D_{j}u_{\sigma}&=2^{jN}u_{\sigma}*h(2^{j}\cdot)\\
&=2^{j(N+\sigma)}u_{\sigma}(2^{j}\cdot)*h(2^{j}\cdot)\\
&=2^{j\sigma}(\D_{0}u_{\sigma})(2^{j}\cdot).
\end{aligned}
\label{invariance}
\end{equation}
Therefore, $2^{j(N-\sigma)}\|\D_{j}u_{\sigma}\|_{L^{1}}=\|\D_{0}u_{\sigma}\|_{L^{1}}$, it reduces the problem to show that $\D_{0}u_{\sigma}$ is in $L^{1}$. As $u_{0}$ is in $L^{1}$, $\D_{0}u_{0}$ is also in $L^{1}$ according the continuity of the operator $\D_{0}$ on Lebesgue spaces. By Bernstein inequalities, we have:
$$\|\D_{0}u_{1}\|_{L^{1}}\leq C_{k}\|D^{k}\D_{0}u_{1}\|_{L^{1}}\leq C_{k}\|D^{k}u_{1}\|_{L^{1}}.$$
Leibniz's formula ensures that $D^{k}u_{1}-(1-\chi)D^{k}u_{\sigma}$ is a smooth compactly supported function. We then complete the proof by choosing $k$ such that $k>N-\sigma$. \hfill {$\Box$}\\
\\
We can now deduce from the previous proposition suitable functions verifying the theorem \ref{ftheo1}.
\begin{corollaire}
Let us consider:
$$u_{0,\e,l_{0}}(x)={\cal S}(\frac{1}{|x|^{1-\e}}),$$
with ${\cal S}$ defined by ${\cal S}\hat{u}(\xi)=1_{\{ \xi\in \R^{N}\sigma B(0,2^{l_{0}})}(\xi)\hat{u}({\xi})$.
Then for all $r\in[1,+\infty[$, for all $M>$ for all $\e_{1}>0$ it exits $\e>0$, it exits $l_{0}>0$ such that:
\begin{equation}
\begin{aligned}
&\|u_{0}\|_{B^{\N-1}_{2,\infty}}\leq\e_{1},\\
&\|u_{0}\|_{B^{\N-1}_{2,r}}\geq M.
\end{aligned}
\end{equation}
\label{corotech}
\end{corollaire}
{\bf Proof:} In the sequel in order to simplify the notation we shall write $u_{0}$ for $u_{0,\e,l_{0}}$. Let us denote by $u_{\sigma}$ the function $\frac{1}{|x|^{\sigma}}$ with $\sigma\in]0,N[$. By (\ref{invariance}) we observe that:
\begin{equation}
\|\D_{l}u_{\sigma}\|_{L^{2}}=2^{l(\sigma-\N)}\|\D_{0}u_{\sigma}\|_{L^{2}}.
\label{techimp}
\end{equation}
It implies that for $l\geq l_{0}$ we have for $M_{0}$ independent on $\e$ when $\e\leq\frac{1}{10}$:
$$2^{l(\N-1+\e)}\|\D_{l}u_{0}\|_{L^{2}}=\|\D_{0}u_{1-\e}\|_{L^{2}}\leq M_{0}.$$
Indeed using the same arguments than in the proof of the proposition \ref{homog} we have by Bernstein inequality for $C>0$:
$$
\begin{aligned}
\|\D_{0}u_{1-\e}\|_{L^{2}}&\leq C \|\D_{0}u_{1-\e}\|_{L^{1}},\\
&\leq C (\|\chi u_{1-\e}\|_{L^{1}}+C^{k}\|D^{k}((1-\chi)u_{1-\e})\|_{L^{1}}),\\
&\leq M_{0}.
\end{aligned}
$$
%In particular we are able to chosen initial data providing global strong solution and large in the energy space. Indeed assume that:
%$$\forall l\geq l_{0},\;2^{l(\N-1+\e}\|\D_{l}u\|_{L^{2}}=M$$
%We set now:
%$$M2^{-l_{0}\e}=\e_{1},$$
According to (\ref{techimp}) we deduce that for all $l\geq l_{0}$:
$$2^{l(\N-1)}\|\D_{l}u_{0}\|_{L^{2}}=2^{-l\e}\|\D_{0}u_{\sigma}\|_{L^{2}}\leq M_{0} 2^{-l_{0}\e}.$$
It implies in particular since $\D_{l}u_{0}=0$ for $l<l_{0}$ that:
\begin{equation}
\|u_{0}\|_{B^{\N-1}_{2,\infty}}\leq M_{0} 2^{-l_{0}\e}.
\label{infiny}
\end{equation}
Let us now estimate the norm of $u_{0}$ in $B^{\N-1}_{2,r}$ with $r\in[1,+\infty[$:
\begin{equation}
\begin{aligned}
\|u_{0}\|_{B^{\N-1}_{2,r}}&=(\sum_{l\in\mathbb{Z}}2^{rl(\N-1)}\|\D_{l}u_{0}\|^{r}_{L^{2}})^{\frac{1}{r}},\\
&=(\sum_{l\geq l_{0}}2^{-rl\e}M^{r}_{\e})^{\frac{1}{r}},\\      
&= M_{\e} 2^{-l_{0}\e}(\frac{1}{1-2^{-r\e}})^{\frac{1}{r}},
\end{aligned}
\label{rbesov}
\end{equation}
It implies that for all $M>0$, for all $r>0$ it exits $\e>0$ small enough such that:
\begin{equation}
(\frac{1}{1-2^{-r\e}})^{\frac{1}{r}}\geq M.
\label{hyp2}
\end{equation}
We fix now $l_{0}$ such that:
\begin{equation}
M_{0} 2^{-l_{0}\e}=\e_{1}.
\label{hyp1}
\end{equation}
Let us prove now that when $\e\leq \frac{1}{10}$ it exist $\alpha>0$ such that $M_{\e}\geq\alpha>0$. Assume by the absurd that this is wrong. It implies that it exists a sequel 
$(\e_{n})_{n\in\mathbb{N}}$ which converges to $0$ when $n$ goes to infinity such that $M_{\e_{n}}=\|\D_{0}u_{1-\e_{n}}\|_{L^{2}}$ goes to $0$ when $n$ goes to infinity. By Plancherel theorem and the fact that we know the Fourier transform of $|\cdot|^{-1+\e}$ (see \cite{BCD} p 23) it implies that:
$$\|{\cal F}\D_{0}u_{1-\e_{n}}\|_{L^{2}}=c_{N,1-\e_{n}}\|\varphi |\cdot|^{1-\e_{n}-N}\|_{L^{2}}\rightarrow _{n\rightarrow +\infty}0.$$
Since we can bound by below $\|\varphi |\cdot|^{1-\e_{n}-N}\|_{L^{2}}$ independently of $n$, it implies that $c_{N,1-\e_{n}}$ goes to $0$ when $n$ goes to infinity and this is absurd.\\
We obtain finally from (\ref{infiny}), (\ref{hyp1}) and (\ref{rbesov}), (\ref{hyp2}) that:
\begin{equation}
\begin{aligned}
&\|u_{0}\|_{B^{\N-1}_{2,\infty}}\leq\e_{1},\\
&\|u_{0}\|_{B^{\N-1}_{2,\infty}}\geq \frac{M\e_{1}\alpha}{M_{0}}.
\end{aligned}
\label{conclu}
\end{equation}
It concludes the proof of the corollary.
%Assume that $l_{0}$ is sufficiently large such that:
\hfill {$\Box$}
%which is large when $\e$ goes to $0$. It proves in particular that $\|u_{0}\|_{B^{\N-1}_{2,1}}$ is large, it consists in choosing a $u_{0}$ which is homogeneous in high frequencies of deter $-1+\e$ with $\e$ very small. For example we can choose:
%$$u_{0}(x)={\cal S}\frac{1}{|x|^{1-\e}},$$
%with ${\cal S}\hat{u}=1_{\xi\in B(0,2^{l_{0}})^{c}}\hat{u}(\xi)$.
\section{Proof of theorems \ref{ftheo1} and \ref{ftheo2}}
\label{section4}
In this part we are interested in proving the theorems \ref{ftheo1} and \ref{ftheo2} of
existence of strong solutions in critical space for the scaling of the equations. We want to point out that in the theorem \ref{ftheo1} the viscosity and the capillarity coefficients are chosen with a very specific structure. This fact will be crucial in the sequel of the proof in order to obtain estimates on the density without
assuming any control on the vacuum or on the $L^{\infty}$ norm of the density for the theorem \ref{ftheo1}. Indeed when $\mu(\rho)=\mu\rho$ and $\kappa(\rho)=\frac{\kappa}{\rho}$ with $\kappa>0$ then the system depends only on the unknown $\ln\rho$ in a linear way and provides a new entropy (see \cite{Hprepa} for more details).\\
In a second time in order to prove the second part of the theorem \ref{ftheo1}, we are going to estimate the vacuum (which means the $L^{\infty}$ norm of $\frac{1}{\rho}$) and the $L^{\infty}$ norm of the density. More precisely we are going to write a solution $(q,u)$ under the form:
$$(q,u)=(q_{L},u_{L})+(\bar{q},\bar{u}),$$
with $(q_{L},u_{L})$ the solution of the linearized part of the system (\ref{NHV1}). We shall combine maximum principle arguments on $q_{L}$ in order to bound $q_{L}$ in $L^{\infty}$ norm and regularizing effect on the third index of Besov space for $\bar{q}$. Indeed we will prove that  $\bar{q}$ is in $\widetilde{L}_{T}^{\infty}(B^{\N}_{2,1})$ for any $T>0$ which is embedded in $L^{\infty}_{T}(L^{\infty})$. Let us mention that in order to prove that $q_{L}$ is bounded in $L^{\infty}$ norm we shall prove a very accurate characterization of the Besov space in term of the semi group associated to the linearized part of the system (\ref{NHV1}) (see the proposition \ref{propmaxfinal}). For more details on this part which is the main difficulty of the proof we refer to the subsection \ref{submaximum} and \ref{controldens}.\\
%(we shall precise more in details this approach in the sequel of this section see the subsection \ref{}).\\
As a first step of the proof of theorems \ref{ftheo1} and \ref{ftheo2}, let us start with studying the
linear part of the system (\ref{NHV1}) (when we are interested in proving the existence of strong solution in finite time with large initial data, in particular it implies that we do not consider the low frequencies terms coming from the pressure) which corresponds to the following system:
$$
\begin{cases}
\p_{t}q+{\rm div}u=F,\\
\p_{t}u-a\D u-b\n{\rm
div}u-c\n\D q=G,
%\p_{t}{\cal T}-{\rm div}(d\n{\cal T})=H,\\
\end{cases}
\leqno{(N)}
$$
\subsection{Study of the linearized equation}
We want to prove a priori estimates in Chemin-Lerner spaces for system $(N)$
with the following hypotheses on $a,b,c$ which are constant:
$$
\begin{aligned}
&0<a<\infty,\;0<
a+b<\infty\;\;\mbox{and}\;\;\;0< c<\infty.
%&0<c_{4}\leq d<M_{4}<\infty.\\
\end{aligned}
$$
This system has been studied by Danchin and Desjardins (see \cite{fDD}) in the framework of the Besov space $B^{s}_{2,1}$, the
following proposition uses exactly the same type of arguments used in \cite{fDD} excepted that we extend the result to general Besov spaces $B^{s}_{2,r}$ with $r\in[1,+\infty]$. By sake of completness we are going to show the following proof. 
\begin{proposition}
\label{flinear3}  Let $1\leq  r\leq+\infty$, $s\in\R$, and we assume that $(q_{0},u_{0})$ belongs in $B_{2,r}^{\N+s}\times (B_{2,r}^{\N-1+s})^{N}$ with the source terms
$(F,G)$ in $\widetilde{L}^{1}_{T}(B^{\N+s}_{2,r})\times
(\widetilde{L}^{1}_{T}(B^{\N-1+s}_{2,r}))^{N}$.\\
Let $(q,u)\in(\widetilde{L}^{1}_{T}(B^{\N+s+2}_{2,r})\cap\widetilde{L}^{\infty}_{T}(B^{\N+s}_{2,r}))\times\big(\widetilde{L}^{1}_{T}(B^{\N+s+1}_{2,r})
\cap \widetilde{L}^{\infty}_{T}(B^{\N+s-1}_{2,r})\big)^{N}$ be a solution of the
system $(N)$, then there exists a universal constant $C$ such that for any $T>0$:
$$
\begin{aligned}
&\|(\n q,u)\|_{\widetilde{L}^{1}_{T}(B^{\N+1+s}_{2,r})\cap \widetilde{L}^{\infty}_{T}(B^{\N-1+s}_{2,r})}\leq C(\|(\n q_{0},u_{0})\|_{B^{\N-1+s}_{2,2}}+\|(\n
F,G)\|_{ \widetilde{L}^{1}_{T}(B^{\N-1+s}_{2,r})}).
%&\hspace{5cm}+\|(\n q,u)\|_{\widetilde{L}^{2}_{T}(B^{\N+s}_{2,r})}
%(\|\n a\|_{\widetilde{L}^{2}(B^{\N}_{2,r})}+\|\n b\|_{\widetilde{L}^{2}(B^{\N}_{2,r})}).
\end{aligned}
$$
\end{proposition}
%L'important est de savoir comment on peut generaliser ce systeme.
{\bf Proof:} As we mentioned previously, we are going to follow the arguments developed in \cite{fDD}. Let us apply to the system $(N)$ the operator $\D_{l}$ which gives:
\begin{eqnarray}
&&\p_{t}q_{l}+{\rm div}u_{l}=F_{l}\label{fl1}\\
&&\p_{t}u_{l}-{\rm div}(a\n u_{l})-\n(b\,{\rm div}u_{l})-c\n \D
q_{l}=G_{l}\label{fl2}
\end{eqnarray}
Performing integrations by parts and using (\ref{fl1}) we have:
$$
\begin{aligned}
-c\int_{\R^{N}}u_{l}\cdot\n\D
q_{l}dx&=c\int_{\R^{N}}{\rm div}u_{l}\,\D
q_{l}dx,\\
&=-c\int_{\R^{N}}\p_{t}q_{l}\,\D
q_{l}dx+c\int_{\R^{N}}F_{l}\,\D
q_{l}dx,\\
&c=\int_{\R^{N}}\p_{t}\n q_{l}\cdot \n
q_{l}dx-c\int_{\R^{N}}\n F_{l}\cdot \n
q_{l}dx,\\
&=\frac{c}{2}\frac{d}{dt}\int_{\R^{N}}|\n
q_{l}|^{2}dx-c\int_{\R^{N}}\n
q_{l}.\n F_{l}\,dx.
\end{aligned}
$$
Next, we take the inner product of (\ref{fl2}) with $u_{l}$
and using the previous equality, it yields:
\begin{equation}
\begin{aligned}
&\frac{1}{2}\frac{d}{dt}\big(\|u_{l}\|_{L^{2}}^{2}+c\int_{\R^{N}}|\n
q_{l}|^{2}dx\big)+\g(a|\n u_{l}|^{2}+b|{\rm
div}u_{l}|^{2})dx\\
&\hspace{7cm}=\g G_{l}.u_{l}\,dx+c\g\n q_{l}.\n F_{l}\,dx\,.\\
\label{f3}
\end{aligned}
\end{equation}
In  order to recover some terms in $\D q_{l}$ we take the
inner product of the gradient of (\ref{fl1}) with $u_{l}$, the inner
product
scalar of (\ref{fl2}) with $\n q_{l}$ and we sum, we obtain then:
\begin{equation}
\begin{aligned}
\frac{d}{dt}\g\n q_{l}.u_{l}dx+c\g (\D q_{l})^{2}dx
=&\g(G_{l}.\n q_{l}+|{\rm div}u_{l}|^{2}+u_{l}.\n F_{l}\\
&\hspace{2cm}-a\n u_{l}:\n ^{2}q_{l}-b\D q_{l}{\rm div}u
_{l})dx.
\label{f4}
\end{aligned}
\end{equation}
Let $\alpha>0$ small enough. We define $k_{l}$ by:
\begin{equation}
k_{l}^{2}=\|u_{l}\|_{L^{2}}^{2}+
c\|\n\q\|^{2}_{L^{2}}+2\alpha\g \n\q.\ui dx\;.\label{f5.39}
\end{equation}
By using (\ref{f3}), (\ref{f4}) and the Young inequalities, we have
by summing and the fact that $\alpha$ is chosen small enough:% and in taking $\alpha$ small:
\begin{equation}
\begin{aligned}
&\frac{1}{2}\frac{d}{dt}k_{l}^{2}+\frac{1}{2}\g(a|\n
u_{l}|^{2}+b|{\rm div}u_{l}|^{2}+2\alpha c|\D\q|^{2})dx
\lesssim\|G_{l}\|_{L^{2}}(2\alpha\|\n\q\|_{L^{2}}+\|\ui\|_{L^{2}})\\
&\hspace{7cm}+\|\n F_{l}\|_{L^{2}}(2\alpha\|\ui\|_{L^{2}}+\|\n\q\|_{L^{2}}).
\label{fl5}
\end{aligned}
\end{equation}
For small enough $\alpha$ and applying the Young inequality, we have according (\ref{f5.39}):
\begin{equation}
\frac{1}{2}k_{l}^{2}\leq\|u_{l}\|^{2}+
c\|\n\q\|^{2}_{L^{2}}\leq\frac{3}{2}k_{l}^{2}\;. \label{fl6}
\end{equation}
Hence according to (\ref{fl5}) and (\ref{fl6}) there exists $K>0$ small enough, $C>0$ such that:
$$
\begin{aligned}
\frac{1}{2}\frac{d}{dt}k_{l}^{2}+K2^{2l}k_{l}^{2} \leq&\,\,
C\,k_{l}\,(\|G_{l}\|_{L^{2}}+\|\n
F_{l}\|_{L^{2}}).\\
\end{aligned}
$$
By integrating with respect to the time, we obtain:
$$
\begin{aligned}
k_{l}(t)\leq &\,e^{-K2^{2l}t}k_{l}(0)+C\int_{0}^{t}
e^{-K2^{2l}(t-\tau)}(\|\n
F_{l}(\tau)\|_{L^{2}}+\|G_{l}(\tau)\|_{L^{2}})d\tau\;.
\end{aligned}
$$
Using convolution inequalities we get for $1\leq\rho_{1}\leq \rho\leq+\infty$:
\begin{equation}
\begin{aligned}
\|k_{l}\|_{L^{\rho}([0,T])}\lesssim&\,\big(2^{-\frac{2l}{\rho}}k_{l}(0)+2^{-2l(1+\frac{1}{\rho}-\frac{1}{\rho_{1}})}\|(\n
F_{l},G_{l})\|_{L^{\rho_{1}}_{T}(L^{2})}\big).
\label{f8}
\end{aligned}
\end{equation}
Moreover since we have:
$$C^{-1}\,k_{l}\leq\|\n\q\|_{L^{2}}+\|\ui\|_{L^{2}}\leq C\,k_{l},$$
multiplying by $2^{(\N-1+s+\frac{2}{\rho})l}$, taking the $l^{r}$ norm and using (\ref{fl6}), we end up with:
$$
\begin{aligned}
\|(\n q,u)&\|_{L^{\rho}_{T}(B^{\N-1+s+\frac{2}{\rho}}_{2,r})}\leq\,\|(\n
F,G)\|_{\widetilde{L}^{\rho_{1}}_{T}(B^{\N-3+s+\frac{2}{\rho_{1}}}_{2,r})}
+\|(\n q_{0},u_{0})\|_{B^{\N-1+s}_{2,r}}.
\label{f27}
\end{aligned}
$$
It conclude the proof of the proposition.
\hfill {$\Box$}\\
\\
Let us extend the result of the proposition \ref{flinear3} to the case where we include the pressure term inside of the linearized system. This is necessary when we are interested in dealing with the existence of global strong solution with small initial data. Indeed in this case it is very important to take into account the behavior in low frequencies of the density and in particular the pressure term which provides the behavior of the density in low frequencies. More precisely we will consider the following linear system:
$$
\begin{cases}
\p_{t}q+{\rm div}u=F,\\
\p_{t}u-a\D u-b\n{\rm
div}u-c\n\D q+d\n q=G,
%\p_{t}{\cal T}-{\rm div}(d\n{\cal T})=H,\\
\end{cases}
\leqno{(N1)}
$$
We now want to prove a priori estimates in Chemin-Lerner spaces for system $(N1)$
with the following hypotheses on $a,b,c,d$ which are constant:
$$
\begin{aligned}
&0< a<\infty,\;0<
a+b<\infty,\;0< c<\infty\;\;\mbox{and}\;\;0< d<\infty.
\end{aligned}
$$
This system has also been studied by Danchin and Desjardins in \cite{fDD} in the framework of Besov spaces of the type $B^{s}_{2,1}$ with $s\in\R$, the
following proposition extends  the results of \cite{fDD} to the case of general Besov spaces $B^{s}_{2,r}$ with $r\in[1,+\infty]$ by using similar arguments.
\begin{proposition}
\label{1flinear3}  
Let $1\leq  r\leq+\infty$, $s\in\R$, and we assume that $(q_{0},u_{0})$ belongs in $\widetilde{B}_{2,r}^{\N-1+s,\N+s}\times (B_{2,r}^{\N-1+s})^{N}$. Furthermore we suppose that the source terms
$(F,G)$ are in $\widetilde{L}^{1}_{T}(\widetilde{B}^{\N-1+s,\N+s}_{2,r})\times
(\widetilde{L}^{1}_{T}(B^{\N-1+s}_{2,r}))^{N}$.\\
Let $(q,u)\in(\widetilde{L}^{1}_{T}(\widetilde{B}^{\N+s+1,\N+s+2}_{2,r})\cap\widetilde{L}^{\infty}_{T}(\widetilde{B}^{\N-1+s,\N+s}_{2,r}))\times\big(\widetilde{L}^{1}_{T}(B^{\N+s+1}_{2,r})
\cap \widetilde{L}^{\infty}_{T}(B^{\N+s-1}_{2,r})\big)^{N}$ be a solution of the
system $(N1)$, then there exists a universal constant $C$ such that for any $T>0$:
$$
\begin{aligned}
&\|q\|_{\widetilde{L}^{1}_{T}(\widetilde{B}^{\N+1+s,\N+2+s}_{2,r})}+\|q\|_{\widetilde{L}^{\infty}_{T}(B^{\N-1+s,\N+s}_{2,r})}
+\|u\|_{\widetilde{L}^{1}_{T}(B^{\N+1+s}_{2,r})}+ \|u\|_{\widetilde{L}^{\infty}_{T}(B^{\N-1+s}_{2,r})}\\
&\hspace{1cm}\leq C\big(\|q_{0}\|_{\widetilde{B}^{\N-1+s,\N+s}_{2,r}}+\|u_{0}\|_{B^{\N+s}_{2,r}}+\|F\|_{\widetilde{L}^{1}_{T}(\widetilde{B}^{\N-1+s,\N+s}_{2,r})}
+\|G\|_{\widetilde{L}^{1}_{T}(B^{\N-1+s}_{2,r})}\big).
%&\hspace{5cm}+\|(\n q,u)\|_{\widetilde{L}^{2}_{T}(B^{\N+s}_{2,r})}
%(\|\n a\|_{\widetilde{L}^{2}(B^{\N}_{2,r})}+\|\n b\|_{\widetilde{L}^{2}(B^{\N}_{2,r})}).
\end{aligned}
$$
\end{proposition}
%L'important est de savoir comment on peut generaliser ce systeme.
{\bf Proof:} It suffices to follow exactly the same lines as the proof of proposition \ref{flinear3} except that we have to consider the following $k_{l}$:
$$k_{l}^{2}=\|u_{l}\|_{L^{2}}^{2}+c\|\n q_{l}\|_{L^{2}}^{2}+d\|q_{l}\|_{L^{2}}^{2}+2\alpha\int_{\R^{N}}\n q_{l}\cdot u_{l}dx.$$
Now choosing $\alpha$ suitably small , we deduce that:
\begin{equation}
\frac{1}{2}k_{l}^{2}\leq \|u_{l}\|_{L^{2}}^{2}+c\|\n q_{l}\|_{L^{2}}^{2}+d\|q_{l}\|_{L^{2}}^{2}\leq \frac{3}{2}k_{l}^{2}.
\label{lienint}
\end{equation}
By combining energy estimates in frequencies space, we show as in \cite{fDD} that it exists $c,C>0$ such that:
$$\frac{1}{2}\frac{d}{dt}k_{l}^{2}+c2^{2l}k_{l}^{2}\leq Ck_{l}(\|G_{l}\|_{L^{2}}+\|(\n F_{l},F_{l})\|_{L^{2}}).$$
As in the proof of proposition \ref{flinear3} routine computations yield proposition \ref{1flinear3}.
\hfill {$\Box$}\\
\\
We are now interested in studying the following system with $\mu>0$ and $\kappa>0$:
\begin{equation}
\begin{cases}
\begin{aligned}
&\p_{t}h_{2}+{\rm div}u_{2}-\mu\n\ln\rho_{1}\cdot\n h_{2}+u_{2}\cdot\n\ln\rho_{1}=F,\\[2mm]
&\p_{t}u_{2}-\mu\D u_{2}-\mu\n{\rm div}u_{2}-\kappa\n\D h_{2}+K\n h_{2}-2\mu\n \ln\rho_{1}\cdot D u_{2}-2\mu\n h_{2}\cdot D u_{1}\\
&\hspace{5cm}+u_{1}\cdot\n u_{2}+u_{2}\cdot\n u_{1}-\mu^{2}\n(\n\ln\rho_{1}\cdot\n h_{2})=G,\\
&(h_{2}(0,\cdot),u_{2}(0,\cdot))=(h_{0}^{2},u_{0}^{2}).
\end{aligned}
\end{cases}
\label{bafond3systeme3}
\end{equation}
Here $(h_{2},u_{2})$ are the unknowns, $\rho_{1}$ and $u_{1}$ corresponds to some functions defined in suitable Chemin Lerner Besov space $\widetilde{L}^{\rho}(B^{s}_{p,r})$ that we will precise below in the proposition \ref{fcrucialprop} and $F$, $G$ are source term (we will also precise their regularity).  In order to prove the theorem \ref{ftheo3} we will need precise estimates on the solution $(h_{2},u_{2})$ of (\ref{bafond3systeme3}) in terms of Chemin Lerner Besov spaces, to do this we are going to prove the following proposition.
\begin{proposition}
Let $(h^{2}_{0},u^{2}_{0})\in \widetilde{B}^{\N-1,\N}_{2,1}\times B^{\N-1}_{2,1}$ and we assume that $\ln\rho_{1}$ belongs in $\widetilde{L}^{\infty}(\R^{+},\widetilde{B}^{\N-1,\N}_{2,1})\cap \widetilde{L}^{1}(\R^{+},\widetilde{B}^{\N+1,\N+2}_{2,1})$ and $u_{1}$ in $\widetilde{L}^{\infty}(\R^{+},B^{\N-1}_{2,1})\cap \widetilde{L}^{1}(\R^{+},B^{\N+1}_{2,1})$.\\
Furthermore $(F,G)$ are in $ \widetilde{L}^{1}(\R^{+},\widetilde{B}^{\N-1,\N}_{2,1})\times  \widetilde{L}^{1}(\R^{+},B^{\N-1}_{2,1}).$ Let $(h_{2},u_{2})$ the solution of the linear system 
(\ref{bafond3systeme3}) then it exists $C>0$ such that $(h_{2},u_{2})$ verify for any $T>0$:
\begin{equation}
\begin{aligned}
&\|h_{2}\|_{\widetilde{L}^{1}_{T}(\widetilde{B}^{\N+1,\N+2}_{2,1})}+\|h_{2}\|_{\widetilde{L}^{\infty}_{T}(B^{\N-1,\N}_{2,1})}
+\|u_{2}\|_{\widetilde{L}^{1}_{T}(B^{\N+1}_{2,1})}+ \|u_{2}\|_{\widetilde{L}^{\infty}_{T}(B^{\N-1}_{2,1})}\\
&\leq C\big(\|h^{2}_{0}\|_{\widetilde{B}^{\N-1,\N}_{2,1}}+\|u^{2}_{0}\|_{B^{\N}_{2,1}}+\|F\|_{\widetilde{L}^{1}_{T}(\widetilde{B}^{\N-1,\N}_{2,1})}
+\|G\|_{\widetilde{L}^{1}_{T}(B^{\N-1}_{2,1})}\big)\\
&\times \exp\biggl(C\int_{0}^{T}\big( \|u_{1}\|_{B^{\N-\frac{1}{2}}_{2,\infty}}^{4}+\|u_{1}\|^{\frac{4}{3}}_{B^{\N+\frac{1}{2}}_{2,\infty}}+\|\n\ln\rho_{1}\|^{4}_{B^{\N-\frac{1}{2}}_{2,\infty}}+\|\n\ln\rho_{1}\|^{\frac{4}{3}}_{B^{\N+\frac{1}{2}}_{2,\infty}}\big)(s)ds\biggl).\\
&
\end{aligned}
\label{raimportant1}
\end{equation}
 \label{fcrucialprop}
 \end{proposition}
 \begin{remarka}
Here to simplify we consider directly $\R^{+}$ than $(0,T)$ for any $T>0$, we would like to point out that in the proposition \ref{fcrucialprop} the main estimates (\ref{raimportant1}) involves a control of $\n\ln\rho_{1}$ in $L^{\frac{4}{3}}(\R^{+},B^{\N+\frac{1}{2}}_{2,\infty})\cap L^{4}(\R^{+},B^{\N-\frac{1}{2}}_{2,\infty})$ and of $u_{1}$ in $L^{\frac{4}{3}}(\R^{+},B^{\N+\frac{1}{2}}_{2,\infty})\cap L^{4}(\R^{+},B^{\N-\frac{1}{2}}_{2,\infty})$. By Minkowski inequality ( see (\ref{Minko})) and since $l^{p}$ is embedded in $l^{\infty}$ for any $p\geq 1$, we know that for example:
 \begin{equation}
 \begin{aligned}
 \|\n\ln\rho_{1}\|_{L^{\frac{4}{3}}(\R^{+},B^{\N+\frac{1}{2}}_{2,\infty})\cap L^{4}(\R^{+},B^{\N-\frac{1}{2}}_{2,\infty})}\leq  \|\n\ln\rho_{1}\|_{\widetilde{L}^{\frac{4}{3}}(\R^{+},B^{\N+\frac{1}{2}}_{2,\frac{4}{3}})\cap \widetilde{L}^{4}(\R^{+},B^{\N-\frac{1}{2}}_{2,\frac{4}{3}})}
 \end{aligned}
 \label{supertechindex}
 \end{equation}
 Of course when we assume that $\ln\rho_{1}$ belongs in $\widetilde{L}^{\infty}(\R^{+},\widetilde{B}^{\N-1,\N}_{2,1})\cap \widetilde{L}^{1}(\R^{+},\widetilde{B}^{\N+1,\N+2}_{2,1})$, it implies by interpolation that $\n\ln\rho_{1}$ is in $\widetilde{L}^{\frac{4}{3}}(\R^{+},B^{\N+\frac{1}{2}}_{2,\frac{4}{3}})\cap \widetilde{L}^{4}(\R^{+},B^{\N-\frac{1}{2}}_{2,\frac{4}{3}})$.
 In particular it proves that with such assumptions we control the right hand side of the inequality (\ref{raimportant1}).\\
 Let us mention that since in the theorem \ref{ftheo3} we will have $u_{1}=-\mu\n\ln\rho_{1}$ with $\rho_{1}$ verifying a heat equation, in particular we could improve the condition on the initial data $\rho_{1}^{0}=h^{0}_{1}+1$ of theorem \ref{ftheo3} by assuming only $h^{1}_{0}\in B^{\N-2}_{2,1}\cap B^{\N}_{2,\frac{4}{3}}\cap L^{\infty}$ and $\rho_{0}^{1}\geq c>0$. Indeed with such condition $\n\ln\rho_{1}$ would verify exactly the quantity on the right hand side of (\ref{supertechindex}) which is sufficient to control the right hand side on
 (\ref{raimportant1}).\\
 In fact a more accurate proof of the proposition \ref{fcrucialprop} by using critical interpolation estimate would show that $h^{1}_{0}\in B^{\N-2}_{2,1}\cap B^{\N}_{2,2-\e}\cap L^{\infty}$ (with $\e>0$) and $\rho_{0}^{1}\geq c>0$ is sufficient for the theorem \ref{ftheo3}.
 \label{importtech44}
\end{remarka}
{\bf Proof:} We observe that $(h_{2},u_{2})$ are solution of the following system:
\begin{equation}
\begin{cases}
\begin{aligned}
&\p_{t}h_{2}+{\rm div}u_{2}=F(h_{2},u_{2}),\\[2mm]
&\p_{t}u_{2}-\mu\D u_{2}-\mu\n{\rm div}u_{2}-\kappa\n\D h_{2}+K\n h_{2}=G(h_{2},u_{2}),\\
&(h_{2}(0,\cdot),u_{2}(0,\cdot))=(h_{0}^{2},u_{0}^{2}).
\end{aligned}
\end{cases}
\label{bafond3systeme3}
\end{equation}
with:
$$
\begin{aligned}
&F(h_{2},u_{2})=F+\mu\n\ln\rho_{1}\cdot\n h_{2}-u_{2}\cdot\n\ln\rho_{1},\\
&G(h_{2},u_{2})=G+2\mu\n \ln\rho_{1}\cdot D u_{2}+2\mu\n h_{2}\cdot D u_{1}-u_{1}\cdot\n u_{2}-u_{2}\cdot\n u_{1}\\
&\hspace{9cm}+\mu^{2}\n(\n\ln\rho_{1}\cdot\n h_{2}).
\end{aligned}
$$
By applying the proposition \ref{1flinear3}, we have for any $T>0$:
\begin{equation}
\begin{aligned}
&\|h_{2}\|_{\widetilde{L}^{1}_{T}(\widetilde{B}^{\N+1,\N+2}_{2,1})}+\|h_{2}\|_{\widetilde{L}^{\infty}_{T}(B^{\N-1,\N}_{2,1})}
+\|u_{2}\|_{\widetilde{L}^{1}_{T}(B^{\N+1}_{2,1})}+ \|u_{2}\|_{\widetilde{L}^{\infty}_{T}(B^{\N-1}_{2,1})}\\
&\leq C\big(\|h^{2}_{0}\|_{\widetilde{B}^{\N-1,\N}_{2,1}}+\|u^{2}_{0}\|_{B^{\N}_{2,1}}+\|F(h_{2},u_{2})\|_{\widetilde{L}^{1}_{T}(\widetilde{B}^{\N-1,\N}_{2,1})}
+\|G(h_{2},u_{2})\|_{\widetilde{L}^{1}_{T}(B^{\N-1}_{2,1})}\big).
\end{aligned}
\label{important}
\end{equation}
We have only to deal with the right hand side of  (\ref{important}), we have in particular:
\begin{equation}
\begin{aligned}
&\|\n\ln\rho_{1}\cdot\n h_{2}\|_{\widetilde{L}^{1}_{T}(\widetilde{B}^{\N-1,\N}_{2,1})}=\int_{0}^{T}\|\n\ln\rho_{1}\cdot\n h_{2}\|_{\widetilde{B}^{\N-1,\N}_{2,1}}(s)ds,\\
& \lesssim \int_{0}^{T}\big(\|\n h_{2}\|_{\widetilde{B}^{\N-\frac{1}{2},\N+\frac{1}{2}}_{2,1}} \|\n\ln\rho_{1}\|_{B^{\N-\frac{1}{2}}_{2,\infty}}\\
&\hspace{5cm}+\|\n h_{2}\|_{\widetilde{B}^{\N-\frac{3}{2},\N-\frac{1}{2}}_{2,1}} \|\n\ln\rho_{1}\|_{B^{\N+\frac{1}{2}}_{2,\infty}}\big)(s)ds.
%& \lesssim \|\n h_{2}\|_{\widetilde{L}^{\frac{4}{3}}(\widetilde{B}^{\N-\frac{1}{2},\N+\frac{1}{2}}_{2,1})} \|\n\ln\rho_{1}\|_{\widetilde{L}^{4}(B^{\N-\frac{1}{2}}_{2,\infty})}\\
%&\hspace{5cm}
%+\|\n h_{2}\|_{\widetilde{L}^{4}(\widetilde{B}^{\N-\frac{3}{2},\N-\frac{1}{2}}_{2,1})} \|\n\ln\rho_{1}\|_{\widetilde{L}^{\frac{4}{3}}(B^{\N+\frac{1}{2}}_{2,\infty})}.
\end{aligned}
\label{etech1}
\end{equation}
By interpolation see the proposition \ref{interpolation} we have:
\begin{equation}
\begin{aligned}
&\|\n h_{2}\|_{\widetilde{B}^{\N-\frac{3}{2},\N-\frac{1}{2}}_{2,1}}\lesssim\|\n h_{2}\|_{\widetilde{B}^{\N-2,\N-1}_{2,1}}^{\frac{3}{4}}\|\n h_{2}\|_{\widetilde{B}^{\N,\N+1}_{2,1}}^{\frac{1}{4}},\\
&\|\n h_{2}\|_{B^{\N-\frac{1}{2},\N+\frac{1}{2}}_{2,1})}\leq\|\n h_{2} \|_{\widetilde{B}^{\N-2,\N-1}_{2,1}}^{\frac{1}{4}}\|\n h_{2}\|_{\widetilde{B}^{\N,\N+1}_{2,1}}^{\frac{3}{4}}.
\end{aligned}
\label{cabinter}
\end{equation}
%\begin{equation}
%\begin{aligned}
%&\|\n h_{2}\|_{\widetilde{L}^{4}_{T}(\widetilde{B}^{\N-\frac{3}{2},\N-\frac{1}{2}}_{2,1})}\lesssim\|\n h_{2}\|_{\widetilde{L}^{\infty}_{T}(\widetilde{B}^{\N-2,\N-1}_{2,\infty})}^{\frac{3}{4}}\|\n h_{2}\|_{\widetilde{L}^{1}_{T}(\widetilde{B}^{\N,\N+1}_{2,\infty})}^{\frac{1}{4}}\leq A_{0}^{\frac{3}{4}}\e^{\frac{1}{4}}+\sqrt{\e},\\
%&\|u^{n-1}\|_{\widetilde{L}^{\frac{4}{3}}_{T}(B^{\N+\frac{1}{2}}_{2,\infty})}\leq\|u^{n-1}\|_{\widetilde{L}^{\infty}_{T}(B^{\N-1}_{2,\infty})}^{\frac{1}{4}}\|u^{n-1}\|_{\widetilde{L}^{1}_{T}(B^{\N+1}_{2,\infty})}^{\frac{3}{4}}\leq A_{0}^{\frac{1}{4}}\e^{\frac{3}{4}}+\sqrt{\e}.
%\end{aligned}
%\label{abinter}
%\end{equation}
By combining (\ref{cabinter}), (\ref{etech1}) and Young inequality we have for any $\e>0$:
\begin{equation}
\begin{aligned}
&\|\n\ln\rho_{1}\cdot\n h_{2}\|_{\widetilde{L}^{1}_{T}(\widetilde{B}^{\N-1,\N}_{2,1})}=\int_{0}^{T}\|\n\ln\rho_{1}\cdot\n h_{2}\|_{\widetilde{B}^{\N-1,\N}_{2,1}}(s)ds,\\
& \lesssim \int_{0}^{T} \big(\|\n h_{2} \|_{\widetilde{B}^{\N-2,\N-1}_{2,1}}^{\frac{1}{4}}\|\n h_{2}\|_{\widetilde{B}^{\N,\N+1}_{2,1}}^{\frac{3}{4}} \|\n\ln\rho_{1}\|_{B^{\N-\frac{1}{2}}_{2,\infty}}\\
&\hspace{4cm}+\|\n h_{2}\|_{\widetilde{B}^{\N-2,\N-1}_{2,1}}^{\frac{3}{4}}\|\n h_{2}\|_{\widetilde{B}^{\N,\N+1}_{2,1}}^{\frac{1}{4}}\|\n\ln\rho_{1}\|_{B^{\N+\frac{1}{2}}_{2,\infty}}\big)(s)ds.\\[2mm]
& \lesssim \int_{0}^{T} \big(2\e\|\n h_{2}\|_{\widetilde{B}^{\N,\N+1}_{2,1}}+C_{\e} \|\n\ln\rho_{1}\|_{B^{\N-\frac{1}{2}}_{2,\infty}}^{4}\|\n h_{2} \|_{\widetilde{B}^{\N-2,\N-1}_{2,1}}\\
&\hspace{6cm}+C_{\e}\|\n h_{2}\|_{\widetilde{B}^{\N-2,\N-1}_{2,\infty}}\|\n\ln\rho_{1}\|^{\frac{4}{3}}_{B^{\N+\frac{1}{2}}_{2,\infty}}\big)(s)ds.
%& \lesssim \|\n h_{2}\|_{\widetilde{L}^{\frac{4}{3}}(\widetilde{B}^{\N-\frac{1}{2},\N+\frac{1}{2}}_{2,1})} \|\n\ln\rho_{1}\|_{\widetilde{L}^{4}(B^{\N-\frac{1}{2}}_{2,\infty})}\\
%&\hspace{5cm}
%+\|\n h_{2}\|_{\widetilde{L}^{4}(\widetilde{B}^{\N-\frac{3}{2},\N-\frac{1}{2}}_{2,1})} \|\n\ln\rho_{1}\|_{\widetilde{L}^{\frac{4}{3}}(B^{\N+\frac{1}{2}}_{2,\infty})}.
\end{aligned}
\label{etech1}
\end{equation}
In a similar way by interpolation we have:
\begin{equation}
\begin{aligned}
&\|u_{2}\cdot\n\ln\rho_{1}\|_{\widetilde{L}^{1}_{T}(\widetilde{B}^{\N-1,\N}_{2,1})}=\int_{\R^{+}}\|u_{2}\cdot\n\ln\rho_{1}\|_{\widetilde{B}^{\N-1,\N}_{2,1}}(s)ds,\\
&\lesssim\int_{0}^{T} \big(\|\n \ln\rho_{1}\|_{\widetilde{B}^{\N-\frac{1}{2},\N+\frac{1}{2}}_{2,\infty}} \|u_{2}\|_{B^{\N-\frac{1}{2}}_{2,1}}
+\|\n \ln\rho_{1}\|_{\widetilde{B}^{\N-\frac{3}{2},\N-\frac{1}{2}}_{2,\infty}} \|u_{2}\|_{B^{\N+\frac{1}{2}}_{2,1}}\big)(s)ds.\\[2mm]
&\lesssim\int_{0}^{T} \big(\|\n \ln\rho_{1}\|_{\widetilde{B}^{\N-\frac{1}{2},\N+\frac{1}{2}}_{2,\infty}} \|u_{2}\|^{\frac{3}{4}}_{B^{\N-1}_{2,1}} \|u_{2}\|^{\frac{1}{4}}_{B^{\N+1}_{2,1}}\\
&\hspace{4cm}
+\|\n \ln\rho_{1}\|_{\widetilde{B}^{\N-\frac{3}{2},\N-\frac{1}{2}}_{2,\infty}} \|u_{2}\|^{\frac{1}{4}}_{B^{\N-1}_{2,1}} \|u_{2}\|^{\frac{3}{4}}_{B^{\N+1}_{2,1}}\big)(s)ds,
\end{aligned}
\end{equation}
Applying Young inequality it yields:
\begin{equation}
\begin{aligned}
&\|u_{2}\cdot\n\ln\rho_{1}\|_{\widetilde{L}^{1}_{T}(\widetilde{B}^{\N-1,\N}_{2,1})} \lesssim \int_{0}^{T} \big(2\e\|u_{2}\|_{B^{\N+1}_{2,1}}+C_{\e} \|\n\ln\rho_{1}\|_{\widetilde{B}^{\N-\frac{3}{2},\N-\frac{1}{2}}_{2,\infty}}^{4}\|u_{2} \|_{B^{\N-1}_{2,1}}\\
&\hspace{6cm}+C_{\e}\|u_{2}\|_{B^{\N-1}_{2,\infty}}\|\n\ln\rho_{1}\|^{\frac{4}{3}}_{\widetilde{B}^{\N-\frac{1}{2},\N+\frac{1}{2}}_{2,\infty}}\big)(s)ds.
\end{aligned}
\label{etech2}
\end{equation}
Let us proceed in a similar way for $\|G(h_{2},u_{2})\|_{\widetilde{L}^{1}_{T}(B^{\N-1}_{2,1})}$, we are going only to treat only two terms (the other one will be left to the reader). In a similar way, we have:
\begin{equation}
\begin{aligned}
&\|\n \ln\rho_{1}\cdot D u_{2}\|_{\widetilde{L}^{1}_{T}(B^{\N-1}_{2,1})}=\int_{0}^{T}\|\n \ln\rho_{1}\cdot D u_{2}\|_{B^{\N-1}_{2,1}}(s)ds,\\
&\lesssim\int_{0}^{T} \big(\|\n \ln\rho_{1}\|_{B^{\N-\frac{1}{2}}_{2,\infty}} \| D u_{2}\|_{B^{\N-\frac{1}{2}}_{2,1}}
+\|\n \ln\rho_{1}\|_{B^{\N+\frac{1}{2}}_{2,\infty}} \|D u_{2}\|_{B^{\N-\frac{3}{2}}_{2,1}}\big)(s)ds,\\[2mm]
%&\lesssim\int_{\R^{+}} \big(\|\n \ln\rho_{1}\|_{\widetilde{B}^{\N-\frac{1}{2},\N+\frac{1}{2}}_{2,\infty}} \|u_{2}\|^{\frac{3}{4}}_{B^{\N-1}_{2,1}} \|u_{2}\|^{\frac{1}{4}}_{B^{\N+1}_{2,1}}\\
%&\hspace{4cm}
%+\|\n \ln\rho_{1}\|_{\widetilde{B}^{\N-\frac{3}{2},\N-\frac{1}{2}}_{2,\infty}} \|u_{2}\|^{\frac{1}{4}}_{B^{\N-1}_{2,1}} \|u_{2}\|^{\frac{3}{4}}_{B^{\N+1}_{2,1}}\big)(s)ds,\\
& \lesssim \int_{0}^{T} \big(2\e\|u_{2}\|_{B^{\N+1}_{2,1}}+C_{\e} \|\n\ln\rho_{1}\|_{\widetilde{B}^{\N-\frac{3}{2},\N-\frac{1}{2}}_{2,\infty}}^{4}\|u_{2} \|_{B^{\N-1}_{2,1}}\\
&\hspace{6cm}+C_{\e}\|u_{2}\|_{B^{\N-1}_{2,\infty}}\|\n\ln\rho_{1}\|^{\frac{4}{3}}_{\widetilde{B}^{\N-\frac{1}{2},\N+\frac{1}{2}}_{2,\infty}}\big)(s)ds.
\end{aligned}
\label{etech3}
\end{equation}
and:
\begin{equation}
\begin{aligned}
&\|u_{1}\cdot\n u_{2}\|_{\widetilde{L}^{1}_{T}(B^{\N-1}_{2,1})}=\int_{0}^{T}\|u_{1}\cdot\n u_{2}\|_{B^{\N-1}_{2,1}}(s)ds,\\
%&\lesssim\int_{\R^{+}} \big(\|\n \ln\rho_{1}\|_{\widetilde{B}^{\N-\frac{1}{2},\N+\frac{1}{2}}_{2,\infty}} \|u_{2}\|_{B^{\N-\frac{1}{2}}_{2,1}}
%+\|\n \ln\rho_{1}\|_{\widetilde{B}^{\N-\frac{1}{2},\N-\frac{1}{2}}_{2,\infty}} \|u_{2}\|_{B^{\N+\frac{1}{2}}_{2,1}}\big)(s)ds.\\[2mm]
%&\lesssim\int_{\R^{+}} \big(\|\n \ln\rho_{1}\|_{\widetilde{B}^{\N-\frac{1}{2},\N+\frac{1}{2}}_{2,\infty}} \|u_{2}\|^{\frac{3}{4}}_{B^{\N-1}_{2,1}} \|u_{2}\|^{\frac{1}{4}}_{B^{\N+1}_{2,1}}\\
%&\hspace{4cm}
%+\|\n \ln\rho_{1}\|_{\widetilde{B}^{\N-\frac{3}{2},\N-\frac{1}{2}}_{2,\infty}} \|u_{2}\|^{\frac{1}{4}}_{B^{\N-1}_{2,1}} \|u_{2}\|^{\frac{3}{4}}_{B^{\N+1}_{2,1}}\big)(s)ds,\\
& \lesssim \int_{0}^{T} \big(2\e\|u_{2}\|_{B^{\N+1}_{2,1}}+C_{\e} \|u_{1}\|_{B^{\N-\frac{1}{2}}_{2,\infty}}^{4}\|u_{2} \|_{B^{\N-1}_{2,1}}+C_{\e}\|u_{2}\|_{B^{\N-1}_{2,1}}\|u_{1}\|^{\frac{4}{3}}_{B^{\N+\frac{1}{2}}_{2,\infty}}\big)(s)ds.
\end{aligned}
\label{etech4}
\end{equation}
Finally combining (\ref{important}), (\ref{etech1}), (\ref{etech2}), (\ref{etech3}) and (\ref{etech4}) we obtain that for $C>0$ and $\e>0$ small enough such that $C\e\leq\frac{1}{2}$ that:
\begin{equation}
\begin{aligned}
&\|h_{2}\|_{\widetilde{L}^{1}_{T}(\widetilde{B}^{\N+1,\N+2}_{2,1})}+\|h_{2}\|_{\widetilde{L}^{\infty}_{T}(B^{\N-1,\N}_{2,1})}
+\|u_{2}\|_{\widetilde{L}^{1}_{T}(B^{\N+1}_{2,1})}+ \|u_{2}\|_{\widetilde{L}^{\infty}_{T}(B^{\N-1}_{2,1})}\\
&\leq C\biggl(\|h^{2}_{0}\|_{\widetilde{B}^{\N-1,\N}_{2,1}}+\|u^{2}_{0}\|_{B^{\N}_{2,1}}+\|F\|_{\widetilde{L}^{1}_{T}(\widetilde{B}^{\N-1,\N}_{2,1})}
+\|G\|_{\widetilde{L}^{1}_{T}(B^{\N-1}_{2,1})}\\
&+\int_{0}^{T} \biggl(\e\big(\|u_{2}\|_{B^{\N+1}_{2,1}}+\|h_{2}\|_{\widetilde{L}^{1}_{T}(\widetilde{B}^{\N+1,\N+2}_{2,1})}\big)+C_{\e}\big(\|u_{2} \|_{B^{\N-1}_{2,1}}+\|h_{2} \|_{\widetilde{B}^{\N-1,\N}_{2,1}}\big)\\
&\hspace{1cm}\times\big( \|u_{1}\|_{B^{\N-\frac{1}{2}}_{2,\infty}}^{4}+\|u_{1}\|^{\frac{4}{3}}_{B^{\N+\frac{1}{2}}_{2,\infty}}+\|\n\ln\rho_{1}\|^{\frac{4}{3}}_{B^{\N+\frac{1}{2}}_{2,\infty}}+\|\n\ln\rho_{1}\|^{4}_{B^{\N-\frac{1}{2}}_{2,\infty}}\big)\biggl)(s)ds\biggl).
\end{aligned}
\label{raimportant}
\end{equation}
By a bootstrap argument and the Gronwall lemma where we use the fact that:
 $$\|u_{2}(s) \|_{B^{\N-1}_{2,1}}+\|h_{2}(s) \|_{\widetilde{B}^{\N-1,\N}_{2,1}}\lesssim
 \|u_{2}\|_{\widetilde{L}^{\infty}_{s}(B^{\N-1}_{2,1})}+\|h_{2} \|_{\widetilde{L}^{\infty}_{s}(\widetilde{B}^{\N-1,\N}_{2,1})}$$
 we get for $C>0$ large enough:
\begin{equation}
\begin{aligned}
&\|h_{2}\|_{\widetilde{L}^{1}_{T}(\widetilde{B}^{\N+1,\N+2}_{2,1})}+\|h_{2}\|_{\widetilde{L}^{\infty}_{T}(B^{\N-1,\N}_{2,1})}
+\|u_{2}\|_{\widetilde{L}^{1}_{T}(B^{\N+1}_{2,1})}+ \|u_{2}\|_{\widetilde{L}^{\infty}_{T}(B^{\N-1}_{2,1})}\\
&\leq C\big(\|h^{2}_{0}\|_{\widetilde{B}^{\N-1,\N}_{2,1}}+\|u^{2}_{0}\|_{B^{\N}_{2,1}}+\|F\|_{\widetilde{L}^{1}_{T}(\widetilde{B}^{\N-1,\N}_{2,1})}
+\|G\|_{\widetilde{L}^{1}_{T}(B^{\N-1}_{2,1})}\big)\\
&\times \exp\biggl(C\int_{0}^{T}\big( \|u_{1}\|_{B^{\N-\frac{1}{2}}_{2,\infty}}^{4}+\|u_{1}\|^{\frac{4}{3}}_{B^{\N+\frac{1}{2}}_{2,\infty}}+\|\n\ln\rho_{1}\|^{\frac{4}{3}}_{B^{\N+\frac{1}{2}}_{2,\infty}}+\|\n\ln\rho_{1}\|^{4}_{B^{\N-\frac{1}{2}}_{2,\infty}}\big)(s)ds\biggl).\\
&
\end{aligned}
\label{raimportant}
\end{equation}
It concludes the proof of the proposition \ref{fcrucialprop}. \hfill {$\Box$}
\subsection{Existence of local solutions for system (\ref{NHV1})}
\label{sub32}
We now are going to prove the existence of strong solutions in finite time with large initial data verifying the hypothesis of theorem \ref{ftheo1} for the system (\ref{NHV1}). More precisely we assume that $(q_{0},u_{0})$ belong in $B^{\N}_{2,\infty}\times B^{\N-1}_{2,\infty}$. 
%In particular we recall that the main interest of theorem \ref{ftheo1} is to allow discontinuous initial data for the density, such that we can authorize discontinuous interfaces.
\subsubsection*{Existence of solutions}
The existence part of the theorem is proved by an iterative method. We define a sequence $(q^{n},u^{n})$ as follows:
$$q^{n}=q_{L}+\bar{q}^{n},\;u^{n}=u_{L}+\bar{u}^{n},$$
where $(q_{L},u_{L})$ stands for the  solution of:
\begin{equation}
\begin{cases}
\p_{t}q_{L}+{\rm div}u_{L}=0,\\
\p_{t}u_{L}-{\cal A}u_{L}-\kappa\n(\D q_{L})=0,
\end{cases}
\label{lineaire}
\end{equation}
supplemented with initial data:
$$q_{L}(0)=q_{0}\;,\;u_{L}(0)=u_{0}.$$
Here ${\cal A}$ define the Lam\'e operator ${\cal A}u=\mu\D u+(\lambda+\mu)\n{\rm div}u$. Using the proposition \ref{flinear3}, we obtain the following estimates on $(q_{L},u_{L})$ for all $T>0$:
$$q_{L}\in\widetilde{C}([0,T],B^{\N}_{2,\infty})\cap\widetilde{L}^{1}_{T}(
B^{\N+2}_{2,\infty})\;\;\mbox{and}\;\;u_{L}\in\widetilde{C}([0,T],B^{\N-1}_{2,\infty})\cap
\widetilde{L}^{1}_{T}
(B^{\N+1}_{2,\infty}).$$
Setting $(\bar{q}^{0},\bar{u}^{0})=(0,0)$ we now define $(\bar{q}_{n},\bar{u}_{n})$ as the solution of the following system:
$$
\begin{cases}
\begin{aligned}
&\p_{t}\bar{q}^{n}+{\rm div}(\bar{u}^{n})=F_{n-1},\\
& \p_{t}\bar{u}_{n}-{\cal A}\bar{u}_{n}-\kappa\n(\D\bar{q}^{n})=G_{n-1},\\
&(\bar{q}_{n},\bar{u}_{n})_{t=0}=(0,0),
\end{aligned}
\end{cases}
\leqno{(N_{1})}
$$
where:
$$
\begin{aligned}
F_{n-1}=&-u^{n-1}\cdot\n q^{n-1},\\
%=&-u_{L}\cdot\n q_{L}-\bar{u}^{n-1}\cdot\n q_{L}-u^{L}\cdot\n \bar{q}^{n-1} -\bar{u}^{n-1}\cdot\n \bar{q}^{n-1}  ,\\
G_{n-1}=&-(u^{n-1})^{*}.\n u^{n-1}+2\mu\n q^{n-1}\cdot D u^{n-1}+\lambda\n q^{n-1}\,{\rm div}u^{n-1}+\frac{\kappa}{2}\n(|\n q^{n-1}|^{2})\\
&\hspace{11cm}-K\n q^{n-1}.
\end{aligned}
$$
%Ajouter les autres termes venant de la viscosite.
\subsubsection*{1) First Step , Uniform Bound}
 Let $\e$ be a small
parameter and  choose $T$ small enough  such that according to the proposition \ref{flinear3} we have:
$$
\begin{aligned}
&\|u_{L}\|_{\widetilde{L}^{1}_{T}(B^{\N+1}_{2,\infty})}+\|q_{L}\|_{\widetilde{L}^{1}_{T}(B^{\N+2}_{2,\infty})}\leq\e,\\
&\|u_{L}\|_{\widetilde{L}^{\infty}_{T}(B^{\N-1}_{2,\infty})}+\|q_{L}\|_{\widetilde{L}^{\infty}_{T}(B^{\N}_{2,\infty})}\leq
C A_{0},
\end{aligned}
\leqno{({\cal{H}}_{\e})}
$$
with $A_{0}=\|q_{0}\|_{B^{\N}_{2,\infty}}+\|u_{0}\|_{B^{\N-1}_{2,\infty}}$. We are going to show by induction that:
$$\|(\bar{q}^{n},\bar{u}^{n})\|_{F_{T}}\leq\sqrt{\e}.\leqno{({\cal{P}}_{n})},$$
for $\e$ small enough with:
$$F_{T}=\big(\widetilde{C}([0,T],B^{\N}_{2,\infty})\cap\widetilde{L}^{1}_{T}(
B^{\N+2}_{2,\infty})\big)\times\big(\widetilde{C}([0,T],B^{\N-1}_{2,\infty})\cap
\widetilde{L}^{1}_{T}
(B^{\N+1}_{2,\infty})\big)^{N}.$$
As $(\bar{q}^{0},\bar{u}^{0})=(0,0)$ the result is true for $n=0$. We now suppose $({\cal P}_{n-1})$ (with $n\geq 1$) true and we are going to show  $({\cal P}_{n})$.
Applying proposition \ref{flinear3}  we have:
\begin{equation}
\begin{aligned}
&\|(\bar{q}^{n},\bar{u}^{n})\|_{F_{T}}\leq C\|(\n
F_{n-1},G_{n-1})\|_{\widetilde{L}^{1}_{T}(B^{\N-1}_{2,\infty})}.
\end{aligned}
\label{fi11a}
\end{equation}
Bounding the right-hand side of (\ref{fi11a}) may be
done by applying
proposition  \ref{produit}, lemma  \ref{composition} and corollary \ref{produit2}. We begin with treating the case of $\|F_{n-1}\|_{\widetilde{L}^{1}_{T}(B^{N/2}_{2,\infty})}$, let us recall that:
\begin{equation}
F_{n-1}=-u_{L}\cdot\n q_{L}-\bar{u}^{n-1}\cdot\n q_{L}-u^{L}\cdot\n \bar{q}^{n-1} -\bar{u}^{n-1}\cdot\n \bar{q}^{n-1}.
\label{Fn}
\end{equation}
We are going to bound each term of we have then:
\begin{equation}
\begin{aligned}
&\|u_{L}\cdot\n q_{L}\|_{\widetilde{L}^{1}_{T}(B^{N/2}_{2,\infty})}\leq \|u_{L}\|_{\widetilde{L}^{1}_{T}(B^{N/2+1}_{2,\infty})}\|q_{L}\|_{\widetilde{L}^{\infty}_{T}(B^{N/2}_{2,\infty})}\\
&\hspace{7cm}
+
\|q_{L}\|_{\widetilde{L}^{1}_{T}(B^{N/2+2}_{2,\infty})} \|u_{L}\|_{\widetilde{L}^{\infty}_{T}(B^{N/2-1}_{2,\infty})}.
\end{aligned}
\label{F1}
\end{equation}
Similarly we obtain:
\begin{equation}
\begin{aligned}
&\|u_{L}\cdot\n \bar{q}^{n-1}\|_{\widetilde{L}^{1}_{T}(B^{N/2}_{2,\infty})}\leq \|u_{L}\|_{\widetilde{L}^{1}_{T}(B^{N/2+1}_{2,\infty})}\|\bar{q}^{n-1}\|_{\widetilde{L}^{\infty}_{T}(B^{N/2}_{2,\infty})}\\
&\hspace{7cm}+\|\n \bar{q}^{n-1}\|_{\widetilde{L}^{\frac{4}{3}}_{T}(B^{\N+\frac{1}{2}}_{2,\infty})}\|u_{L}\|_{\widetilde{L}^{4}_{T}(B^{\N-\frac{1}{2}}_{2,\infty})},
\end{aligned}
\label{F2}
\end{equation}
\begin{equation}
\begin{aligned}
&\|\bar{u}^{n-1}\cdot\n q_{L}\|_{\widetilde{L}^{1}_{T}(B^{N/2}_{2,\infty})}\leq \|\bar{u}^{n-1}\|_{\widetilde{L}^{\frac{4}{3}}_{T}(B^{\N+\frac{1}{2}}_{2,\infty})}\|\n q_{L}\|_{\widetilde{L}^{4}_{T}(B^{\N-\frac{1}{2}}_{2,\infty})}\\
&\hspace{7cm}+\|q_{L}\|_{\widetilde{L}^{1}_{T}(B^{N/2+2}_{2,\infty})} \|\bar{u}^{n-1}\|_{\widetilde{L}^{\infty}_{T}(B^{N/2-1}_{2,\infty})},
\end{aligned}
\label{F3}
\end{equation}
and:
\begin{equation}
\begin{aligned}
&\|\bar{u}^{n-1}\cdot\n \bar{q}^{n-1}\|_{\widetilde{L}^{1}_{T}(B^{N/2}_{2,\infty})}\leq \|\bar{u}^{n-1}\|_{\widetilde{L}^{1}_{T}(B^{N/2+1}_{2,\infty})}\|\bar{q}^{n-1}\|_{\widetilde{L}^{\infty}_{T}(B^{N/2}_{2,\infty})}\\
&\hspace{6cm}+\|\bar{q}^{n-1}\|_{\widetilde{L}^{1}_{T}(B^{N/2+2}_{2,\infty})} \|\bar{u}^{n-1}\|_{\widetilde{L}^{\infty}_{T}(B^{N/2-1}_{2,\infty})}.
\end{aligned}
\label{F4}
\end{equation}
By using the previous inequalities (\ref{F1}), (\ref{F2}), (\ref{F3}), (\ref{F4}), $({\cal P}_{n-1})$ and by interpolation, we obtain that for $C>0$ large enough:
\begin{equation}
\|F_{n}\|_{\widetilde{L}^{1}_{T}(B^{N/2}_{2,\infty})}\leq C\sqrt{\e}(A_{0}^{\frac{3}{4}}\e^{\frac{1}{4}}+\sqrt{\e}(1+A_{0})+\e).
\label{Fn}
\end{equation}
Next we want to control  $\|G_{n}\|_{\widetilde{L}^{1}(B^{\N-1}_{2,\infty})}$. According to
propositions \ref{produit}, corollary \ref{produit2} and \ref{flinear3}, we have:
\begin{equation}
\begin{aligned}
&\|(u^{n-1})^{*}.\n
u^{n-1}\|_{\widetilde{L}^{1}_{T}(B^{\N-1}_{2,\infty})}\lesssim\|u^{n-1}\|_{\widetilde{L}^{\frac{4}{3}}_{T}(B^{\N+\frac{1}{2}}_{2,\infty})}\|u^{n-1}\|_{\widetilde{L}^{4}_{T}(B^{\N-\frac{1}{2}}_{2,\infty})},\\[2mm]
&\|\n(|\n q^{n-1}|^{2})\|_{\widetilde{L}^{1}_{T}(B^{\N-1}_{2,\infty})}\lesssim\||\n q^{n-1}|^{2}
\|_{\widetilde{L}^{1}_{T}(B^{\N}_{2,\infty})},\\
&\hspace{3,6cm}\lesssim \|\n q^{n-1}
\|_{\widetilde{L}^{\frac{4}{3}}_{T}(B^{\N+\frac{1}{2}}_{2,\infty})}\|\n q^{n-1}
\|_{\widetilde{L}^{4}_{T}(B^{\N-\frac{1}{2}}_{2,\infty})},\\
&\hspace{3,6cm}\lesssim \|q^{n-1}
\|_{\widetilde{L}^{\frac{4}{3}}_{T}(B^{\N+\frac{3}{2}}_{2,\infty})}\|q^{n-1}
\|_{\widetilde{L}^{4}_{T}(B^{\N+\frac{1}{2}}_{2,\infty})},\\[2mm]
&\|\n q^{n-1}\cdot D u^{n-1}\|_{\widetilde{L}^{1}_{T}(B^{\N-1}_{2,\infty})}\lesssim \|\n q^{n-1}\|_{\widetilde{L}^{\frac{4}{3}}_{T}(B^{\N+\frac{1}{2}}_{2,\infty})}\|u^{n-1}\|_{\widetilde{L}^{4}_{T}(B^{\N-\frac{1}{2}}_{2,\infty})}\\
&\hspace{6cm}+\|u^{n-1}\|_{\widetilde{L}^{\frac{4}{3}}_{T}(B^{\N+\frac{1}{2}}_{2,\infty})}\|\n q^{n-1}\|_{\widetilde{L}^{4}_{T}(B^{\N-\frac{1}{2}}_{2,\infty})},\\[2mm]
&\|\n q^{n-1}{\rm div} u^{n-1}\|_{\widetilde{L}^{1}_{T}(B^{\N-1}_{2,\infty})}\lesssim \|\n q^{n-1}\|_{\widetilde{L}^{\frac{4}{3}}_{T}(B^{\N+\frac{1}{2}}_{2,\infty})}\|u^{n-1}\|_{\widetilde{L}^{4}_{T}(B^{\N-\frac{1}{2}}_{2,\infty})}\\
&\hspace{6cm}+\|u^{n-1}\|_{\widetilde{L}^{\frac{4}{3}}_{T}(B^{\N+\frac{1}{2}}_{2,\infty})}\|\n q^{n-1}\|_{\widetilde{L}^{4}_{T}(B^{\N-\frac{1}{2}}_{2,\infty})},\\[2mm]
&\|\n q^{n-1}\|_{\widetilde{L}^{1}_{T}(B^{\N-1}_{2,\infty})}\leq T\|q^{n-1}\|_{\widetilde{L}^{\infty}_{T}(B^{\N}_{2,\infty})}.
\end{aligned}
\label{aGn}
\end{equation}
Let us recall that we have by interpolation and the condition $({\cal P}_{n-1})$:
\begin{equation}
\begin{aligned}
&\|u^{n-1}\|_{\widetilde{L}^{4}_{T}(B^{\N-\frac{1}{2}}_{2,\infty})}\leq\|u^{n-1}\|_{\widetilde{L}^{\infty}_{T}(B^{\N-1}_{2,\infty})}^{\frac{3}{4}}\|u^{n-1}\|_{\widetilde{L}^{1}_{T}(B^{\N+1}_{2,\infty})}^{\frac{1}{4}}\leq A_{0}^{\frac{3}{4}}\e^{\frac{1}{4}}+\sqrt{\e},\\
&\|u^{n-1}\|_{\widetilde{L}^{\frac{4}{3}}_{T}(B^{\N+\frac{1}{2}}_{2,\infty})}\leq\|u^{n-1}\|_{\widetilde{L}^{\infty}_{T}(B^{\N-1}_{2,\infty})}^{\frac{1}{4}}\|u^{n-1}\|_{\widetilde{L}^{1}_{T}(B^{\N+1}_{2,\infty})}^{\frac{3}{4}}\leq A_{0}^{\frac{1}{4}}\e^{\frac{3}{4}}+\sqrt{\e}.
\end{aligned}
\label{abinter}
\end{equation}
Using (\ref{fi11a}), (\ref{Fn}), (\ref{aGn}) and (\ref{abinter}) we obtain for a certain $C>0$ large enough:
$$
\begin{aligned}
\|(\bar{q}^{n},\bar{u}^{n})\|_{F_{T}}&\leq C\sqrt{\e}(A_{0}^{\frac{3}{4}}\e^{\frac{1}{4}}+\sqrt{\e}(1+A_{0})+\e)+C( A_{0}^{\frac{3}{4}}\e^{\frac{1}{4}}+\sqrt{\e})(A_{0}^{\frac{1}{4}}\e^{\frac{3}{4}}+\sqrt{\e})\\
&\hspace{9cm}+T (A_{0}+\sqrt{\e}),\\
&\leq C\sqrt{\e}(A_{0}^{\frac{3}{4}}\e^{\frac{1}{4}}+2\sqrt{\e}(1+A_{0})+A_{0}^{\frac{3}{4}}\e^{\frac{1}{4}}+A_{0}^{\frac{1}{4}}\e^{\frac{3}{4}}+\e)+T (A_{0}+\sqrt{\e}).
\end{aligned}
$$
By choosing $T$ and $\e$ small enough  the property $({\cal{P}}_{n})$ is verified, so we
have shown by induction that $(q^{n},u^{n})$ is bounded
in $F_{T}$.%\\
%\\
%In the case where we want to obtain global strong solution with small initial data, we need to take in account the low frequencies. That is why we will include the pressure term in the linear part and we will use proposition \ref{1flinear3} to obtain estimate on $(q_{L},u_{L})$ where $(q_{L},u_{L})$ here stands for the  solution of:
%\begin{equation}
%\begin{cases}
%\p_{t}q_{L}+{\rm div}u_{L}=0,\\
%\p_{t}u_{L}-{\cal A}u_{L}+K\n q_{L}-\kappa\n(\D q_{L})=0,
%\end{cases}
%\label{lineaire}
%\end{equation}
%In particular it explains why we need of additional regularity on the initial density, i.e $q_{0}\in B^{\N-1}_{2,\infty}$. The rest of the proof for this case is similar to the case with large initial data.
\subsubsection*{Second Step: Convergence of the
sequence}
 We will show
that $(q^{n},u^{n})$ is a Cauchy sequence in the Banach
space $F_{T}$, hence converges to some
$(q,u)\in F_{T}$. Let:
$$\delta q^{n}=q^{n+1}-q^{n},\;\delta u^{n}=u^{n+1}-u^{n}.$$
The system verified by $(\de q^{n},\de u^{n})$ reads:
$$
\begin{cases}
\begin{aligned}
&\p_{t}\delta q^{n}+{\rm div}\delta u^{n}=F_{n}-F_{n-1},\\
&\p_{t}\delta u^{n}-\mu\D\delta u^{n}-(\lambda+\mu)\n{\rm div}\delta u^{n} -\kappa\n \D \delta q^{n}=G_{n}-G_{n-1},\\
&\delta q^{n}(0)=0\;,\;\delta u^{n}(0)=0,
\end{aligned}
\end{cases}
$$
Applying propositions \ref{flinear3}, we obtain:
\begin{equation}
\begin{aligned}
\|(\de q^{n},\de u^{n})\|_{F_{T}}\leq\;&
C(\|F_{n}-F_{n-1}\|_{\widetilde{L}^{1}_{T}(B^{N/2}_{2,\infty})}+\|G_{n}-G_{n-1}\|_{\widetilde{L}^{1}_{T}(B^{N/2-1}_{2,\infty})}).
\end{aligned}
\label{cauchy}
\end{equation}
Tedious calculus ensure that:
$$
\begin{aligned}
&F_{n}-F_{n-1}=-\delta u^{n-1}\cdot\n q^{n}-u^{n-1}\cdot\n\delta q^{n-1},\\[2mm]
&G_{n}-G_{n-1}=-u^{n}\cdot\n\delta u^{n-1}-\delta u^{n-1}\cdot\n u^{n-1}+\mu\n q^{n}\cdot D\delta u^{n-1}+\mu\n\delta q^{n-1}\cdot D u^{n-1}\\
&+\lambda\n q^{n} {\rm div}\delta u^{n-1}+\lambda \n\delta q^{n-1}{\rm div} u^{n-1}-K\n\delta q^{n-1}+\n\big(\n q^{n}\cdot\n \delta q^{n-1}+\n  \delta q^{n-1}\cdot\n q^{n-1}\big).
\end{aligned}
$$
It remains only to estimate the terms on the right hand side of (\ref{cauchy}) by using the same type of estimates than in the previous section and the property $({\cal P}_{n})$. More precisely we have via the proposition \ref{produit} and $({\cal P}_{n})$, it exists $C>0$ such that:
\begin{equation}
\begin{aligned}
&\|F_{n}-F_{n-1}\|_{\widetilde{L}^{1}_{T}(B^{N/2}_{2,\infty})}\leq \|\delta u^{n-1}\|_{\widetilde{L}^{\frac{4}{3}}_{T}(B^{\N+\frac{1}{2}}_{2,\infty})}\|\n q^{n}\|_{\widetilde{L}^{4}_{T}(B^{\N-\frac{1}{2}}_{2,\infty})}\\
&\hspace{0,5cm}+ \|\delta u^{n-1}\|_{\widetilde{L}^{\infty}_{T}(B^{\N-1}_{2,\infty})}\|\n q^{n}\|_{\widetilde{L}^{1}_{T}(B^{\N+1}_{2,\infty})}+\|u^{n-1}\|_{\widetilde{L}^{1}_{T}(B^{\N+1}_{2,\infty})}\|\n \delta q^{n-1}\|_{\widetilde{L}^{\infty}_{T}(B^{\N-1}_{2,\infty})}\\
&\hspace{6,5cm}+\|\n \delta q^{n-1}\|_{\widetilde{L}^{\frac{4}{3}}_{T}(B^{\N+\frac{1}{2}}_{2,\infty})}\|u^{n-1}\|_{\widetilde{L}^{4}_{T}(B^{\N-\frac{1}{2}}_{2,\infty})}\\[1,5mm]
&\leq C (A_{0}^{\frac{3}{4}}\e^{\frac{1}{4}}+A_{0}^{\frac{1}{4}}\e^{\frac{3}{4}}+\sqrt{\e}+\e)\|(\delta q^{n-1},\delta u^{n-1})\|_{F_{T}}.
\end{aligned}
\label{cFn}
\end{equation}
In a similar way we show that it exists $C>0$ large enough such that:
\begin{equation}
\begin{aligned}
&\|G_{n}-G_{n-1}\|_{\widetilde{L}^{1}_{T}(B^{N/2-1}_{2,\infty})}\leq C (A_{0}^{\frac{3}{4}}\e^{\frac{1}{4}}+A_{0}^{\frac{1}{4}}\e^{\frac{3}{4}}+\sqrt{\e}+\e+T)\|(\delta q^{n-1},\delta u^{n-1})\|_{F_{T}}.
\end{aligned}
\label{cGn}
\end{equation}
By combining (\ref{cauchy}), (\ref{cFn}) and (\ref{cGn}) , we get for $C>0$ large enough:
$$\|(\de q^{n},\de u^{n})\|_{F_{T}}\leq C (A_{0}^{\frac{3}{4}}\e^{\frac{1}{4}}+A_{0}^{\frac{1}{4}}\e^{\frac{3}{4}}+\sqrt{\e}+\e+T)\|(\delta q^{n-1},\delta u^{n-1})\|_{F_{T}}.$$ 
It implies that choosing $\e$ and $T$ small enough 
$(q^{n},u^{n})$ is a Cauchy sequence in $F_{T}$ which is a Banach. It provides that $(q^{n},u^{n})$ converges to
$(q,u)$ is in $F_{T}$. The verification that the limit $(q,u)$ is solution of $(NHV1)$ in the sense of distributions is a straightforward application of proposition \ref{produit}.
\subsubsection*{Third step: Uniqueness}
Now, we are going to prove the uniqueness of the solution in $F_{T}$% the following space:
%$$\widetilde{F}^{\N}_{T}=\big(\widetilde{C}([0,T],B^{\N}_{2,\infty})\cap\widetilde{L}^{2}_{T}(
%B^{\N+1}_{2,\infty})\big)\times\big(\widetilde{C}([0,T],B^{\N-1}_{2,\infty})\cap
%\widetilde{L}^{2}_{T}
%(B^{\N}_{2,\infty})\big).$$
Suppose that $(q_{1},u_{1})$ and $(q_{2},u_{2})$ are solutions with the same initial conditions and belonging in $F_{T}$ where
$(q_{1},u_{1})$ corresponds to the previous
solution. We set:
$$\de q=q_{2}-q_{1}\;\;\;\mbox{and}\;\;\;\de u=u_{2}-u_{1}.$$
We deduce that $(\de q,\de u)$ satisfy the following system:
$$
\begin{cases}
\begin{aligned}
&\p_{t}\delta q+{\rm div}\delta
u=F_{2}-F_{1},\\
&\p_{t}\delta u-\mu\D\delta u-(\lambda+\mu)\n{\rm div}\delta u -\kappa\n \D \delta q=G_{1}-G_{2},\\
&\delta q(0)=0\;,\;\delta u(0)=0.
\end{aligned}
\end{cases}
$$
We now apply proposition \ref{flinear3} to the previous system, and by using the same type of estimates than in the previous part, we
show that:
$$
\begin{aligned}
&\|(\de q,\de u)\|_{\widetilde{F}^{\N}_{T_{1}}}\lesssim(\|q_{1}\|_{\widetilde{L}^{2}{T_{1}}(B^{\N+1}_{2,\infty})}+\|q_{2}\|_{\widetilde{L}{T_{1}}^{2}(B^{\N+1}_{2,\infty})}+
\|u_{1}\|_{\widetilde{L}^{2}{T_{1}}(B^{\N+1}_{2,\infty})}+\|u_{2}\|_{\widetilde{L}{T_{1}}^{2}(B^{\N}_{2,\infty})})\\
&\hspace{11cm}\times\|(\de
q,\de u)\|_{\widetilde{F}^{\N}_{T_{1}}}.
\end{aligned}
$$
We have then for $T_{1}$ small enough: $(\de q,\de u)=(0,0)$ on $[0,T_{1}]$ and by connectivity we finally
conclude that:
$$q_{1}=q_{2},\;u_{1}=u_{2}\;\;\mbox{on}\;\;[0,T].$$ \hfill {$\Box$}
\subsection{Global solution near equilibrium for system (\ref{NHV1})}
\label{infini}
We are now interested in proving the existence of global strong solution with small initial data for the system (\ref{NHV1}). The main difference with the previous proof consists
essentially in dealing with the behavior of the density in low frequencies, to do this we shall use the proposition \ref{flinear3}. More precisely we are going to use a contracting mapping argument for the function $\psi$ defined as follows:
\begin{equation}
\psi(q,u)=W(t,\cdot)*\left(\begin{array}{c}
q_{0}\\
u_{0}\\
\end{array}
\right)+\int_{0}^{t}W(t-s)\left(\begin{array}{c}
F(q,u)\\
G(q,u)\\
\end{array}
\right)\ ds\; .\label{f1a1}
\end{equation}
where $W$ is the semi group associated to the linear system $(N1)$ with $a=\mu$, $b=\lambda+\mu$, $c=\kappa$ and $d=K$.
%In what follows we set:
%$$\rho=\bar{\rho}(1+q)\;,\;\theta=\bar{\theta}+{\cal T}\;,\;\widetilde{T}=\Psi^{-1}(\theta).$$
The non linear terms  $F,G$ are defined as follows:
\begin{equation}
\begin{aligned}
F(q,u)=&-u\cdot\n q,\\[2mm]
G(q,u)=&-u.\n u+2\mu\n q\cdot D u+\lambda{\rm div}uÊ\,\n q+\frac{\kappa}{2}\n(|\n q|^{2}).
\end{aligned}
\label{f1a3}
\end{equation}
We are going to check that we can apply a fixed point theorem for the function $\psi$ in $E^{\N}$ defined below, the proof is divided in two step the stability of $\psi$ for a ball $B(0,R)$ in $E^{\N}$ and the contraction property. We define $E^{\N}$ by:
$$E^{\N}=\big(\widetilde{L}^{\infty}(\widetilde{B}^{\N-1,\N}_{2,\infty})\cap \widetilde{L}^{1}(\widetilde{B}^{\N+1,\N+2}_{2,\infty})\big)\times\big(\widetilde{L}^{\infty}(B^{\N-1}_{2,\infty})\cap \widetilde{L}^{1}(B^{\N+1}_{2,\infty})\big)^{N}.$$
\subsubsection*{1) First step, stability of $B(0,R)$:}
Let:
$$\eta=\| q_{0}\|_{\widetilde{B}^{\N-1,\N}_{2,\infty}}+\|u_{0}\|_{B^{\N-1}_{2,\infty}}\;.$$
We are going to show that $\psi$  maps the ball $B(0,R)$ into itself if $R$ is small enough.
According to proposition \ref{1flinear3}, we have:
\begin{equation}
\|W(t,\cdot)*\left(\begin{array}{c}
q_{0}\\
u_{0}\\
\end{array}
\right)\|_{E^{\N}}\leq C\eta \;.
\label{f1a5}
\end{equation}
According to the proposition \ref{1flinear3} it implies also that:
\begin{equation}
\begin{aligned}
\|\psi(q,u)\|_{E^{\N}}\leq&\;C\big(\eta+
\|F(q,u)\|_{ \widetilde{L}^{1}(\widetilde{B}^{\N-1,\N}_{2,\infty})}+\|G(q,u)\|_{ \widetilde{L}^{1}(B^{\N-1}_{2,\infty})}\big).
\end{aligned}
\label{f1a6}
\end{equation}
The main task consists in using the propositions \ref{produit} and corollary \ref{produit2} to obtain estimates on
$$\|F(q,u)\|_{ \widetilde{L}^{1}(\widetilde{B}^{\N-1,\N}_{2,\infty})},\;\|G(q,u)\|_{ \widetilde{L}^{1}(B^{\N-1}_{2,\infty})}.$$
Let us first estimate $\|F(q,u)\|_{ \widetilde{L}^{1}(\widetilde{B}^{\N-1,\N}_{2,\infty})}$. According to proposition
\ref{produit},
we have:
\begin{equation}
\begin{aligned}
&\|u\cdot\n q\|_{
\widetilde{L}^{1}(\widetilde{B}^{\N-1,\N}_{2,\infty})}\lesssim \|\n q\|_{\widetilde{L}^{\frac{4}{3}}(\widetilde{B}^{\N-\frac{1}{2},\N+\frac{1}{2}}_{2,\infty})} \|u\|_{\widetilde{L}^{4}(B^{\N-\frac{1}{2}}_{2,\infty})}\\
&\hspace{5cm}
+\|\n q\|_{\widetilde{L}^{4}(\widetilde{B}^{\N-\frac{3}{2},\N-\frac{1}{2}}_{2,\infty})} \|u\|_{\widetilde{L}^{\frac{4}{3}}(B^{\N+\frac{1}{2}}_{2,\infty})}.
\end{aligned}
\label{estim1}
\end{equation}
%Because $\widetilde{B}^{\N,\N+1}\h B^{\N}$ and
%$\widetilde{B}^{\N,\N+1}\h B^{\N+1}$ (from proposition
%\ref{fimportant}), we get:
%$$\|{\rm div}(qu)\|_{L^{1}(\widetilde{B}^{\N-1,\N})}\leq\|q\|_{L^{\infty}(\widetilde{B}^{\N,\N+1})}
%\|u\|_{L^{1}(B^{\N+1})}+\|q\|_{L^{2}(\widetilde{B}^{\N,\N+1})}\|u\|_{L^{2}(B^{\N})}.$$
%\\
%2) We have  to 
Let us now estimate $\|G(q,u)\|_{\widetilde{L}^{1}(B^{\N-1}_{2,\infty})}$. Hence by proposition \ref{produit}
it yields:
\begin{equation}
\begin{aligned}
\|u\cdot \n u\|
_{\widetilde{L}^{1}(B^{\N-1}_{2,\infty})}\lesssim&\|u\|_{\widetilde{L}^{\frac{4}{3}}(B^{\N+\frac{1}{2}}_{2,\infty})} \|u\|_{\widetilde{L}^{4}_{T}(B^{\N-\frac{1}{2}}_{2,\infty})}.
\end{aligned}
\label{estim2}
\end{equation}
In the same way we have by proposition \ref{produit} and the fact that $\widetilde{B}^{\N-\frac{1}{2},\N+\frac{1}{2}}_{2,\infty}\h B^{\N+\frac{1}{2}}_{2,\infty}$,  $\widetilde{B}^{\N-\frac{3}{2},\N-\frac{1}{2}}_{2,\infty}\h B^{\N-\frac{1}{2}}_{2,\infty}$ (see the remark \ref{r13}):
\begin{equation}
\begin{aligned}
&\|\n q\cdot D u\|
_{\widetilde{L}^{1}(B^{\N-1}_{2,\infty})}\lesssim \|\n q\|_{\widetilde{L}^{\frac{4}{3}}(\widetilde{B}^{\N-\frac{1}{2},\N+\frac{1}{2}}_{2,\infty})} \|D u\|_{\widetilde{L}^{4}(B^{\N-\frac{3}{2}}_{2,\infty})}\\
&\hspace{5cm}
+\|\n q\|_{\widetilde{L}^{4}(\widetilde{B}^{\N-\frac{3}{2},\N-\frac{1}{2}}_{2,\infty})} \|D u\|_{\widetilde{L}^{\frac{4}{3}}(B^{\N-\frac{1}{2}}_{2,\infty})}.\\[2mm]
&\|\n q\,{\rm div}u\|
_{\widetilde{L}^{1}(B^{\N-1}_{2,\infty})}\lesssim \|\n q\|_{\widetilde{L}^{\frac{4}{3}}(\widetilde{B}^{\N-\frac{1}{2},\N+\frac{1}{2}}_{2,\infty})} \|{\rm div} u\|_{\widetilde{L}^{4}(B^{\N-\frac{3}{2}}_{2,\infty})}\\
&\hspace{5cm}
+\|\n q\|_{\widetilde{L}^{4}(\widetilde{B}^{\N-\frac{3}{2},\N-\frac{1}{2}}_{2,\infty})} \|{\rm div} u\|_{\widetilde{L}^{\frac{4}{3}}(B^{\N-\frac{1}{2}}_{2,\infty})}.
\end{aligned}
\label{estim3}
\end{equation}
It now remains only to deal with the capillary terms:
\begin{equation}
\begin{aligned}
\|\n(|\n q|^{2})\|
_{\widetilde{L}^{1}(B^{\N-1}_{2,\infty})}\lesssim &\||\n q|^{2}\|
_{\widetilde{L}^{1}(B^{\N}_{2,\infty})},\\
\lesssim&  \|\n q\|_{\widetilde{L}^{\frac{4}{3}}(\widetilde{B}^{\N-\frac{1}{2},\N+\frac{1}{2}}_{2,\infty})} \|\n q\|_{\widetilde{L}^{4}(\widetilde{B}^{\N-\frac{3}{2},\N-\frac{1}{2}}_{2,\infty})}.
\end{aligned}
\label{estim4}
\end{equation}
We have previously used the fact that $\widetilde{B}^{\N-\frac{1}{2},\N+\frac{1}{2}}_{2,\infty}\h B^{\N+\frac{1}{2}}_{2,\infty}$ and $\widetilde{B}^{\N-\frac{3}{2},\N-\frac{1}{2}}_{2,\infty}\h B^{\N-\frac{1}{2}}_{2,\infty}$.
We are now going to assume that $(q,u)$ belongs in the ball $B(0,R)$ of $E^{\N}$ with $R>0$. Combining the estimates (\ref{estim1}), (\ref{estim2}), (\ref{estim3}) and  (\ref{estim4}) 
we get:
\begin{equation}
\|\psi(q,u)\|_{E^{\N}}\leq C((C+1)\eta+R)^{2}. \label{f1a7}
\end{equation}
By choosing $R$ and $\eta$ small enough we have:
\begin{equation}
C((C+1)\eta+R)^{2}\leq R.
\label{Rcrucial}
\end{equation}
%such that:
%$$R\leq \inf((3C)^{-1},c,1),\;\mbox{and}\;\;\eta\leq \frac{\inf(R,c)}{C+1}\;.$$
It implies that the ball $B(0,R)$ of $E^{\N}$ is stable under $\psi$, indeed we have:
$$\psi(B(0,R))\subset B(0,R)\;,$$
\subsubsection*{2) Second step: Property of contraction}
We consider $(q_{1},u_{1}),\,
(q_{2},u_{2})$ in $B(0,R)$ and we are interested in verifying that $\psi$ is a contraction.
According to the proposition \ref{1flinear3} we have:
\begin{equation}
\begin{aligned}
&\|\psi(q_{2},u_{2})-\psi(q_{1},u_{1})\|_{E^{\N}}\leq C\big(
\|F(q_{2},u_{2})-F(q_{1},u_{1})\|_{\widetilde{L}^{1}(\widetilde{B}^{\N-1,\N}_{2,\infty})}\\
&\hspace{6cm}+\|G(q_{2},u_{2})-G(q_{1},u_{1})\|_{\widetilde{L}^{1}(B^{\N-1}_{2,\infty})}\big).
\end{aligned}
\label{f1a8}
\end{equation}
We set:
$$\delta q=q_{2}-q_{1}\;\;\;\mbox{and}\;\;\;\delta u=u_{2}-u_{1}.$$
We have:
$$
\begin{aligned}
&F(q_{2},u_{2})-F(q_{1},u_{1})=-\delta u\cdot\n q _{2}-u_{1}\cdot\n\delta q.\\[2mm]
&G(q_{2},u_{2})-G(q_{1},u_{1})=-u_{2}\cdot\n\delta u-\delta u \cdot\n u_{1}+\mu\n q_{2}\cdot D\delta u+\mu\n\delta q\cdot D u_{1}\\
&\hspace{3cm}+\lambda\n q_{2} {\rm div}\delta u+\lambda \n\delta q {\rm div} u_{1}+\n\big(\n q_{2}\cdot\n \delta q+\n  \delta q\cdot\n q_{1}\big).
\end{aligned}
$$
Let us first estimate $\|F(q_{2},u_{2})-F(q_{1},u_{1})\|_{\widetilde{L}^{1}(\widetilde{B}^{\N-1,\N}_{2,\infty})}$. We have by proposition \ref{produit} and the remark \label{r13}:
\begin{equation}
\begin{aligned}
&\|F(q_{2},u_{2})-F(q_{1},u_{1})\|_{\widetilde{L}^{1}(\widetilde{B}^{\N-1,\N}_{2,\infty})}\lesssim
 \|\delta u\|_{\widetilde{L}^{\frac{4}{3}}_{T}(B^{\N+\frac{1}{2}}_{2,\infty})}\|\n q_{2}\|_{\widetilde{L}^{4}_{T}(\widetilde{B}^{\N-\frac{3}{2},\N-\frac{1}{2}}_{2,\infty})}\\
&\hspace{0,5cm}+ \|\delta u\|_{\widetilde{L}^{\infty}_{T}(B^{\N-1}_{2,\infty})}\|\n q_{2}\|_{\widetilde{L}^{1}_{T}(\widetilde{B}^{\N,\N+1}_{2,\infty})}+\|u_{1}\|_{\widetilde{L}^{1}_{T}(B^{\N+1}_{2,\infty})}\|\n \delta q\|_{\widetilde{L}^{\infty}_{T}(\widetilde{B}^{\N-2,\N-1}_{2,\infty})}\\
&\hspace{6,5cm}+\|\n \delta q\|_{\widetilde{L}^{\frac{4}{3}}_{T}(\widetilde{B}^{\N-\frac{1}{2},\N+\frac{1}{2}}_{2,\infty})}\|u_{1}\|_{\widetilde{L}^{4}_{T}(B^{\N-\frac{1}{2}}_{2,\infty})}\\[1,5mm]
&\leq C \big(2\|(q_{2},u_{2})\|_{E^{\N}}+2\|(q_{1},u_{1})\|_{E^{\N}}\big)\|(\delta q,\delta u)\|_{E{\N}}.
\end{aligned}
\label{acFn}
\end{equation}
\\
Next, we have to bound $\|G(q_{2},u_{2})-G(q_{1},u_{1})\|_{\widetilde{L}^{1}(B^{\N-1}_{2,\infty})}$. We treat
only one typical term, the others are of the same form.\\
\begin{equation}
\begin{aligned}
&\|-u_{2}\cdot\n\delta u-\delta u \cdot\n u_{1}\|_{\widetilde{L}^{1}(B^{\N-1}_{2,\infty})}\lesssim
 \|\n\delta u\|_{\widetilde{L}^{\frac{4}{3}}_{T}(B^{\N-\frac{1}{2}}_{2,\infty})}\|u_{2}\|_{\widetilde{L}^{4}_{T}(B^{\N-\frac{1}{2}}_{2,\infty})}\\
&\hspace{0,5cm}+ \|\n \delta u\|_{\widetilde{L}^{\infty}_{T}(B^{\N-2}_{2,\infty})}\|u_{2}\|_{\widetilde{L}^{1}_{T}(B^{\N+1}_{2,\infty})}+\|\n u_{1}\|_{\widetilde{L}^{1}_{T}(B^{\N}_{2,\infty})}\| \delta u\|_{\widetilde{L}^{\infty}_{T}(B^{\N-1}_{2,\infty})}\\
&\hspace{6,5cm}+\| \delta u\|_{\widetilde{L}^{\frac{4}{3}}_{T}(B^{\N+\frac{1}{2}}_{2,\infty})}\|\n u_{1}\|_{\widetilde{L}^{4}_{T}(B^{\N-\frac{3}{2}}_{2,\infty})}\\[1,5mm]
&\leq C \big(2\|(q_{2},u_{2})\|_{E^{\N}}+2\|(q_{1},u_{1})\|_{E^{\N}}\big)\|(\delta q,\delta u)\|_{E{\N}}.
\end{aligned}
\label{acFn1}
\end{equation}
We can bound the other terms of  $\|G(q_{2},u_{2})-G(q_{1},u_{1})\|_{\widetilde{L}^{1}(B^{\N-1}_{2,\infty})}$. in a similar way and this work  is left to the reader.
Finally by combining (\ref{f1a8}), (\ref{acFn}) and  (\ref{acFn1}) we obtain for $C>0$ large enough:
$$
\begin{aligned}
\|\psi(q_{2},u_{2})-
\psi(q_{1},u_{1})\|_{E^{\N}} \leq
&C\,\|(\de q,\de u)\|_{E^{\N}}\,\big(\|(q_{1},u_{1})\|_{E^{\N}}+\|(q_{2},u_{2})\|_{E^{\N}}\|_{E^{\N}}\big).
\end{aligned}
$$
If one chooses $R$ small enough such that $RC\leq\frac{3}{4}$, we end up with using the previous estimate which yields:
$$\|\psi(q_{2},u_{2})-\Psi(q_{1},u_{1})\|_{E^{\N}}
\leq\frac{3}{4}\,\|(\de q,\de u)\|_{E^{\N}}.$$ We thus
have the property of contraction and so by the fixed point theorem,
we have the existence of a global solution $(q,u)$ to the system $(NHV1)$. Indeed we can see easily that
$E^{\N}$ is a Banach space.\\
Concerning the uniqueness of this solution, it suffices to apply the arguments of the third step of the result of uniqueness in finite time (see the previous section). More precisely if $(q,u)$ is the previous solution and $(q_{1},u_{1})$ an other solution in $E^{\N}$ then by setting $\delta q=q-q_{1}$ and $\delta u=u-u_{1}$ we show that for $T$ small enough:
$$\de q=0\;\;\;\mbox{and}\;\;\;\de u=0.$$
 on $[0,T]$ and we conclude after by
connectivity to get the uniqueness on $\R^{+}$.
\hfill {$\Box$}
\subsection{Global solution near equilibrium for system (\ref{3systeme})}
\label{vraisysteme}
We are now interested in proving the existence of global strong solution for the original system (\ref{3systeme}), as explained in the remarka \ref{remcrucim} the main difficulty consists in observing that a priori the system (\ref{3systeme}) is not equivalent to the system (\ref{NHV1}). Indeed we need at least to control the $L^{\infty}$ norm on the density $\rho$ and on the vacuum $\frac{1}{\rho}$ in order to propagate on the density $\rho$ the regularity proved in theorem \ref{ftheo1} for the unknown $q=\ln\rho$. Furthermore we shall easily verify that with such regularity on $(\rho,u)$ then $(\rho,u)$ verify the system (\ref{3systeme}) and is a unique solution of (\ref{3systeme}). To do this we are going to assume additional hypothesis on $(q_{0},u_{0})$ since now $(q_{0},u_{0})$ is in $\widetilde{B}^{\N-1,\N}_{2,2}$ with $q_{0}=\ln\rho_{0}\in L^{\infty}$ (this last condition implies in particular that $\rho_{0}$ and $\frac{1}{\rho_{0}}$ belong in $L^{\infty}$).\\
The first part of the proof consists in getting $L^{\infty}$ estimates for the linear system $(N)$, and the second part consists in splitting the solution $(q,u)$ under the following form:
$$(q,u)=(q_{L},u_{L})+(\bar{q},\bar{u}),$$
with $(q_{L},u_{L})$ solution of $(N)$ with $F=G=0$ and $(q_{L}(0,\cdot),u_{L}(0,\cdot)=(q_{0},u_{0})$. The key points will be to show that $q_{L}$ is belonging in $L^{\infty}_{T}(L^{\infty})$ for any $T>0$ and that $\bar{q}$ is more regular that $q_{L}$. More precisely by a regularizing effect on the third index of the Besov space we shall prove that $\bar{q}$ is in $\widetilde{L}_{T}^{\infty}(B^{\N}_{2,1})$ for any $T>0$ with is embedded in $L^{\infty}_{T}(L^{\infty})$. It will then be sufficient to deduce a control of $\ln\rho$ in  $L^{\infty}_{T}(L^{\infty})$ and to propagate the regularity of $q$ on $\rho$ via the proposition \ref{composition}.
%is a strong solution of system (\ref{3systeme}).
\subsubsection{Result of maximum principle type for the linear system $(N)$}
\label{submaximum}
Let us start with studying the following system:
\begin{equation}
\begin{cases}
\begin{aligned}
&\p_{t}q+{\rm div}u=0,\\
&\p_{t}u-\mu\D u-(\lambda+\mu)\n{\rm div}u-\kappa\n\D q=0,\\
&(q(0,\cdot),u(0,\cdot))=(q_{0},u_{0}),
\end{aligned}
\end{cases}
\label{systemelin}
\end{equation}
with $\mu>0$ and $\lambda+2\mu>0$. We are interesting in characterizing the Besov space in term of the semi group $B(t)$ associated to the system (\ref{systemelin}), it shall be useful in order to obtain $L^{\infty}$ estimates for $q$. More precisely we have the following proposition.
\begin{proposition}
Let $s$ be a positive real number and $(p,r)\in[1,+\infty]^{2}$. Let $(q,u)$ the solution of (\ref{systemelin}) with $(q,u)(t)=e^{B(t)}(q_{0},u_{0})$ and with the following notation:
$$(\n q,u)(t)=e^{B(t)}(\n q_{0},u_{0}).$$
Then there exists a constant $C>0$ which satisfies:
%\begin{equation}
%C^{-1}\|(\n q, u)\|_{B^{-2s}_{p,r}}\leq \|\;\|t^{s}e^{B(t)}u\|_{L^{p}}\;\|_{L^{r}(\R^{+},\frac{dt}{t})}\leq C\|(\n q, u)\|_{B^{-2s}_{p,r}}\;\;\;\forall (\n q, u)\in B^{-2s}_{p,r}. 
%\label{propimpin}
%\end{equation}
\begin{equation}
 \|\;\|t^{s}e^{B(t)}(\n q_{0},u_{0})\|_{L^{p}}\;\|_{L^{r}(\R^{+},\frac{dt}{t})}\leq C\|(\n q_{0}, u_{0})\|_{B^{-2s}_{p,r}}\;\;\;\forall (\n q_{0}, u_{0})\in B^{-2s}_{p,r}. 
\label{propimpin}
\end{equation}
\label{propimp}
\end{proposition}
{\bf Proof} Apply operator $\D$ to the first equation of (\ref{systemelin}) and operators ${\rm div}$ and ${\rm curl}$ to the second one; we obtain the following system
with $\nu=2\mu+\lambda$:
\begin{equation}
\begin{cases}
\begin{aligned}
&\p_{t}\D q+\D {\rm div}u=0,\\
&\p_{t}{\rm div}u-\nu\D {\rm div}u-\kappa\D^{2} q=0,\\
&\p_{t}{\rm curl}u-\mu\D {\rm curl}u=0,\\
&(q(0,\cdot),u(0,\cdot))=(q_{0},u_{0}).
\end{aligned}
\end{cases}
\label{65a}
\end{equation}
We observe that the third equation is a heat equation and we know via lemma 2.4 p 54 in \cite{BCD} that it exists $C,c>0$ such that:
\begin{equation}
\|e^{\mu t\D}\D_{l}{\rm curl}u_{0}\|_{L^{p}}\leq Ce^{-c2^{2l}\mu}\|\D_{l}{\rm curl}u_{0}\|_{L^{p}}\;\;\forall p\in[1,+\infty].
\label{curl}
\end{equation}
Let us study now the following system:
\begin{equation}
\begin{cases}
\begin{aligned}
&\p_{t}c+\D v=0,\\
&\p_{t}v-\mu\D v-\kappa\D c=0,\\
%&\p_{t}{\rm curl}u-\mu\D {\rm curl}u=0,\\
&(c(0,\cdot),v(0,\cdot))=(c_{0},v_{0}).
\end{aligned}
\end{cases}
\label{65}
\end{equation}
where $c=\D q$ and $v={\rm div}u$.  Denoting by $U(t)$ the semi-group associated to (\ref{65}), we deduce from Duhamel's formula that:
 $$
\p_{t}\left(\begin{array}{c}
\hat{q}(t,\xi)\\
\hat{d}(t,\xi)\\
%\hat{{\cal{T}}}(t,\xi)\\
\end{array}
\right)=U(t)
\left(\begin{array}{c}
q_{0}\\
u_{0}\\
%\hat{{\cal{T}}}(t,\xi)\\
\end{array}
\right)
$$
We are now interested in proving estimate of the same form than (\ref{curl}) for $e^{U(t)}(\D_{l}c,\D_{l}v)$, and to do this we are going to prove the following lemma which is a direct consequence of the lemma 3 of \cite{fDD}. For the sake of completeness we are going to recall its proof.
\begin{lemme}
For any $p\in[1,+\infty]$ and any $t>0$ we have:
\begin{equation}
\|e^{U(t)}(\D_{q}c_{0},\D_{q}v_{0})\|_{L^{p}}\leq  C e^{-c\min(1,\frac{4\kappa}{\nu^{2}})2^{2q}\nu t}(\|\D_{q}c_{0}\|_{L^{p}}+\|\D_{q}v_{0}\|_{L^{p}}).
\label{div}
\end{equation}
\label{tech}
\end{lemme}
%A(\xi)\left(\begin{array}{c}
%\hat{q}(t,\xi)\\
%\hat{d}(t,\xi)\\
%\hat{{\cal{T}}}(t,\xi)\\
%\end{array}
%\right)=\left(\begin{array}{c}
%\hat{F}(t,\xi)\\
%\Lambda^{-1}{\rm div}\,\hat{G}(t,\xi)\\
%\hat{H}(t,\xi)\\
%\end{array}
%\right) \leqno{(M^{'}_{2})}
%$$
{\bf Proof:} Simple calculus show that $U(t)=e^{-t A(D)}$ with:
$$A(\xi)=\left(\begin{array}{ccc}
0&-|\xi|^{2}\\
\kappa|\xi|^{2} &\mu|\xi|^{2}\\
\end{array}
\right)\,.
$$
%\begin{proposition}
%\label{defBesov}
%\end{proposition}
Following \cite{fDD} we show that:
$$e^{-t A(\xi)}=e^{-\frac{t\nu|\xi|^{2}}{2}}\left(\begin{array}{ccc}
h_{1}(t,\xi)+\frac{\nu}{2}h_{2}(t,\xi)&h_{2}(t,\xi)\\
-\kappa h_{2}(t,\xi)&h_{1}(t,\xi)-\frac{\nu}{2}h_{2}(t,\xi)\\
\end{array}
\right)\,,
$$
with:
$$
\begin{aligned}
&h_{1}(t,\xi)=\cos (\nu^{'}|\xi|^{2} t),\;\;h_{2}(t,\xi)=\frac{\sin (\nu^{'}|\xi|^{2} t)}{\nu^{'}},\mbox{if}\;\;\mu^{2}<4\kappa.\\
&h_{1}(t,\xi)=1,\;\;h_{2}(t,\xi)=t|\xi|^{2};\mbox{if}\;\;\nu^{2}=4\kappa.\\
&h_{1}(t,\xi)=\cos (\nu^{'}|\xi|^{2} t),\;\;h_{2}(t,\xi)=\frac{\sin (\nu^{'}|\xi|^{2} t)}{\mu^{'}},\mbox{if}\;\;\nu^{2}>4\kappa.
\end{aligned}
$$
and $\nu^{'}=\sqrt{|\kappa-\frac{\nu^{2}}{4}|}$.
Let $\varphi$ defined as in the definition of Littlewood-Paley theory, then denoting $a_{ij}(t,\xi)$
the coefficients of the matrix $e^{-tA(\xi)}$ and:
$$\D_{q}b_{ij}(t,x)=\D_{q}\big({\cal F}^{-1}(a_{ij})(t,x)\big)=(2\pi)^{-N}\int e^{ix\cdot\xi}a_{ij}(t,\xi)\varphi(2^{-q}\xi)d\xi.$$
Let us show now that:
\begin{equation}
\|\D_{q}(b_{ij})(t,\cdot)\|_{L^{1}}\leq C e^{-c\min (1,\frac{4\kappa}{\nu^{2}})2^{2q}\nu t}.
\label{66}
\end{equation}
where $c$ depends only on $\mu$, $\kappa$ and $c$ is a universal constant. We first remark that $\|\D_{q}b_{ij}\|_{L^{1}}=\|h_{ijq}\|_{L^{1}}$ with:
$$h_{ijq}(t,y)=(2\pi)^{-N}\int e^{iy\cdot\eta}a_{ij}(t,2^{q}\eta)\varphi(\eta)d\eta.$$
Let us observe that the functions $h_{ijq}$ can be rewritten under the form:
\begin{equation}
h_{q}(t,x)=\int e^{iy\cdot\xi}f(2^{q}|xi|^{2}t)\varphi(\xi)d\xi.
\label{67}
\end{equation}
with $f\in C^{\infty}(\R^{+})$. By integration by parts and Leibniz' formula, we obtain for all $\alpha\in\mathbb{N}^{N}$,
\begin{equation}
(-i x)^{\alpha}h_{q}(x)=\sum_{\beta\leq\alpha}(^{\alpha}_{\beta})\int e^{ix\cdot\xi}\p^{\beta}f(2^{q}|xi|^{2}t)\p^{\alpha-\beta}\varphi(\xi)d\xi.
\label{68}
\end{equation}
Next, from Fa\`a-di-Bruno's formula, we deduce that:
\begin{equation}
\p^{\beta}f(2^{q}|xi|^{2}t)=\sum_{\gamma_{1}+\cdots+\gamma_{m}=\beta, |\gamma_{i}|\geq 1}f^{(m)}(2^{q}|xi|^{2}t)(2^{2q}t)^{m}(\\Pi^{m}_{j=1}\p^{\gamma_{j}}(|xi|^{2})).
\label{69}
\end{equation}
Let us suppose first that $\nu^{2}>4\kappa$. Then, we just have to prove that:
\begin{equation}
\|h_{q}\|_{L^{1}}\leq C e^{-c2^{2q}\nu t},
\label{70}
\end{equation}
for $f(u)=e^{i\mu^{'}u}e^{-\mu u/2}$. We have:
$$f^{(m)}(u)=(i\nu^{'}-\frac{\nu}{2})^{m}e^{i\nu^{'}u}e^{-\nu u/2},$$
so that $|f^{(m)}(u)|\leq (\nu^{'}+\frac{\nu}{2})^{m}e^{-\mu u/2}$. Using (\ref{68}), (\ref{69}), we prove the existence of constant $C_{\alpha,\beta,m}$ such that:
$$|x^{\alpha}h_{q}(x)|\leq \sum_{\beta\leq\alpha}\sum^{|\beta|}_{m=1}C_{\alpha,\beta,m}(2^{2q}t)^{m}e^{-\nu t2^{2q}/8}.$$
For any constant $c<1$ and $m\in\mathbb{N}$, there exists $C_{m}$ such that $u^{m}e^{-u}\leq C_{m}e^{-cu}$. This clearly yields (\ref{70}).\\
When $\nu^{2}=4\kappa$, it suffices to verify (\ref{70}) for $f(u)=u e^{-\nu u/2}$ and $f(u)=e^{-\nu u/2}$. This is a direct consequence of (\ref{69}) and Leibniz' formula. When $\nu^{2}>4\kappa$, we have to check (\ref{70}) for:
$$f(u)=\exp(Ð\frac{\nu}{2}(1+\sqrt{1-\frac{4\kappa}{\nu^{2}}})u)\;\;\mbox{and}\;\;f(u)=\exp(Ð\frac{\nu}{2}(1-\sqrt{1-\frac{4\kappa}{\nu^{2}}})u).$$
Using again (\ref{69}) we have:
$$|x^{\alpha}h_{q}(x)|\leq C \max(e^{-c\nu t 2^{2q}(1+\sqrt{1-\frac{4\kappa}{\nu^{2}} })},e^{-c\nu t 2^{2q}(1-\sqrt{1-\frac{4\kappa}{\nu^{2}} })})\leq C e^{-c(\frac{\kappa}{\nu})2^{2q}t}$$
and we conclude to (\ref{66}).\\
We obtain finally by using (\ref{66}) and the Young inequality for the convolution:
$$\|e^{U(t)}(\D_{q}c_{0},\D_{q}v_{0})\|_{L^{p}}\leq  C e^{-c\min(1,\frac{4\kappa}{\nu^{2}})2^{2q}\mu t}(\|\D_{q}c_{0}\|_{L^{p}}+\|\D_{q}v_{0}\|_{L^{p}}).$$
It proves the lemma \ref{tech}.  \hfill {$\Box$}\\
\\
Via the Bernstein lemma, the lemma \ref{singuliere}, (\ref{curl}) and (\ref{div}) we obtain that $\forall l\in\mathbb{Z}$:
\begin{equation}
\|e^{B(t)}(\D_{l}\n q_{0},\D_{l}u_{0})\|_{L^{p}}\leq  C e^{-c\min(1,\frac{4\kappa}{\nu^{2}})2^{2l}\nu t}(\|\D_{l}\n q_{0}\|_{L^{p}}+\|\D_{l}u_{0}\|_{L^{p}}).
\label{crucial}
\end{equation}
According the estimate (\ref{crucial}) we have by setting $V=(\n q_{0},u_{0})$ for $c,C>0$:
$$\|t^{s}\D_{l}e^{B(t)}V\|_{L^{p}}\leq Ct^{s}2^{2ls}e^{-ct 2^{2l}}2^{-2ls}\|\D_{l}V\|_{L^{p}}.$$
We are now going to adapt a criterion to define some Besov space in terms of the heat kernel (see theorem 2.34 in \cite{BCD}) to our case, it means to the semi-group $B(t)$.\\
Since $V$ belongs in ${\cal S}^{'}_{h}$ and the definition of the homogeneous Besov semi norm we have:
$$
\begin{aligned}
\|t^{s}e^{tB(t)}V\|_{L^{p}}&\leq\sum_{l\in\mathbb{Z}}\|t^{s}\D_{l}e^{B(t)}V\|_{L^{p}},\\
&\leq C\|V\|_{B^{-2s}_{p,r}}\sum_{l\in\mathbb{Z}}t^{s}2^{2ls}e^{-ct 2^{2l}}c_{r,l}
\end{aligned}
$$
where $(c_{rl})_{l\in\mathbb{Z}}$ is a element of the unit sphere of $l^{r}(\mathbb{Z})$. If $r=+\infty$, we easily show (\ref{div}) by using the following lemma (we left to the reader the proof of this last one).
\begin{lem}
For any $s$, we have:
$$\sup_{t>0} \sum_{l\in\mathbb{Z}}t^{s}2^{2ls}e^{-ct 2^{2l}}<+\infty.$$
\label{lemmepetittech}
\end{lem}
Let us deal now with the case  $r<+\infty$, combining H\" older's inequality with the weight $2^{2ls}e^{-ct 2^{2l}}$ and lemma \ref{lemmepetittech}, we obtain:
$$
\begin{aligned}
\int^{+\infty}_{0}t^{rs}\|e^{B(t)}V\|_{L^{p}}^{r}\frac{dt}{t}&\leq C\|V\|^{r}_{B^{-2s}_{p,r}}\int^{+\infty}_{0}(\sum_{l\in\mathbb{Z}}t^{s}2^{2ls}e^{-ct 2^{2l}}c_{rl})^{r}\frac{dt}{t},\\
&\leq C\|V\|^{r}_{B^{-2s}_{p,r}}\int^{+\infty}_{0}(\sum_{l\in\mathbb{Z}}t^{s}2^{2ls}e^{-ct 2^{2l}})^{r-1}(\sum_{l\in\mathbb{Z}}t^{s}2^{2ls}e^{-ct 2^{2l}}c_{rl}^{r})
\frac{dt}{t},\\
&\leq C\|V\|^{r}_{B^{-2s}_{p,r}}\int^{+\infty}_{0}(\sum_{l\in\mathbb{Z}}t^{s}2^{2ls}e^{-ct 2^{2l}}c_{rl}^{r})\frac{dt}{t}
\end{aligned}
$$
Using Fubini's theorem and the change of variable $u=ct2^{j}$, we have:
$$
\begin{aligned}
\int^{+\infty}_{0}t^{rs}\|e^{B(t)}V\|_{L^{p}}^{r}\frac{dt}{t}&\leq  C\|V\|^{r}_{B^{-2s}_{p,r}}\sum_{l\in\mathbb{Z}}c_{rl}^{r} \int^{+\infty}_{0}t^{s}2^{2ls}e^{-ct 2^{2l}}\frac{dt}{t},\\
&\leq C\Gamma(s)\|V\|^{r}_{B^{-2s}_{p,r}},
\end{aligned}
$$
with $\Gamma(s)=\int^{+\infty}_{0}t^{s-1}e^{-t}dt$. %It achieves the proof of the right inequality in (\ref{propimpin}).\\[2mm]
The proposition \ref{propimp} is thus proved. \hfill {$\Box$}\\
\\
Let us prove now $L^{\infty}$ estimate for the density solution $q$ of the system (\ref{systemelin}). We recall that in the sequel $q$ will correspond roughly speaking to $\ln\rho$, in particular $L^{\infty}$ estimate on $q$ shall also provide $L^{\infty}$ estimate on $\rho$ and $\frac{1}{\rho}$. %, more precisely we are dealing with the unknown $(q,(\D)^{-1}{\rm div}u)$.
\begin{proposition}
Let $q_{0}\in B^{\N}_{2,2}$, $u_{0}\in B^{N-1}_{2,2}$ and $q_{0}\in L^{\infty}$. Let $(q,u)$ the solution of the system \ref{systemelin}, then we have for all $T>0$:
\begin{equation}
\begin{aligned}
&\sup_{t\in[0,T]}(\sqrt{t}\|(\n q,u)(t,\cdot)\|_{L^{\infty}})\leq C\|(\n q_{0},u_{0})\|_{B^{\N-1}_{2,2}}.\\
&\|q\|_{L^{\infty}_{T}(L^{\infty})}\leq \|q_{0}\|_{L^{\infty}}+C\|(\n q_{0},u_{0})\|_{B^{\N-1}_{2,2}}.
\end{aligned}
\label{maxfinal}
\end{equation}
\label{propmaxfinal}
\end{proposition}
{\bf Proof:} The first estimate in (\ref{maxfinal}) is a direct application of the proposition \ref{propimp} applied to $p=+\infty$, $r=+\infty$, $s=\frac{1}{2}$ and using the fact that $B^{\N-1}_{2,2}$ is embedded in $B^{-1}_{\infty,\infty}$.\\
Let ${\cal E}$ the fundamental solution of the Laplacian operator, and we define the operator $(\D)^{-1}$ by the convolution operator $(\D)^{-1}f={\cal E}*f$ with $f\in{\cal S}^{'}(\R^{N})$. By applying the operator $(\D)^{-1}{\rm div}$ to the second equation of (\ref{systemelin}) and using the fact that $\D c={\rm div}u$, we obtain the following system with $c=(\D)^{-1}{\rm div}u$:
$$
\begin{cases}
\begin{aligned}
&\p_{t}q-\frac{\kappa}{\mu}\D q=-\frac{1}{\mu}\p_{t}c,\\
&\p_{t}c-\mu\D c-\kappa\D q=0.
\end{aligned}
\end{cases}
$$
Let us prove that $q
$ belongs in $L^{\infty}_{T}(L^{\infty})$ for any $T>0$; by Duhamel formula we have:
\begin{equation}
q(t,x)=e^{\frac{\kappa}{\mu}t\D}q_{0}-\frac{1}{\mu}\int^{t}_{0}e^{\frac{\kappa}{\mu}(t-s)\D}\p_{t}c(s)ds.
\label{Duhamel}
\end{equation}
By the maximum principle for the heat equation, we deduce that:
\begin{equation}
\|e^{\frac{\kappa}{\mu}t\D}q_{0}\|_{L^{\infty}(L^{\infty})}\leq \|q_{0}\|_{L^{\infty}(\R^{N})}.
\label{bmax}
\end{equation}
Next we are going to consider $\int^{t}_{0}e^{\frac{\kappa}{\mu}(t-s)\D}\p_{s}c(s)ds$, we recall that:
\begin{equation}
\begin{aligned}
e^{\frac{\kappa}{\mu}(t-s)\D}\p_{s}c(s)&=K(\frac{\cdot}{\sqrt{t-s}})*_{x}\p_{s}(\D)^{-1}{\rm div}u(s,\cdot),\\
&=K(\frac{\cdot}{\sqrt{t-s}})*_{x}\big(\p_{s}[({\cal E}*_{x}{\rm div}u(s,\cdot)]\big),\\
&=K(\frac{x}{\sqrt{t-s}})*_{x}\big([\sum_{i}\p_{i}{\cal E}*_{x}\p_{s}u_{i}(s,\cdot)]\big),\\
&=\sum_{i}\big(K(\frac{x}{\sqrt{t-s}})*_{x}\p_{i}{\cal E}\big)*_{x}\p_{s}u_{i}(s,\cdot)]\big),\\
&=\sum_{i}\big(\p_{i}[K(\frac{\cdot}{\sqrt{t-s}})]*_{x}{\cal E}\big)*_{x}\p_{s}u_{i}(s,\cdot)\\
&=\sum_{i}\p_{i}[K(\frac{\cdot}{\sqrt{t-s}})]*_{x}\big({\cal E}*_{x}\p_{s}u_{i}(s,\cdot)\big)\\
&=\sum_{i}\p_{i}[K(\frac{\cdot}{\sqrt{t-s}})]*_{x}(\D)^{-1}\p_{s}u_{i}(s,\cdot),
\end{aligned}
\label{petitcalc}
\end{equation}
with $*_{x}$ the convolution in space and setting $\bar{\mu}=\frac{\kappa}{\mu}$:
$$K(\frac{x}{\sqrt{t}})=\frac{1}{(4\pi \bar{\mu}  t)^{\N}}e^{-\frac{|x|^{2}}{4 \bar{\mu}t}}.$$
We deduce that:
$$\p_{i}[K(\frac{\cdot}{\sqrt{t-s}})](x)=\frac{-2x_{i}}{\pi^{\N}}\frac{1}{(4\bar{\mu}(t-s))^{\N+1}}e^{-\frac{|x|^{2}}{4 \bar{\mu}(t-s)}}$$
We easily check by a change of variable $u=\frac{x}{\sqrt{4 \bar{\mu}(t-s)}}$ that:
$$\|\p_{i}[K(\frac{\cdot}{\sqrt{t-s}})]\|_{L^{1}}\leq \frac{C}{\sqrt{t-s}}.$$
Then by Young's inequality we have for $0<s<t$:
\begin{equation}
\|\p_{i}[K(\frac{\cdot}{\sqrt{t-s}})]*_{x}(\D)^{-1}\p_{s}u_{i}(s,\cdot)\|_{L^{\infty}_{x}}\leq \frac{C}{\sqrt{t-s}}\frac{1}{\sqrt{s}}\sup_{0<s<t}\sqrt{s}\|(\D)^{-1}\p_{s}u(s)\|_{L^{\infty}}
\label{petitcalc2}
\end{equation}
%\leq \frac{C^{'} }{\sqrt{t-s}}\frac{1}{\sqrt{s}},$$
with $C>0$. By applying the operator $(\D)^{-1}$ to the second equation of (\ref{systemelin}) we show that:
\begin{equation}
\p_{t}(\D)^{-1}u=\mu u+\kappa\n q.
\label{petitcalc1}
\end{equation}
We conclude as follows by using (\ref{petitcalc}), (\ref{petitcalc2}) and (\ref{petitcalc1}), then there exists $C>0$ such that:
\begin{equation}
\begin{aligned}
\|\int^{t}_{0}e^{\frac{\kappa}{\mu}(t-s)\D}\p_{s}c(s)ds\|_{L^{\infty}(L^{\infty})}&\leq(\sup_{0<s<t}\sqrt{s}\|\D^{-1}\p_{s} u(s)\|_{L^{\infty}}) \int^{t}_{0}\frac{1}{\sqrt{t-s}}\frac{1}{\sqrt{s}}ds,\\
&\leq\pi \sup_{0<s<t}\sqrt{s}\|\mu u(s,\cdot)+\kappa\n q(s,\cdot)\|_{L^{\infty}},\\
&\leq C\pi \sup_{0<s<t}\sqrt{s}\|(u(s,\cdot),\n q(s,\cdot))\|_{L^{\infty}}.
\end{aligned}
\label{bv1}
\end{equation}
because:
$$\int^{t}_{0}\frac{1}{\sqrt{t-s}}\frac{1}{\sqrt{s}}ds=\pi.$$
By using the first estimate in (\ref{maxfinal}) we deduce that by using (\ref{bv1}) and (\ref{bmax}) we obtain that for any $T>0$, it exists $C>0$:
$$\|q\|_{L^{\infty}_{T}(L^{\infty})}\leq \|q_{0}\|_{L^{\infty}}+C\|(\n q_{0},u_{0})\|_{B^{\N-1}_{2,2}}.$$
It achieves the proof of the proposition \ref{maxfinal}.
 \hfill {$\Box$}
 \subsubsection{The unique solution of (\ref{NHV1}) verifies $\rho\in L^{\infty}_{T}(L^{\infty})$}
 \label{controldens}
 \subsubsection*{Regularizing effect on the third index of the Besov spaces}
 Let us start by recalling that we are concerned now with initial data such that $(q_{0},u_{0})$ belong in $(B^{\N}_{2,2}\cap L^{\infty})\times B^{\N-1}_{2,2}$. This additional regularity assumptions on the initial data will be crucial in order to prove $L^{\infty}$ estimates on $\ln\rho$. We are going to begin with proving additional regularity
 assumption on the solution $(q,u)$ of the system (\ref{NHV1}) which verify $(q,u)\in E$ with $E$ defined as follows:
 $$E=\big(\widetilde{L}^{\infty}(\R^{+},\widetilde{B}^{\N-1,\N}_{2,\infty})\cap \widetilde{L}^{1}(\R^{+},\widetilde{B}^{\N+1,\N+2}_{2,\infty}\big)\times
 \big(\widetilde{L}^{\infty}(\R^{+},B^{\N-1}_{2,\infty})\cap  \widetilde{L}^{1}(\R^{+},B^{\N+1}_{2,\infty})\big)^{N}.$$
Indeed we are interested in splitting the unique solution $(q,u)$ constructed in the subsection \ref{infini} as the following sum:
$$(q,u)=(q_{L},u_{L})+(\bar{q},\bar{u}).$$
with $(q_{L},u_{L})$ the solution of the system $(N1)$ with $F=G=0$ $a=\mu$, $b=\lambda+\mu$, $c=\kappa$, $d=K$ and with initial data $(\ln\rho_{0},u_{0})$. Compared with the subsection \ref{infini}, we are going to show that the remainder term $(\bar{q},\bar{u})$ is more regular than $(q_{L},u_{L})$ and is in the space $F$ defined as follows:
$$F=\big(\widetilde{L}^{\infty}(\R^{+},\widetilde{B}^{\N-1,\N}_{2,1})\cap \widetilde{L}^{1}(\R^{+},\widetilde{B}^{\N+1,\N+2}_{2,1})\big)\times
 \big(\widetilde{L}^{\infty}(\R^{+},B^{\N-1}_{2,1}\cap  \widetilde{L}^{1}(\R^{+},B^{\N+1}_{2,1})\big)^{N}.$$
This type of regularizing effect on the remainder $(\bar{q},\bar{u})$ in term of the third index on the Besov spaces has been observed for the first time by Cannone and Planchon in \cite{CP} for the incompressible Navier-Stokes equations. We are going to use similar ideas in your case than \cite{CP}, to do this let us use the proposition \ref{1flinear3} which ensures:
\begin{equation}
\begin{aligned}
&\|q_{L}\|_{\widetilde{L}^{\infty}(\R^{+},B^{\N-1}_{2,2}\cap B^{\N}_{2,2})}+\|q_{L}\|_{\widetilde{L}^{1}(\R^{+},B^{\N+1}_{2,2}\cap B^{\N+2}_{2,2})}
+\|u_{L}\|_{\widetilde{L}^{\infty}(\R^{+},B^{\N-1}_{2,2})}\\
&\hspace{5cm}+\|u_{L}\|_{\widetilde{L}^{1}(\R^{+},B^{\N+1}_{2,2})}\lesssim \|q_{0}\|_{\widetilde{B}^{\N-1,\N}_{2,2}}+\|u_{0}\|_{B^{\N-1}_{2,2}}.
\end{aligned}
\label{rinitial}
\end{equation}
%In particular by applying exactly the same arguments than in the subsection \ref{infini}, we can easily verify that it exists a global
As in the subsection \ref{infini}, we know that $(\bar{q},\bar{u})$ verify the following system:
\begin{equation}
\begin{cases}
\begin{aligned}
&\p_{t}\bar{q}+{\rm div}\bar{u}=F(\bar{q},\bar{u}),\\
&\p_{t}\bar{u}-\mu\D\bar{u}-(\lambda+\mu)\n{\rm div}\bar{u}-\kappa\n\D\bar{q}+K\n q=G(\bar{q},\bar{u}),\\
&(\bar{q},\bar{u})(0,\cdot)=(0,0),
\end{aligned}
\end{cases}
\label{sysreste}
\end{equation}
with:
$$
\begin{aligned}
&F(\bar{q},\bar{u})=-u\cdot\n q,\\
&G(\bar{q},\bar{u})=-u\cdot\n u+2\mu\n q\cdot Du+\lambda{\rm div}u\,\n q+\frac{\kappa}{2}\n (|\n q|^{2}).
\end{aligned}
$$
We are going to apply the proposition \ref{1flinear3} to $(\bar{q},\bar{u})$, it implies that:
\begin{equation}
\begin{aligned}
&\|\bar{q}\|_{\widetilde{L}^{\infty}(\R^{+},\widetilde{B}^{\N-1,\N}_{2,1})}+\|\bar{q}\|_{\widetilde{L}^{1}(\R^{+},\widetilde{B}^{\N+1,\N+2}_{2,1})}
+\|\bar{u}\|_{\widetilde{L}^{\infty}(\R^{+},B^{\N-1}_{2,1})}+\|\bar{u}\|_{\widetilde{L}^{1}(\R^{+},B^{\N+1}_{2,1})}\\
&\hspace{4cm}\lesssim \|F(\bar{q},\bar{u})\|_{\widetilde{L}^{1}(\R^{+},\widetilde{B}^{\N-1,\N}_{2,1})}+\|G(\bar{q},\bar{u})\|_{\widetilde{L}^{1}(\R^{+},B^{\N-1}_{2,1})}.
\end{aligned}
\label{rinitiala}
\end{equation}
It remains only to bound the terms $F(\bar{q},\bar{u})$ and $G(\bar{q},\bar{u})$ on the right hand side of (\ref{rinitiala}). Let us start with $F(\bar{q},\bar{u})$%; according to proposition \ref{produit} we have:% 
and $u\cdot\n q$ which that we rewrite under the following form:
$$u\cdot\n q=\bar{u}\cdot\n q+u_{L}\cdot\n\bar{q}+u_{L}\cdot\n q_{L}.$$
According to proposition \ref{produit} we have:
\begin{equation}
\begin{aligned}
&\|\bar{u}\cdot\n q\|_{\widetilde{L}^{1}_{T}(\widetilde{B}^{\N-1,\N}_{2,1})}\lesssim \|\n q
\|_{\widetilde{L}^{\frac{4}{3}}_{T}(\widetilde{B}^{\N-\frac{1}{2},\N+\frac{1}{2}}_{2,\infty})}\|\bar{u}
\|_{\widetilde{L}^{4}_{T}(B^{\N-\frac{1}{2}}_{2,1})}+\|\n q
\|_{\widetilde{L}^{4}_{T}(\widetilde{B}^{\N-\frac{3}{2},\N-\frac{1}{2}}_{2,\infty})}\\
&\hspace{10cm}\times\|\bar{u}
\|_{\widetilde{L}^{\frac{4}{3}}_{T}(B^{\N+\frac{1}{2}}_{2,1})}
,\\[2mm]
&\|u_{L}\cdot\n \bar{q}\|_{\widetilde{L}^{1}_{T}(B^{\N}_{2,1})}\lesssim \|\n \bar{q}
\|_{\widetilde{L}^{\frac{4}{3}}_{T}(\widetilde{B}^{\N-\frac{1}{2},\N+\frac{1}{2}}_{2,1})}\|u_{L}
\|_{\widetilde{L}^{4}_{T}(B^{\N-\frac{1}{2}}_{2,\infty})}+\|\n\bar{ q}
\|_{\widetilde{L}^{4}_{T}(\widetilde{B}^{\N-\frac{3}{2},\N-\frac{1}{2}}_{2,1})}\\
&\hspace{10cm}\times\|u_{L}
\|_{\widetilde{L}^{\frac{4}{3}}_{T}(B^{\N+\frac{1}{2}}_{2,\infty})}
,\\[2mm]
&\|u_{L}\cdot\n q_{L}\|_{\widetilde{L}^{1}_{T}(B^{\N}_{2,1})}\lesssim \|\n q_{L}
\|_{\widetilde{L}^{\frac{4}{3}}_{T}(\widetilde{B}^{\N-\frac{1}{2},\N+\frac{1}{2}}_{2,2})}\|u_{L}
\|_{\widetilde{L}^{4}_{T}(B^{\N-\frac{1}{2}}_{2,2})}+\|\n q_{L}
\|_{\widetilde{L}^{4}_{T}(\widetilde{B}^{\N-\frac{3}{2},\N-\frac{1}{2}}_{2,2})}\\
&\hspace{10cm}\times\|u_{L}
\|_{\widetilde{L}^{\frac{4}{3}}_{T}(B^{\N+\frac{1}{2}}_{2,2})}.
%&\hspace{3,6cm}\lesssim \|q^{n-1}
%\|_{\widetilde{L}^{\frac{4}{3}}_{T}(B^{\N+\frac{3}{2}}_{2,2})}\|q^{n-1}
%\|_{\widetilde{L}^{4}_{T}(B^{\N+\frac{1}{2}}_{2,2})}.
\end{aligned}
\label{tech1}
\end{equation}
Let us deal now with the term $G(\bar{q},\bar{u})$ and in particular the term $u\cdot\n u$, we have then:
$$u\cdot\n u=\bar{u}\cdot\n u+u_{L}\cdot\n\bar{u}+u_{L}\cdot\n u_{L}.$$
We have then by proposition \ref{produit}:
\begin{equation}
\begin{aligned}
&\|\bar{u}\cdot\n u\|_{\widetilde{L}^{1}_{T}(B^{\N-1}_{2,1})}\lesssim \|\bar{u}\|_{\widetilde{L}^{\infty}(B^{\N-1}_{2,1})}\|\n u\|_{\widetilde{L}^{1}(B^{\N}_{2,\infty})}+\|\n u\|_{\widetilde{L}^{2}(B^{\N-1}_{2,\infty})}\|\bar{u}\|_{\widetilde{L}^{2}(B^{\N}_{2,1})}
\\[2mm]
&\|u_{L}\cdot\n \bar{u}\|_{\widetilde{L}^{1}_{T}(B^{\N-1}_{2,1})}\lesssim \|u_{L}\|_{\widetilde{L}^{\infty}(B^{\N-1}_{2,\infty})}\|\n \bar{u}\|_{\widetilde{L}^{1}(B^{\N}_{2,1})}+\|\n \bar{u}\|_{\widetilde{L}^{2}(B^{\N-1}_{2,1})}\|u_{L}\|_{\widetilde{L}^{2}(B^{\N}_{2,\infty})}
\\[2mm]
&\|u_{L}\cdot\n u_{L}\|_{\widetilde{L}^{1}_{T}(B^{\N-1}_{2,1})}\lesssim \|u_{L}\|_{\widetilde{L}^{\infty}(B^{\N-1}_{2,2})}\|\n u_{L}\|_{\widetilde{L}^{1}(B^{\N}_{2,2})}+\|\n u_{L}\|_{\widetilde{L}^{2}(B^{\N-1}_{2,2})}\|u_{L}\|_{\widetilde{L}^{2}(B^{\N}_{2,2})}.
\\%\|_{\widetilde{L}^{\frac{4}{3}}_{T}(B^{\N+\frac{3}{2}}_{2,2})}\|q^{n-1}
%\|_{\widetilde{L}^{4}_{T}(B^{\N+\frac{1}{2}}_{2,2})}.
\end{aligned}
\label{tech3}
\end{equation}
We can proceed similarly for the terms $\n q\cdot Du$ and ${\rm div}u\,\n q$. Let us treat the last term $\n|\n q|^{2}$ we have then:
$$|\n q|^{2}=\n q\cdot\n \bar{q}+\n \bar{q}\cdot \n q_{L}+|\n q_{L}|^{2}.$$
By proposition \ref{produit}, we get:
\begin{equation}
\begin{aligned}
&\|\n (|\n q|^{2})\|_{\widetilde{L}^{1}_{T}(B^{\N-1}_{2,1})}\lesssim \||\n q|^{2}\|_{\widetilde{L}^{1}_{T}(B^{\N}_{2,1})},\\[2mm]
&\||\n q_{L}|^{2}\|_{\widetilde{L}^{1}_{T}(B^{\N-1}_{2,1})}\lesssim \|\n q_{L}\|_{\widetilde{L}^{\frac{4}{3}}_{T}(B^{\N+\frac{1}{2}}_{2,2})}\|\n q_{L}
\|_{\widetilde{L}^{4}_{T}(B^{\N-\frac{1}{2}}_{2,2})},\\
&\|\n q\cdot\n\bar{q}\|_{\widetilde{L}^{1}_{T}(B^{\N-1}_{2,1})}\lesssim \|\n q\|_{\widetilde{L}^{\frac{4}{3}}_{T}(B^{\N+\frac{1}{2}}_{2,\infty})}\|\n \bar{q}
\|_{\widetilde{L}^{4}_{T}(B^{\N-\frac{1}{2}}_{2,1})}+\|\n \bar{q}\|_{\widetilde{L}^{\frac{4}{3}}_{T}(B^{\N+\frac{1}{2}}_{2,1})}\|\n q
\|_{\widetilde{L}^{4}_{T}(B^{\N-\frac{1}{2}}_{2,\infty})},\\
\end{aligned}
\label{tech4}
\end{equation}
By collecting the estimates (\ref{tech1}), (\ref{tech3}) and (\ref{tech4}) and by interpolation we obtain that it exists $C>0$ such that:
\begin{equation}
\begin{aligned}
&\|\bar{q}\|_{\widetilde{L}^{\infty}(\R^{+},B^{\N-1}_{2,1}\cap B^{\N}_{2,1})}+\|\bar{q}\|_{\widetilde{L}^{1}(\R^{+},B^{\N+1}_{2,1}\cap B^{\N+2}_{2,1})}
+\|\bar{u}\|_{\widetilde{L}^{\infty}(\R^{+},B^{\N-1}_{2,1})}+\|\bar{u}\|_{\widetilde{L}^{1}(\R^{+},B^{\N+1}_{2,1})}\\
&\leq C\big( \|u_{L}\|_{\widetilde{L}^{\infty}(B^{\N-1}_{2,\infty})}+\|q_{L}\|_{\widetilde{L}^{\infty}(B^{\N-1,\N}_{2,\infty})}+ \|u_{L}\|_{\widetilde{L}^{1}(B^{\N+1}_{2,\infty})}+\|q_{L}\|_{\widetilde{L}^{1}(B^{\N+1,\N+2}_{2,\infty})}\\
&+\|u\|_{\widetilde{L}^{\infty}(B^{\N-1}_{2,\infty})}+\|q\|_{\widetilde{L}^{\infty}(B^{\N-1,\N}_{2,\infty})}+ \|u\|_{\widetilde{L}^{1}(B^{\N+1}_{2,\infty})}+\|q\|_{\widetilde{L}^{\infty}(B^{\N+1,\N+2}_{2,\infty})}\big)\\
&\times (\|\bar{q}\|_{\widetilde{L}^{\infty}(\R^{+},B^{\N-1,\N}_{2,1})}+\|\bar{q}\|_{\widetilde{L}^{1}(\R^{+},B^{\N+1,\N+2}_{2,1})}
+\|\bar{u}\|_{\widetilde{L}^{\infty}(\R^{+},B^{\N-1}_{2,1})}+\|\bar{u}\|_{\widetilde{L}^{1}(\R^{+},B^{\N+1}_{2,1})})\\
&+C\big( \|u_{L}\|_{\widetilde{L}^{\infty}(B^{\N-1}_{2,2})}+\|q_{L}\|_{\widetilde{L}^{\infty}(B^{\N-1,\N}_{2,2})}+ \|u_{L}\|_{\widetilde{L}^{1}(B^{\N+1}_{2,2})}+\|q_{L}\|_{\widetilde{L}^{\infty}(B^{\N+1,\N+2}_{2,2})})^{2}.
\end{aligned}
\label{crucialfin}
\end{equation}
Let us recall that via the first part of the theorem \ref{ftheo1} and the proposition \ref{1flinear3} (see the proof in the subsection \ref{infini}) we have for a $C>0$ and $M>0$ large enough:
$$
\begin{aligned}
&\|u_{L}\|_{\widetilde{L}^{\infty}(B^{\N-1}_{2,\infty})}+\|q_{L}\|_{\widetilde{L}^{\infty}(B^{\N-1,\N}_{2,\infty})}+ \|u_{L}\|_{\widetilde{L}^{1}(B^{\N+1}_{2,\infty})}+\|q_{L}\|_{\widetilde{L}^{\infty}(B^{\N+1,\N+2}_{2,\infty})}\leq C\e_{0},\\[2mm]
&\|u\|_{\widetilde{L}^{\infty}(B^{\N-1}_{2,\infty})}+\|q\|_{\widetilde{L}^{\infty}(B^{\N-1,\N}_{2,\infty})}+ \|u\|_{\widetilde{L}^{1}(B^{\N+1}_{2,\infty})}+\|q\|_{\widetilde{L}^{\infty}(B^{\N+1,\N+2}_{2,\infty})}\leq M\e_{0},
\end{aligned}
$$
with:
\begin{equation}
\e_{0}=\|q_{0}\|_{\widetilde{B}^{\N-1,\N}_{2,\infty}}+\|u_{0}\|_{B^{\N-1}_{2,\infty}}.
\label{definitep}
\end{equation}
By choosing $\e_{0}$ small enough we can apply a bootstrap argument in (\ref{crucialfin}) which implies that it exists $C>0$ such that
\begin{equation}
\begin{aligned}
&\|\bar{q}\|_{\widetilde{L}^{\infty}(\R^{+},B^{\N-1}_{2,1}\cap B^{\N}_{2,1})}+\|\bar{q}\|_{\widetilde{L}^{1}(\R^{+},B^{\N+1}_{2,1}\cap B^{\N+2}_{2,1})}
+\|\bar{u}\|_{\widetilde{L}^{\infty}(\R^{+},B^{\N-1}_{2,1})}+\|\bar{u}\|_{\widetilde{L}^{1}(\R^{+},B^{\N+1}_{2,1})}\\
&\leq C\big( \|u_{L}\|_{\widetilde{L}^{\infty}(B^{\N-1}_{2,2})}+\|q_{L}\|_{\widetilde{L}^{\infty}(B^{\N-1,\N}_{2,2})}+ \|u_{L}\|_{\widetilde{L}^{1}(B^{\N+1}_{2,2})}+\|q_{L}\|_{\widetilde{L}^{\infty}(B^{\N+1,\N+2}_{2,2})})^{2}.
\end{aligned}
\label{crucialfin1}
\end{equation}
It proves finally that $(\bar{q},\bar{u})$ is in $F$ and achieves this subsection.
\subsubsection*{$q$ is bounded in $L^{\infty}_{T}(L^{\infty}(\R^{N}))$ for any $T>0$}
In order to show that $q$ is bounded in $L^{\infty}_{T}(L^{\infty}(\R^{N}))$ for any $T>0$ we shall use the proposition \ref{maxfinal}. Indeed we have seen in the previous subsection that:
$$(q,u)=(q_{L},u_{L})+(\bar{q},\bar{u}),$$
with $(q_{L},u_{L})$ solution of the system $(N1)$ with $F=G=0$ $a=\mu$, $b=\lambda+\mu$, $c=\kappa$, $d=K$ and with initial data $(\ln\rho_{0},u_{0})$. In particular it implies via the definition of the semi group $B(t)$ in the proposition \ref{propimp} and the Duhamel formula that:
$$(q_{L}(t),u_{L}(t))=e^{B(t)}(q_{0},u_{0})+\int^{t}_{0}e^{B(t-s)}(0,K\n q_{L})(s)ds.$$
By using propositions \ref{flinear3}, \ref{maxfinal} and the embedding of $B^{\N}_{2,1}$ in $L^{\infty}$ we deduce that for any $T>0$ it exists $C>0$ independent on $T$ such that (here $[\cdot]_{1}$ defines the first coordinate of the vector field):
\begin{equation}
\begin{aligned}
\|q_{L}\|_{L^{\infty}_{T}(L^{\infty})}&\leq \|[e^{B(t)}(q_{0},u_{0})]_{1}\|_{L^{\infty}_{T}(L^{\infty})}+\|[\int^{t}_{0}e^{B(t-s)}(0,K\n q_{L})(s)ds]_{1}\|_{\widetilde{L}^{\infty}_{T}(B^{\N}_{2,1})},\\
&\leq C( \|\ln\rho_{0}\|_{L^{\infty}}+\|(\n q_{0},u_{0})\|_{B^{\N-1}_{2,2}}+\|\n q_{L}\|_{\widetilde{L}^{1}_{T}(B^{\N-1}_{2,1})}).
\end{aligned}
\label{crucomfin}
\end{equation}
From proposition \ref{1flinear3}, we know that it exists $C>0$ such that:
\begin{equation}
\|q_{L}\|_{\widetilde{L}^{\infty}(\widetilde{B}^{\N-1,\N}_{2,2})}+\|q_{L}\|_{\widetilde{L}^{1}(\widetilde{B}^{\N+1,\N+2}_{2,2})}\leq C(\|q_{0}\|_{\widetilde{B}^{\N-1,\N}_{2,2}}+\|u_{0}\|_{B^{\N-1}_{2,2}}).
\label{tecitan}
\end{equation}
Let us deal now with the term $\|\n q_{L}\|_{\widetilde{L}^{1}_{T}(B^{\N-1}_{2,1})}$ on the right hand side of (\ref{crucomfin}), we have by interpolation for a constant $C>0$:
$$\|q_{L}(t)\|_{B^{\N}_{2,1}}\leq C\|q_{L}(t)\|^{\frac{1}{2}}_{B^{\N-1}_{2,2}}\|q_{L}(t)\|^{\frac{1}{2}}_{B^{\N+1}_{2,2}}.$$
It implies that by H\"older's inequality and (\ref{Minko}) that it exists $C,C_{1}>0$ large enough such that:
\begin{equation}
\begin{aligned}
\|q_{L}(t)\|_{L^{1}_{T}(B^{\N}_{2,1})}&\leq C\|q_{L}\|^{\frac{1}{2}}_{L^{1}_{T}(B^{\N-1}_{2,2})}\|q_{L}\|^{\frac{1}{2}}_{L^{1}_{T}(B^{\N+1}_{2,2})},\\
&\leq CT^{\frac{3}{4}}\|q_{L}\|^{\frac{1}{2}}_{L^{\infty}_{T}(B^{\N-1}_{2,2})}\|q_{L}\|^{\frac{1}{2}}_{L^{2}_{T}(B^{\N+1}_{2,2})},\\
&\leq C_{1} T^{\frac{3}{4}}(\|q_{0}\|_{\widetilde{B}^{\N-1,\N}_{2,2}}+\|u_{0}\|_{B^{\N-1}_{2,2}}).
\end{aligned}
\label{pdernier}
\end{equation}
Combining (\ref{crucomfin}) and (\ref{pdernier}) we get:
\begin{equation}
\begin{aligned}
&\|q_{L}\|_{L^{\infty}_{T}(L^{\infty})}\leq \|[e^{B(t)}(q_{0},u_{0})]_{1}\|_{L^{\infty}_{T}(L^{\infty})}+\|[\int^{t}_{0}e^{B(t-s)}(0,K\n q_{L})(s)ds]_{1}\|_{\widetilde{L}^{\infty}_{T}(B^{\N}_{2,1})},\\
&\hspace{1cm}\leq C( \|\ln\rho_{0}\|_{L^{\infty}}+\|(\n q_{0},u_{0})\|_{B^{\N-1}_{2,2}}+ C_{1} T^{\frac{3}{4}}(\|q_{0}\|_{\widetilde{B}^{\N-1,\N}_{2,2}}+\|u_{0}\|_{B^{\N-1}_{2,2}}).
\end{aligned}
\label{crucomfin1}
\end{equation}
We have then proved that $q_{L}$ belongs in $L^{\infty}_{lot}(L^{\infty})$. Since we have seen that $(\bar{q},\bar{u})$ is in $F$, it implies in particular that $\bar{q}$ belongs in $\widetilde{L}^{\infty}(B^{\N}_{2,1})$ which is embedded in $L^{\infty}(L^{\infty})$. We deduce that $q=q_{L}+\bar{q}$ belongs in $L^{\infty}_{loc}(L^{\infty})$.\\
 It remains only to prove that $(\rho,u)=(\exp(q),u)$ is a global strong solution of (\ref{3systeme}) and we define $h=\exp(q)-1$ with $\rho=1+h$. By proposition \ref{composition} and the fact that $\rho$ and $\frac{1}{\rho}$  belong in $L^{\infty}_{loc}(L^{\infty})$ we easily show that $h$ is for any $T>0$ in:
 $$H_{T}=\big(\widetilde{L}^{\infty}_{T}(\widetilde{B}^{\N-1,\N}_{2,\infty})\cap \widetilde{L}^{1}_{T}(\widetilde{B}^{\N+1,\N+2}_{2,\infty}\big).$$
With such regularity on $h$ and $u$, the verification that $(\rho,u)$ is a solution of (\ref{3systeme}) in the sense of distribution is a straightforward application of propositions \ref{produit} and \ref{composition}. The uniqueness follows also the same line than in the proof of subsection \ref{sub32}.
\hfill {$\Box$}
\subsection{Proof of the theorem \ref{ftheo2}}
In the theorem \ref{ftheo2} we are interested in extending the results of theorem \ref{ftheo1} to the case of general pressure $P$ and also to the case of constant viscosity and capillary coefficients. It implies in particular that the eulerian form of the system \ref{3systeme} does not only depend on the unknown $q=\ln\rho$. Typically the pressure term $\frac{\n P(\rho)}{\rho}$ is non linear in terms of $q=\ln\rho$, that is why we need to control the  $L^{\infty}$ norm of the density in order to estimate this term in Besov space, indeed otherwise  we could have a loss of regularity on this term. For this reason it seems impossible to hope a result of global strong solution involving only a smallness hypothesis on $\|q_{0}\|_{B^{\N}_{2,\infty}}$ and $\|u_{0}\|_{B^{\N-1}_{2,\infty}}$ as it is the case under some specific physical condition in the theorem \ref{ftheo1}.\\
We are going now to explain how to adapt the previous arguments of the proof of theorem \ref{ftheo1}  to this new situation. We are only dealing with the case $\mu(\rho)=\mu$, $\lambda(\rho)=\lambda$, $\kappa(\rho)=\kappa$ and $P$ a regular function such that $P^{'}(1)>0$, here $\mu>0$ $2\mu+\lambda>0$ and $ \kappa>0$ (the other case are similar to treat). The system \ref{3systeme} is equivalent to the following system where we set in this section $\rho=1+h$. Let us mention that in this case a straighforward  calculus gives:
$${\rm div}K=\kappa\rho\n\D\rho.$$
We can then rewrite the system (\ref{3systeme}) as follows:
\begin{equation}
\begin{aligned}
\begin{cases}
&\p_{t}h+{\rm div}u=F(h,u),\\
&\p_{t}u-\mu\D u-(\mu+\lambda)
\n{\rm div}u+F^{'}(1)\n h-\kappa\n\D h=G(h,u),\\
&(h,u)_{\ t=0}=(h_{0},u_{0}),
\end{cases}
\end{aligned}
\label{0.6}
\end{equation}
with $F^{'}(\rho)=\frac{P^{'}(\rho)}{\rho}$ and:
\begin{equation}
\begin{aligned}
&F(h,u)=-u\cdot\n h-h{\rm div}u,\\
&G(h,u)=(F^{'}(1)-F^{'}(\rho))\n h+(\frac{\mu}{\rho}-\mu)\D u+(\frac{\mu+\lambda}{\rho}-\mu-\lambda)\n{\rm div} u-u\cdot\n u.
\end{aligned}
\label{reste}
\end{equation}
Let $(h_{L},u_{L})$ the solution of $(N1)$ with $a=\mu$, $b=\lambda+\mu$, $c=\kappa$, $d=P^{'}(1)$, $F=G=0$ and the initial data $(h_{0},u_{0})$. In the sequel we will consider solution under the form:
$$(h,u)=(h_{L},u_{L})+(\bar{h},\bar{u}).$$
In order to prove theorem \ref{ftheo2} we shall use a fixed point theorem, more precisely we are going to consider the following functional:
\begin{equation}
\psi(\bar{h},\bar{u})=\int_{0}^{t}W(t-s)\left(\begin{array}{c}
F(h,u)\\
G(h,u)\\
\end{array}
\right)\ ds\; .\label{f1a1}
\end{equation}
$W(t)$ is the semi group associated to $(N1)$ with the following conditions $a=\mu$, $b=\lambda+\mu$, $c=\kappa$, $d=P^{'}(1)$. The proof is divided in two step the stability of $\psi$ for a ball $B(0,R)$ in $E^{\N}$ defined below and the contraction property. We consider $E^{\N}$ defined as follows:
$$E^{\N}=\big(\widetilde{L}^{\infty}(\widetilde{B}^{\N-1,\N}_{2,1})\cap \widetilde{L}^{1}(\widetilde{B}^{\N+1,\N+2}_{2,1})\big)\times\big(\widetilde{L}^{\infty}(B^{\N-1}_{2,1})\cap \widetilde{L}^{1}(B^{\N+1}_{2,1})\big)^{N}.$$
\subsubsection*{1) First step, stability of $B(0,R)$:}
By using the proposition \ref{1flinear3} we have: for $C>0$:
%$$E_{T}=\widetilde{C}_{T}(B^{\N}_{2,1})\times\big(\widetilde{C}_{T}(B^{\fd-1}_{2,1})\cap\widetilde{L}^{1}_{T}(
%B^\fd_{2,1})\big).$$
%{\bf GENERALISER AU SECOND TERME AVEC PRESSION} Before doing this, let us recall that thanks to proposition \ref{flinear3} and \ref{maxfinal}  we have for all $T>0$:
\begin{equation}
\begin{aligned}
&\|h_{L}\|_{\widetilde{L}_{T}^{\infty}(\widetilde{B}^{\N-1,\N}_{2,2})}+\|h_{L}\|_{\widetilde{L}^{1}_{T}(\widetilde{B}^{\N+1,\N+2}_{2,2})}+\|u_{L}\|_{\widetilde{L}_{T}^{\infty}(B^{\N-1}_{2,2})}+\|u_{L}\|_{\widetilde{L}^{1}_{T}(B^{\N+1}_{2,2})}\\
&\hspace{7cm}\leq C(\|h_{0}\|_{\widetilde{B}^{\N-1,\N}_{2,2}}+\|u_{0}\|_{B^{\N-1}_{2,2}}).
\end{aligned}
\label{lineaire1d}
\end{equation}
By combining the proposition \ref{maxfinal} and \ref{1flinear3}, we have for any $T>0$ and $C>0$ independent on $T$:
\begin{equation}
\|h_{l}\|_{L^{\infty}_{T}(L^{\infty})}\leq C(\|h_{0}\|_{L^{\infty}}++\|h_{0}\|_{B^{\N-2}_{2,1}}+\|u_{0}\|_{\widetilde{B}^{\N-2}_{2,1}}).
\label{lineaired2}
\end{equation}
\begin{remarka}
\label{remtecg}
Let us mention that the condition $(h_{0},u_{0})\in B^{\N-2}_{2,1}\times (B^{\N-2}_{2,1})^{N}$ plays only a role in the previous estimate (\ref{lineaired2}) in order to bound the term $P^{'}(1)\n h_{L}$ that we consider as a remainder term for the system $(N)$. Indeed we are interested in applying proposition \ref{maxfinal} and the Duhamel formula. It would be possible to avoid this additional regularity by extending the proposition \ref{maxfinal} to the system $(N1)$.
\end{remarka}
It remains only to apply the proposition \ref{1flinear3} in order to get a priori estimates on $(\bar{h},\bar{u})$, indeed we have for $C>0$:
\begin{equation}
\begin{aligned}
&\|\psi(\bar{h},\bar{u})\|_{E^{\N}}\leq C(\|F(h,u)\|_{\widetilde{L}_{T}^{1}(\widetilde{B}^{\N-1,\N}_{2,1})}+\|G(h,u)\|_{\widetilde{L}_{T}^{1}(\widetilde{B}^{\N-1,\N}_{2,1})}).
\end{aligned}
\label{lineaire1a}
\end{equation}
By applying proposition \ref{produit} and lemma \ref{composition} (indeed we have seen via the estimate (\ref{lineaire1d}) and the fact that $\bar{h}$ is in $\widetilde{L}^{\infty}(B^{\N}_{2,1})$ that $h$ belongs in $L^{\infty}$ and $\rho=1+h\geq c>0$) we have as in (\ref{crucialfin}) for a function $C>0$:
\begin{equation}
\begin{aligned}
&\|\psi(\bar{h},\bar{u})\|_{E^{\N}}\leq C(\|h\|_{L^{\infty}})\big( \|u_{L}\|_{\widetilde{L}^{\infty}(B^{\N-1}_{2,\infty})}+\|h_{L}\|_{\widetilde{L}^{\infty}(B^{\N-1,\N}_{2,\infty})}+ \|h_{L}\|_{\widetilde{L}^{1}(B^{\N+1}_{2,\infty})}\\
&+\|h_{L}\|_{\widetilde{L}^{1}(B^{\N+1,\N+2}_{2,\infty})}+\|u\|_{\widetilde{L}^{\infty}(B^{\N-1}_{2,\infty})}+\|h\|_{\widetilde{L}^{\infty}(B^{\N-1,\N}_{2,\infty})}+ \|u\|_{\widetilde{L}^{1}(B^{\N+1}_{2,\infty})}+\|h\|_{\widetilde{L}^{\infty}(B^{\N+1,\N+2}_{2,\infty})}\big)\\
&\times (\|\bar{h}\|_{\widetilde{L}^{\infty}(\R^{+},B^{\N-1,\N}_{2,1})}+\|\bar{h}\|_{\widetilde{L}^{1}(\R^{+},B^{\N+1,\N+2}_{2,1})}
+\|\bar{u}\|_{\widetilde{L}^{\infty}(\R^{+},B^{\N-1}_{2,1})}+\|\bar{u}\|_{\widetilde{L}^{1}(\R^{+},B^{\N+1}_{2,1})})\\
&+C(\|h\|_{L^{\infty}})\big( \|u_{L}\|_{\widetilde{L}^{\infty}(B^{\N-1}_{2,2})}+\|h_{L}\|_{\widetilde{L}^{\infty}(B^{\N-1,\N}_{2,2})}+ \|u_{L}\|_{\widetilde{L}^{1}(B^{\N+1}_{2,2})}+\|h_{L}\|_{\widetilde{L}^{\infty}(B^{\N+1,\N+2}_{2,2})})^{2}.
\end{aligned}
\label{bcrucialfin}
\end{equation}
We set
\begin{equation}
\e_{0}=\|h_{0}\|_{\widetilde{B}^{\N-1,\N}_{2,2}}+\|u_{0}\|_{B^{\N-1}_{2,\infty}}+\|h_{0}\|_{L^{\infty}}+\|h_{0}\|_{B^{\N-2}_{2,1}}+\|u_{0}\|_{B^{\N-2}_{2,1}}.
\label{definitep}
\end{equation}
By choosing $R=M\e_{0}$ with $M>2$, it suffices to choose $\e_{0}$ small enough such that (\ref{bcrucialfin}) ensures that:
$$\|\psi(\bar{h},\bar{u})\|_{E^{\N}}\leq R.$$
The proof of the contraction follows the same lines as in the proof of the theorem \ref{ftheo1} and is left to the reader. It achieves the proof of the theorem \ref{ftheo2}.
\hfill {$\Box$}
\section{Proof of the theorem \ref{ftheo3}}
\label{section4}
We are now interested in proving the theorem \ref{ftheo3} where we assume that $\kappa(\rho)=\frac{\mu^{2}}{\rho}$, $\mu(\rho)=\mu\rho$ and $\lambda(\rho)=0$. We have observed that it exists quasi solution, it means a solution of the form $(\rho_{1},-\mu\n\ln\rho_{1})$ with:
$$
\begin{cases}
\begin{aligned}
&\p_{t}\rho_{1}-\mu\D\rho_{1}=0,\\
&(\rho_{1})_{t=0}=(\rho_{1})_{0}.
\end{aligned}
\end{cases}
$$
This quasi solution $(\rho_{1},-\mu\n\ln\rho_{1})$ verifies the following system:
\begin{equation}
\begin{cases}
\begin{aligned}
&\frac{\p}{\p t}\rho+{\rm div}(\rho u)=0,\\
&\frac{\p}{\p t}(\rho u)+{\rm div}(\rho
u\otimes u)-\rm div(2\mu\rho D (u))-\mu^{2}\rho\n\D\ln\rho-\frac{\mu^{2}}{2}\n(|\n\ln\rho|^{2})=0,\\
&(\rho,u)_{t=0}=\big((\rho_{1})_{0},-\mu\n\ln (\rho_{1})_{0}\big) .
\end{aligned}
\end{cases}
\label{3systemequasi}
\end{equation}
We are now interested in working around this quasi solution, more precisely we search global solution of (\ref{3systeme}) under the form:
$$q=\ln\rho=\ln\rho_{1}+h_{2}\;\;\mbox{with}\;\rho_{1}=1+h_{1}\;\;\mbox{and}\;\;u=-\mu\n\ln\rho_{1}+u_{2}.$$ 
We deduce from (\ref{3systeme2ab}) that $(h_{2},u_{2})$ verifies the following system:
\begin{equation}
\begin{cases}
\begin{aligned}
&\p_{t}h_{2}+{\rm div}u_{2}-\mu\n\ln\rho_{1}\cdot\n h_{2}+u_{2}\cdot\n\ln\rho_{1}=F(h_{2},u_{2}),\\[2mm]
&\p_{t}u_{2}-\mu\D u_{2}-\mu\n{\rm div}u_{2}-\kappa\n\D h_{2}+K\n h_{2}-2\mu\n \ln\rho_{1}\cdot D u_{2}-2\mu\n h_{2}\cdot D u_{1}\\
&\hspace{4cm}+u_{1}\cdot\n u_{2}+u_{2}\cdot\n u_{1}-\mu^{2}\n(\n\ln\rho_{1}\cdot\n h_{2})=G(h_{2},u_{2}),\\
&(h_{2}(0,\cdot),u_{2}(0,\cdot))=(h_{0}^{2},u_{0}^{2}).
\end{aligned}
\end{cases}
\label{afond3systeme2}
\end{equation}
with:
\begin{equation}
\begin{aligned}
&F(h_{2},u_{2})=-u_{2}\cdot\n h_{2},\\
&G(h_{2},u_{2})=-u_{2}\cdot\n u_{2}+2\mu\n h_{2}\cdot D u_{2}-K\n\ln\rho_{1}+\frac{\mu^{2}}{2}\n(|\n h_{2}|^{2}).
\end{aligned}
\label{restetech}
\end{equation}
%Here we define $v$, $v_{1}$ and $v_{2}$ by:
%$$
%\begin{aligned}
%&v=v_{1}+v_{2}=u+\mu\n\ln\rho,\;v_{1}=0\;\;\mbox{and}\;\;v_{2}=u_{2}+\mu\n h_{2}.
%\end{aligned}
%$$
We have now to solve the previous system (\ref{afond3systeme2}) and to do this we are going to apply a fixed point theorem. We start with defining the following map $\psi$:
\begin{equation}
\psi(h,u)=W_{1}(t,\cdot)*\left(\begin{array}{c}
h_{0}^{2}\\
u_{0}^{2}\\
\end{array}
\right)+\int_{0}^{t}W_{1}(t-s)\left(\begin{array}{c}
F(h,u)\\
G(h,u)\\
\end{array}
\right)\ ds\; .\label{f1a1}
\end{equation}
where $W_{1}$ is the semi group associated to the following linear system (\ref{afond3systeme3}):
\begin{equation}
\begin{cases}
\begin{aligned}
&\p_{t}h_{2}+{\rm div}u_{2}-\mu\n\ln\rho_{1}\cdot\n h_{2}+u_{2}\cdot\n\ln\rho_{1}=0,\\[2mm]
&\p_{t}u_{2}-\mu\D u_{2}-\mu\n{\rm div}u_{2}-\kappa\n\D h_{2}+K\n h_{2}-2\mu\n \ln\rho_{1}\cdot D u_{2}-2\mu\n h_{2}\cdot D u_{1}\\
&\hspace{5cm}+u_{1}\cdot\n u_{2}+u_{2}\cdot\n u_{1}-\mu^{2}\n(\n\ln\rho_{1}\cdot\n h_{2})=0.\\
&(h_{2}(0,\cdot),u_{2}(0,\cdot))=(h_{0}^{2},u_{0}^{2}).
\end{aligned}
\end{cases}
\label{afond3systeme3}
\end{equation}
%In what follows we set:
%$$\rho=\bar{\rho}(1+q)\;,\;\theta=\bar{\theta}+{\cal T}\;,\;\widetilde{T}=\Psi^{-1}(\theta).$$
The non linear terms  $F,G$ are defined in (\ref{restetech}). We are going to check that we can apply a fixed point theorem for the function $\psi$, the proof is divided in two step the stability of $\psi$ for a ball $B(0,R)$ in $E^{\N}$ and the contraction property. We define now $E^{\N}$ by:
$$E^{\N}=\big(\widetilde{L}^{\infty}(\widetilde{B}^{\N-1,\N}_{2,1})\cap \widetilde{L}^{1}(\widetilde{B}^{\N+1,\N+2}_{2,1})\big)\times\big(\widetilde{L}^{\infty}(B^{\N-1}_{2,1})\cap \widetilde{L}^{1}(B^{\N+1}_{2,1})\big)^{N}.$$
\subsubsection*{1) First step, stability of $B(0,R)$:}
Let:
$$\eta=\| h^{2}_{0}\|_{\widetilde{B}^{\N-1,\N}_{2,\infty}}+\|u^{2}_{0}\|_{B^{\N-1}_{2,\infty}}\;.$$
We are going to show that $\psi$  maps the ball $B(0,R)$ into itself if $R$ is small enough. By using the proposition \ref{fcrucialprop} we have:
%According to proposition \ref{1flinear4}, we have:
%\begin{equation}
%\|W(t,\cdot)*\left(\begin{array}{c}
%q_{0}\\
%u_{0}\\
%\end{array}
%\right)\|_{E^{\N}}\leq C\eta \;.
%\label{f1a5}
%\end{equation}
%According to the proposition \ref{1flinear3} it implies also that:
\begin{equation}
\begin{aligned}
&\|\psi(h,u)\|_{E^{\N}}\leq%
%\;\begin{equation}
%\begin{aligned}
%&\|h_{2}\|_{\widetilde{L}^{1}_{T}(\widetilde{B}^{\N+1,\N+2}_{2,1})}+\|h_{2}\|_{\widetilde{L}^{\infty}_{T}(B^{\N-1,\N}_{2,1})}
%+\|u_{2}\|_{\widetilde{L}^{1}_{T}(B^{\N+1}_{2,1})}+ \|u_{2}\|_{\widetilde{L}^{\infty}_{T}(B^{\N-1}_{2,1})}\\
%&\leq C\big(\|h^{2}_{0}\|_{\widetilde{B}^{\N-1,\N}_{2,1}}+\|u^{2}_{0}\|_{B^{\N}_{2,1}}+\|F\|_{\widetilde{L}^{1}_{T}(\widetilde{B}^{\N-1,\N}_{2,1})}
%+\|G\|_{\widetilde{L}^{1}_{T}(B^{\N-1}_{2,1})}\big)\\
%&\times \exp\biggl(C\int_{\R^{+}}\big( \|u_{1}\|_{B^{\N-\frac{1}{2}}_{2,\infty}}^{4}+\|u_{1}\|^{\frac{4}{3}}_{B^{\N+\frac{1}{2}}_{2,\infty}}+\|\n\ln\rho_{1}\|^{\frac{4}{3}}_{B^{\N+\frac{1}{2}}_{2,\infty}}+\|\n\ln\rho_{1}\|^{\frac{4}{3}}_{B^{\N+\frac{1}{2}}_{2,\infty}}\big)(s)ds\biggl).
%\end{aligned}
%\label{raimportant}
%\end{equation}
C\big(\eta+\|F(h,u)\|_{\widetilde{L}^{1}_{T}(\widetilde{B}^{\N-1,\N}_{2,1})}
+\|G(h,u)\|_{\widetilde{L}^{1}_{T}(B^{\N-1}_{2,1})}\big)\\
&\times \exp\biggl(C\int_{\R^{+}}\big( \|u_{1}\|_{B^{\N-\frac{1}{2}}_{2,\infty}}^{4}+\|u_{1}\|^{\frac{4}{3}}_{B^{\N+\frac{1}{2}}_{2,\infty}}+\|\n\ln\rho_{1}\|^{\frac{4}{3}}_{B^{\N+\frac{1}{2}}_{2,\infty}}+\|\n\ln\rho_{1}\|^{\frac{4}{3}}_{B^{\N+\frac{1}{2}}_{2,\infty}}\big)(s)ds\biggl).
%C\big(\eta+
%\|F(q,u)\|_{ \widetilde{L}^{1}(\widetilde{B}^{\N-1,\N}_{2,\infty})}+\|G(q,u)\|_{ \widetilde{L}^{1}(B^{\N-1}_{2,\infty})}\big).
\end{aligned}
\label{bf1a6}
\end{equation}
The main task consists in using the propositions \ref{produit} and corollary \ref{produit2} to obtain estimates on
$$\|F(h,u)\|_{ \widetilde{L}^{1}(\widetilde{B}^{\N-1,\N}_{2,\infty})},\;\|G(h,u)\|_{ \widetilde{L}^{1}(B^{\N-1}_{2,\infty})}.$$
Let us first estimate $\|F(h,u)\|_{ \widetilde{L}^{1}(\widetilde{B}^{\N-1,\N}_{2,\infty})}$. According to proposition
\ref{produit},
we have:
\begin{equation}
\begin{aligned}
&\|u\cdot\n h\|_{
\widetilde{L}^{1}(\widetilde{B}^{\N-1,\N}_{2,1})}\lesssim \|\n h\|_{\widetilde{L}^{\frac{4}{3}}(\widetilde{B}^{\N-\frac{1}{2},\N+\frac{1}{2}}_{2,1})} \|u\|_{\widetilde{L}^{4}(B^{\N-\frac{1}{2}}_{2,1})}\\
&\hspace{5cm}
+\|\n h\|_{\widetilde{L}^{4}(\widetilde{B}^{\N-\frac{3}{2},\N-\frac{1}{2}}_{2,1})} \|u\|_{\widetilde{L}^{\frac{4}{3}}(B^{\N+\frac{1}{2}}_{2,1})}.
\end{aligned}
\label{bestim1}
\end{equation}
Similarly for $\|G(h,u)\|_{\widetilde{L}^{1}_{T}(B^{\N-1}_{2,1})}$ we have:
\begin{equation}
\begin{aligned}
&\|u\cdot\n u\|_{
\widetilde{L}^{1}(B^{\N-1}_{2,1})}\lesssim \|u\|_{\widetilde{L}^{\frac{4}{3}}(B^{\N+\frac{1}{2}}_{2,1})} \|\n u\|_{\widetilde{L}^{4}(B^{\N-\frac{3}{2}}_{2,1})}\\
&\hspace{5cm}
+\| u\|_{\widetilde{L}^{4}(B^{\N-\frac{1}{2}}_{2,1})} \|\n u\|_{\widetilde{L}^{\frac{4}{3}}(B^{\N-\frac{1}{2}}_{2,1})}.
\end{aligned}
\label{bestim2}
\end{equation}
The most important term is certainly $K\n\ln\rho_{1}$ which belongs in $\widetilde{L}^{1}(B^{\N-1}_{2,1})$ since we have assume that $\ln\rho_{1}^{0}$ belongs in $B^{\N-2}_{2,1}\cap L^{\infty}$ and the fact that $\rho_{1}$ verifies a heat equation (\ref{eqchaleur}). Using propositions \ref{composition}, \ref{chaleur} and the maximum principle it exists a regular function $g$ such that:
\begin{equation}
\begin{aligned}
\|\n\ln\rho_{1}\|_{\widetilde{L}^{1}(B^{\N-1}_{2,1})}&\lesssim \|\ln \rho_{1}\|_{
\widetilde{L}^{1}(B^{\N}_{2,1})},\\
&\lesssim g(\|(\rho_{1},\frac{1}{\rho_{1}})\|_{L^{\infty}})\|h_{1}\|_{
\widetilde{L}^{1}(B^{\N}_{2,1})},\\
&\lesssim g(\|(\rho_{1}^{0},\frac{1}{\rho_{1}^{0}})\|_{L^{\infty}})\|h^{0}_{1}\|_{B^{\N-2}_{2,1}}.
\end{aligned}
\label{bestim3}
\end{equation}
We deal with the others terms in a similar way and combining (\ref{bf1a6}), (\ref{bestim1}), (\ref{bestim2}) and (\ref{bestim3}) we obtain for a $C>0$ large enough and using the fact that  $(h,u)$ is in the ball $B(0,R)$ of $E^{\N}$:
\begin{equation}
\begin{aligned}
&\|\psi(h,u)\|_{E^{\N}}\leq%
%\;\begin{equation}
%\begin{aligned}
%&\|h_{2}\|_{\widetilde{L}^{1}_{T}(\widetilde{B}^{\N+1,\N+2}_{2,1})}+\|h_{2}\|_{\widetilde{L}^{\infty}_{T}(B^{\N-1,\N}_{2,1})}
%+\|u_{2}\|_{\widetilde{L}^{1}_{T}(B^{\N+1}_{2,1})}+ \|u_{2}\|_{\widetilde{L}^{\infty}_{T}(B^{\N-1}_{2,1})}\\
%&\leq C\big(\|h^{2}_{0}\|_{\widetilde{B}^{\N-1,\N}_{2,1}}+\|u^{2}_{0}\|_{B^{\N}_{2,1}}+\|F\|_{\widetilde{L}^{1}_{T}(\widetilde{B}^{\N-1,\N}_{2,1})}
%+\|G\|_{\widetilde{L}^{1}_{T}(B^{\N-1}_{2,1})}\big)\\
%&\times \exp\biggl(C\int_{\R^{+}}\big( \|u_{1}\|_{B^{\N-\frac{1}{2}}_{2,\infty}}^{4}+\|u_{1}\|^{\frac{4}{3}}_{B^{\N+\frac{1}{2}}_{2,\infty}}+\|\n\ln\rho_{1}\|^{\frac{4}{3}}_{B^{\N+\frac{1}{2}}_{2,\infty}}+\|\n\ln\rho_{1}\|^{\frac{4}{3}}_{B^{\N+\frac{1}{2}}_{2,\infty}}\big)(s)ds\biggl).
%\end{aligned}
%\label{raimportant}
%\end{equation}
C\big(\eta+R^{2}+g(\|(\rho_{1}^{0},\frac{1}{\rho_{1}^{0}})\|_{L^{\infty}})\|h^{0}_{1}\|_{B^{\N-2}_{2,1}}\big)\\
&\times \exp\biggl(C\int_{\R^{+}}\big( \|u_{1}\|_{B^{\N-\frac{1}{2}}_{2,\infty}}^{4}+\|u_{1}\|^{\frac{4}{3}}_{B^{\N+\frac{1}{2}}_{2,\infty}}+\|\n\ln\rho_{1}\|^{\frac{4}{3}}_{B^{\N+\frac{1}{2}}_{2,\infty}}+\|\n\ln\rho_{1}\|^{\frac{4}{3}}_{B^{\N+\frac{1}{2}}_{2,\infty}}\big)(s)ds\biggl).
\end{aligned}
\label{supcrucial}
\end{equation}
Let us recall that $\rho_{1}$ verifies a heat equation (\ref{eqchaleur}), we deduce from propositions \ref{chaleur} and the maximum principle that it exists a regular function $g_{1}$ such that:
\begin{equation}
\begin{aligned}
&\int_{\R^{+}}\big( \|u_{1}\|_{B^{\N-\frac{1}{2}}_{2,\infty}}^{4}+\|u_{1}\|^{\frac{4}{3}}_{B^{\N+\frac{1}{2}}_{2,\infty}}+\|\n\ln\rho_{1}\|^{\frac{4}{3}}_{B^{\N+\frac{1}{2}}_{2,\infty}}+\|\n\ln\rho_{1}\|^{\frac{4}{3}}_{B^{\N+\frac{1}{2}}_{2,\infty}}\big)(s)ds\\
&\hspace{8cm}\lesssim g_{1}(\|(\rho^{0}_{1},\frac{1}{\rho^{0}_{1}})\|_{L^{\infty}})\|h^{0}_{1}\|_{
B^{\N}_{2,1}}.
\end{aligned}
\label{supimp}
\end{equation}
From (\ref{supimp}) we deduce that for a $C$ large enough we have:
\begin{equation}
\begin{aligned}
&g(\|(\rho^{0}_{1},\frac{1}{\rho^{0}_{1}})\|_{L^{\infty}})\|h^{0}_{1}\|_{
B^{\N-2}_{2,1}}\times\exp\biggl(C\int_{\R^{+}}\big( \|u_{1}\|_{B^{\N-\frac{1}{2}}_{2,\infty}}^{4}+\|u_{1}\|^{\frac{4}{3}}_{B^{\N+\frac{1}{2}}_{2,\infty}}\\
&\hspace{4cm}+\|\n\ln\rho_{1}\|^{\frac{4}{3}}_{B^{\N+\frac{1}{2}}_{2,\infty}}+\|\n\ln\rho_{1}\|^{\frac{4}{3}}_{B^{\N+\frac{1}{2}}_{2,\infty}}\big)(s)ds\biggl)\\
&\leq g(\|(\rho^{0}_{1},\frac{1}{\rho^{0}_{1}})\|_{L^{\infty}})\|h^{0}_{1}\|_{
B^{\N-2}_{2,1}}\exp\big(C g_{1}( \|(\rho^{0}_{1},\frac{1}{\rho^{0}_{1}})\|_{L^{\infty}})\|h^{0}_{1}\|_{
B^{\N}_{2,1}}\big).
\end{aligned}
\label{tropimportant}
\end{equation}
In particular it implies that:
\begin{equation}
\begin{aligned}
&C\eta \exp\biggl(C\int_{\R^{+}}\big( \|u_{1}\|_{B^{\N-\frac{1}{2}}_{2,\infty}}^{4}+\|u_{1}\|^{\frac{4}{3}}_{B^{\N+\frac{1}{2}}_{2,\infty}}+\|\n\ln\rho_{1}\|^{\frac{4}{3}}_{B^{\N+\frac{1}{2}}_{2,\infty}}+\|\n\ln\rho_{1}\|^{\frac{4}{3}}_{B^{\N+\frac{1}{2}}_{2,\infty}}\big)(s)ds\biggl)\\
&\hspace{7cm}\leq C\eta \exp (Cg_{1}(\|(\rho^{0}_{1},\frac{1}{\rho^{0}_{1}})\|_{L^{\infty}})\|h^{0}_{1}\|_{
B^{\N}_{2,1}}).
\end{aligned}
\label{imp11}
\end{equation}
Let us prove now the stability of the functional $\psi$, from label (\ref{supcrucial}) let us choose:
 \begin{equation}
 R=4 C\eta \exp (Cg_{1}(\|(\rho^{0}_{1},\frac{1}{\rho^{0}_{1}})\|_{L^{\infty}})\|h^{0}_{1}\|_{
B^{\N}_{2,1}}).
 \label{valeurR}
 \end{equation}
 and suppose that:
\begin{equation}
\begin{aligned}
& g(\|(\rho^{0}_{1},\frac{1}{\rho^{0}_{1}})\|_{L^{\infty}})\|h^{0}_{1}\|_{
B^{\N-2}_{2,1}}\leq C\eta.
\end{aligned}
\label{tropimportant1}
\end{equation}
%Let us point out that this condition corresponds exactly to the condition (\ref{crucinitial}) of the theorem \ref{ftheo3}. 
Now combining (\ref{supcrucial}), (\ref{valeurR}) and (\ref{tropimportant1}) it yields:
\begin{equation}
\begin{aligned}
&\|\psi(h,u)\|_{E^{\N}}\leq%
%\;\begin{equation}
%\begin{aligned}
%&\|h_{2}\|_{\widetilde{L}^{1}_{T}(\widetilde{B}^{\N+1,\N+2}_{2,1})}+\|h_{2}\|_{\widetilde{L}^{\infty}_{T}(B^{\N-1,\N}_{2,1})}
%+\|u_{2}\|_{\widetilde{L}^{1}_{T}(B^{\N+1}_{2,1})}+ \|u_{2}\|_{\widetilde{L}^{\infty}_{T}(B^{\N-1}_{2,1})}\\
%&\leq C\big(\|h^{2}_{0}\|_{\widetilde{B}^{\N-1,\N}_{2,1}}+\|u^{2}_{0}\|_{B^{\N}_{2,1}}+\|F\|_{\widetilde{L}^{1}_{T}(\widetilde{B}^{\N-1,\N}_{2,1})}
%+\|G\|_{\widetilde{L}^{1}_{T}(B^{\N-1}_{2,1})}\big)\\
%&\times \exp\biggl(C\int_{\R^{+}}\big( \|u_{1}\|_{B^{\N-\frac{1}{2}}_{2,\infty}}^{4}+\|u_{1}\|^{\frac{4}{3}}_{B^{\N+\frac{1}{2}}_{2,\infty}}+\|\n\ln\rho_{1}\|^{\frac{4}{3}}_{B^{\N+\frac{1}{2}}_{2,\infty}}+\|\n\ln\rho_{1}\|^{\frac{4}{3}}_{B^{\N+\frac{1}{2}}_{2,\infty}}\big)(s)ds\biggl).
%\end{aligned}
%\label{raimportant}
%\end{equation}
2 C\eta \exp (Cg_{1}(\|(\rho^{0}_{1},\frac{1}{\rho^{0}_{1}})\|_{L^{\infty}})\|h^{0}_{1}\|_{
B^{\N}_{2,1}})+C R^{2} \exp (Cg_{1}(\|(\rho^{0}_{1},\frac{1}{\rho^{0}_{1}})\|_{L^{\infty}})\|h^{0}_{1}\|_{
B^{\N}_{2,1}}).
\end{aligned}
\label{supcrucial1}
\end{equation}
From (\ref{supcrucial1}) we need to assume that:
\begin{equation}
R^{2}\leq 2 C\eta \exp (Cg_{1}(\|(\rho^{0}_{1},\frac{1}{\rho^{0}_{1}})\|_{L^{\infty}})\|h^{0}_{1}\|_{
B^{\N}_{2,1}}).
\label{tropimportant2}
\end{equation}
Indeed combining (\ref{tropimportant2}), (\ref{supcrucial1}) and (\ref{valeurR}) shows that:
\begin{equation}
\|\psi(h,u)\|_{E^{\N}}\leq R.
\label{stabilitefinal}
\end{equation}
This concludes the proof of the stability except that it remains to choose $\eta$ and to verify that it implies the condition \ref{crucinitial} of theorem \ref{ftheo3}.
The condition \ref{tropimportant2} via (\ref{valeurR}) is equivalent to:
\begin{equation}
\begin{aligned}
16 C^{2}\eta^{2} \exp (2Cg_{1}(\|(\rho^{0}_{1},\frac{1}{\rho^{0}_{1}})\|_{L^{\infty}})\|h^{0}_{1}\|_{
B^{\N}_{2,1}})\leq 2 C\eta \exp (Cg_{1}(\|(\rho^{0}_{1},\frac{1}{\rho^{0}_{1}})\|_{L^{\infty}})\|h^{0}_{1}\|_{
B^{\N}_{2,1}}).
\end{aligned}
\end{equation}
and:
\begin{equation}
\begin{aligned}
8 C\eta \exp (Cg_{1}(\|(\rho^{0}_{1},\frac{1}{\rho^{0}_{1}})\|_{L^{\infty}})\|h^{0}_{1}\|_{
B^{\N}_{2,1}})\leq1.
\end{aligned}
\label{tropimportant3}
\end{equation}
From (\ref{tropimportant1}) let us choose:
\begin{equation}
\eta=\frac{1}{C}g(\|(\rho^{0}_{1},\frac{1}{\rho^{0}_{1}})\|_{L^{\infty}})\|h^{0}_{1}\|_{
B^{\N-2}_{2,1}}
\label{tropimportant4}
\end{equation}
The condition  (\ref{tropimportant3}) implies by using (\ref{tropimportant4}) that:
\begin{equation}
8 g(\|(\rho^{0}_{1},\frac{1}{\rho^{0}_{1}})\|_{L^{\infty}})\|h^{0}_{1}\|_{
B^{\N-2}_{2,1}} \exp (Cg_{1}(\|(\rho^{0}_{1},\frac{1}{\rho^{0}_{1}})\|_{L^{\infty}})\|h^{0}_{1}\|_{
B^{\N}_{2,1}})\leq1.
\label{tropimportant5}
\end{equation}
Let us point out that this condition corresponds exactly to the condition (\ref{crucinitial}) of the theorem \ref{ftheo3}. It concludes the proof of the stability of the functional $\psi$.
Let us prove now some contraction properties of the functional $\psi$.
\subsubsection*{2) Second step, contraction properties:}
Consider two element $(h,u)$ and  $(h^{'},u^{'})$ in $B(0,R)$, according to proposition \ref{fcrucialprop} we have:
\begin{equation}
\begin{aligned}
&\|\psi(h,u)-\psi(h^{'},u^{'})\|_{E^{\N}}\leq C \big(\|F(h,u)-F(h^{'},u^{'})\|_{\widetilde{L}^{1}_{T}(\widetilde{B}^{\N-1,\N}_{2,1})}
+\|G(h,u)-G(h^{'},u^{'})\|_{\widetilde{L}^{1}_{T}(B^{\N-1}_{2,1})}\big)\\
&\times \exp\biggl(C\int_{\R^{+}}\big( \|u_{1}\|_{B^{\N-\frac{1}{2}}_{2,\infty}}^{4}+\|u_{1}\|^{\frac{4}{3}}_{B^{\N+\frac{1}{2}}_{2,\infty}}+\|\n\ln\rho_{1}\|^{\frac{4}{3}}_{B^{\N+\frac{1}{2}}_{2,\infty}}+\|\n\ln\rho_{1}\|^{\frac{4}{3}}_{B^{\N+\frac{1}{2}}_{2,\infty}}\big)(s)ds\biggl).
\end{aligned}
\label{38}
\end{equation}
Let us deal with the term $\|F(h,u)-F(h^{'},u^{'})\|_{\widetilde{L}^{1}_{T}(\widetilde{B}^{\N-1,\N}_{2,1})}$, we have then by proposition \ref{produit} and by denoting $\delta h=h-h^{'}$ and $\delta u=u-u^{'}$:
\begin{equation}
\begin{aligned}
&\|F(h^{'},u^{'})-F(h,u)\|_{\widetilde{L}^{1}(\widetilde{B}^{\N-1,\N}_{2,\infty})}\lesssim
 \|\delta u\|_{\widetilde{L}^{\frac{4}{3}}_{T}(B^{\N+\frac{1}{2}}_{2,\infty})}\|\n h^{'}\|_{\widetilde{L}^{4}_{T}(\widetilde{B}^{\N-\frac{3}{2},\N-\frac{1}{2}}_{2,\infty})}\\
&\hspace{0,5cm}+ \|\delta u\|_{\widetilde{L}^{\infty}_{T}(B^{\N-1}_{2,\infty})}\|\n h^{'}\|_{\widetilde{L}^{1}_{T}(\widetilde{B}^{\N,\N+1}_{2,\infty})}+\|u\|_{\widetilde{L}^{1}_{T}(B^{\N+1}_{2,\infty})}\|\n \delta h\|_{\widetilde{L}^{\infty}_{T}(\widetilde{B}^{\N-2,\N-1}_{2,\infty})}\\
&\hspace{6,5cm}+\|\n \delta h\|_{\widetilde{L}^{\frac{4}{3}}_{T}(\widetilde{B}^{\N-\frac{1}{2},\N+\frac{1}{2}}_{2,\infty})}\|u\|_{\widetilde{L}^{4}_{T}(B^{\N-\frac{1}{2}}_{2,\infty})}\\[1,5mm]
&\leq C \big(2\|(h,u)\|_{E^{\N}}+2\|(h^{'},u^{'})\|_{E^{\N}}\big)\|(\delta h,\delta u)\|_{E{\N}}.
\end{aligned}
\label{bacFn}
\end{equation}
We proceed similarly for the term $\|G(h,u)-G(h^{'},u^{'})\|_{\widetilde{L}^{1}_{T}(B^{\N-1}_{2,1})}$ (let us mention only that the delicate term $K\n\ln\rho_{1}$ disappears by the subtraction). We get finally for a $C>0$ large enough by using (\ref{supimp}):
\begin{equation}
\begin{aligned}
&\|\psi(h,u)-\psi(h^{'},u^{'})\|_{E^{\N}}\leq C  \big(2\|(h,u)\|_{E^{\N}}+2\|(h^{'},u^{'})\|_{E^{\N}}\big)\|(\delta h,\delta u)\|_{E{\N}}\\
&\times \exp\biggl(C\int_{\R^{+}}\big( \|u_{1}\|_{B^{\N-\frac{1}{2}}_{2,\infty}}^{4}+\|u_{1}\|^{\frac{4}{3}}_{B^{\N+\frac{1}{2}}_{2,\infty}}+\|\n\ln\rho_{1}\|^{\frac{4}{3}}_{B^{\N+\frac{1}{2}}_{2,\infty}}+\|\n\ln\rho_{1}\|^{\frac{4}{3}}_{B^{\N+\frac{1}{2}}_{2,\infty}}\big)(s)ds\biggl),\\[2mm]
&\leq 4C R\exp\big(Cg_{1}(\|(\rho^{0}_{1},\frac{1}{\rho^{0}_{1}})\|_{L^{\infty}})\|h^{0}_{1}\|_{
B^{\N}_{2,1}})\big)\|(\delta h,\delta u)\|_{E{\N}},\\
\end{aligned}
\label{39}
\end{equation}
From (\ref{valeurR}) and (\ref{tropimportant4}) we have with $C$large enough:
\begin{equation}
\begin{aligned}
&4C R\exp\big(Cg_{1}(\|(\rho^{0}_{1},\frac{1}{\rho^{0}_{1}})\|_{L^{\infty}})\|h^{0}_{1}\|_{
B^{\N}_{2,1}})\big)=16C\eta \exp\big(Cg_{1}(\|(\rho^{0}_{1},\frac{1}{\rho^{0}_{1}})\|_{L^{\infty}})\|h^{0}_{1}\|_{
B^{\N}_{2,1}})\big)\\
&=16C g(\|(\rho^{0}_{1},\frac{1}{\rho^{0}_{1}})\|_{L^{\infty}})\|h^{0}_{1}\|_{
B^{\N-2}_{2,1}} \exp\big(Cg_{1}(\|(\rho^{0}_{1},\frac{1}{\rho^{0}_{1}})\|_{L^{\infty}})\|h^{0}_{1}\|_{
B^{\N}_{2,1}})\big).
\end{aligned}
\label{hyperimp}
\end{equation}
In particular we ensure the contraction property via (\ref{39}) if:
\begin{equation}
4C R\exp\big(Cg_{1}(\|(\rho^{0}_{1},\frac{1}{\rho^{0}_{1}})\|_{L^{\infty}})\|h^{0}_{1}\|_{
B^{\N}_{2,1}})\big)\leq\frac{1}{2}.
\label{hyperimp1}
\end{equation}
By using (\ref{hyperimp}), the previous condition (\ref{hyperimp1}) is realized if:
\begin{equation}
16C g(\|(\rho^{0}_{1},\frac{1}{\rho^{0}_{1}})\|_{L^{\infty}})\|h^{0}_{1}\|_{
B^{\N-2}_{2,1}} \exp\big(Cg_{1}(\|(\rho^{0}_{1},\frac{1}{\rho^{0}_{1}})\|_{L^{\infty}})\|h^{0}_{1}\|_{
B^{\N}_{2,1}})\big)\leq\frac{1}{2}.
\label{hyperimp2}
\end{equation}
Let us mention that the condition (\ref{hyperimp2}) is exactly the condition (\ref{crucinitial}) of the theorem \ref{ftheo3}. It achieves the proof of the contraction property of $\psi$ and by the fixed point theorem we have proved the existence of global strong solution under the hypothesis (\ref{crucinitial}). It concludes the proof of theorem \ref{ftheo3}.
\section{Proof of corollary \ref{cor3}}
\label{section5}
It suffices to apply the same proof than for the theorem \ref{ftheo3}, in particular to apply a fixed point theorem for the function $\psi$ previously defined. The main difference concerns the way to deals with the remainder term $K\n\ln\rho_{1}$. Following the same arguments than for the estimate (\ref{bestim3}) we have for a regular function $g$ and $C>0$:
\begin{equation}
\begin{aligned}
\|K\n\ln\rho_{1}\|_{\widetilde{L}^{1}(B^{\N-1}_{2,1})}&\lesssim \|\ln \rho_{1}\|_{
\widetilde{L}^{1}(B^{\N}_{2,1})},\\
&\leq CK g(\|(\rho_{1},\frac{1}{\rho_{1}})\|_{L^{\infty}})\|h_{1}\|_{
\widetilde{L}^{1}(B^{\N}_{2,1})},\\
&\leq CK g(\|(\rho_{1}^{0},\frac{1}{\rho_{1}^{0}})\|_{L^{\infty}})\|h^{0}_{1}\|_{B^{\N-2}_{2,1}}.
\end{aligned}
\label{fbestim3}
\end{equation}
By using exactly the same arguments than in the previous proof we need the following smallness hypothesis with $C>0$ and $g$, $g_{1}$ regular function in order to ensure the stability of the functional $\psi$ for a ball $B(0,R)$ with $R$ defined as previously:
\begin{equation}
16K C g(\|(\rho^{0}_{1},\frac{1}{\rho^{0}_{1}})\|_{L^{\infty}})\|h^{0}_{1}\|_{
B^{\N-2}_{2,1}} \exp\big(Cg_{1}(\|(\rho^{0}_{1},\frac{1}{\rho^{0}_{1}})\|_{L^{\infty}})\|h^{0}_{1}\|_{
B^{\N}_{2,1}})\big)\leq\frac{1}{2}.
\label{hyperimp2}
\end{equation}
If we choose $K$ small enough this last condition will be verified, the rest of the proof follows the same lines as in the proof of theorem \ref{ftheo3} which achieves the proof of corollary \ref{cor3}.
%\subsection{Proof of theorem \ref{theo3}}
%The first step consist in preserving the regularity $B^{\N+\e}_{2,\infty}$ on the density (A void un principe de multiplicateur ou de persŽvŽrance!!!). The second step comes to prove a maximum principle for the density and to conclude via uniqueness result.
\section{Appendix}
\label{appendix}
In this appendix, we just give a technical lemma on the computation of the capillarity tensor.
\begin{lem}
When $\kappa(\rho)=\frac{\kappa}{\rho}$ with $\kappa>0$ then:
$${\rm div}K=\kappa\rho(\n\D\ln\rho+\frac{1}{2}\n(|\n\ln\rho|^{2})).$$
and:
$$
{\rm div}K=\kappa{\rm div}(\rho\n\n\ln\rho).
$$
\end{lem}
{\bf Proof:} We recall that:
$${\rm div}K
=\n\big(\rho\kappa(\rho)\D\rho+\frac{1}{2}(\kappa(\rho)+\rho\kappa^{'}(\rho))|\n\rho|^{2}\big)
-{\rm div}\big(\kappa(\rho)\n\rho\otimes\n\rho\big).$$
When $\kappa(\rho)=\frac{\kappa}{\rho}$, we have:
\begin{equation}
{\rm div}K=\kappa\n\D\rho-\kappa{\rm div}(\frac{1}{\rho}\n\rho\otimes\n\rho).
\label{1a1}
\end{equation}
But as:
$$\D\rho=\rho\D\ln\rho+\frac{1}{\rho}|\n\rho|^{2},$$
we have by injecting this expression in (\ref{1a1}):
\begin{equation}
{\rm div}K=\kappa\rho\n\D\ln\rho+\kappa\n\rho\D\ln\rho+\kappa\n(\frac{1}{\rho}|\n\rho|^{2})-\kappa{\rm div}(\frac{1}{\rho}\n\rho\otimes\n\rho).
\label{a2}
\end{equation}
As we have:
$$\kappa{\rm div}(\frac{1}{\rho}\n\rho\otimes\n\rho)=\kappa\D\ln\rho\n\rho+\n(\frac{1}{\rho}|\n\ln\rho|^{2})-\frac{\kappa}{2}\rho\n(|\n\ln\rho|^{2}).$$
It concludes the first part of the lemma.\\
\\
We now want to prove that we can rewrite the capillarity tensor under the form of a viscosity tensor. To see this, we have:
$$
\begin{aligned}
{\rm div}(\rho\n(\n \ln\rho))_{j}&=\sum_{i}\p_{i}(\rho\p_{ij}\ln\rho),\\
&=\sum_{i}[\p_{i}\rho\p_{ij}\ln\rho+\rho\p_{iij}\ln\rho],\\
&=\rho(\D\n \ln\rho)_{j}+\sum_{i}\rho\p_{i}\ln\rho\p_{j}\p_{i}\ln\rho),\\
&=\rho(\D\n \ln\rho)_{j}+\frac{\rho}{2}(\n(|\n\ln\rho|^{2}))_{j},\\
&={\rm div}K.
\end{aligned}
$$
We have then:
$${\rm div}K=\kappa{\rm div}(\rho\n\n\ln\rho)=\kappa{\rm div}(\rho D(\n\ln\rho)).$$
\hfill {$\Box$}
%\subsection{Ladyzenskaya et al method and maximum principle}
%We now are going to deal with the following parabolic system:
%$$
%\begin{cases}
%\begin{aligned}
%&\p_{t}\rho-\mu\D\rho+{\rm div}(\rho v)={\rm div}F+f,\\
%&\rho(0,\cdot)=\rho_{0}.
%\end{aligned}
%\end{cases}
%$$
%Let us assume that:
%\begin{equation}
%\|(|v|^{2}, |F|^{2}, f)\|_{L^{r}_{T}(L^{q}(\R^{N})}\leq M,
%\end{equation}
%where $r$ and $q$ are arbitrary positive numbers satisfying the condition (Lions semblerait etre capable de faire le cadre critique a voir!!! a-t-il besoin de controller la partie $div b$ negative)
%\begin{equation}
%\frac{1}{q}+\frac{N}{2r}=1-\kappa_{1},
%\end{equation}
%with $0<\kappa_{1}=1$.
%\begin{theorem}
%Let
%\end{theorem}
%{\bf Proof:} Easily by multiplying the equation by $(\rho-k)_{+}$ with $\rho_{0}<k$ and integrating over $(0,t)\times\R^{N}$ we have:
%\begin{equation}
%\begin{aligned}
%&
%\end{aligned}
%\end{equation}

\end{document}